\documentclass[11pt,a4paper]{article}

\usepackage{url,amsmath,amssymb,latexsym,pstricks,mathrsfs,comment,amsthm,graphicx,tikz,tikz-cd,enumerate,accents,pgffor,cite,wrapfig,multicol,float,cases,calc,bibspacing,geometry,bbm,arydshln}
\usepackage[colorlinks]{hyperref}
\usepackage[T1]{fontenc}
\usepackage{ wasysym }
\usepackage{stmaryrd}

\usetikzlibrary{calc}

\allowdisplaybreaks

\newcommand{\nc}{\newcommand}
\nc{\rnc}{\renewcommand}

\let\oldproofname=\proofname
\rnc{\proofname}{\rm\bf{\oldproofname}}

\DeclareMathSymbol{\widehatsym}{\mathord}{largesymbols}{"62}
\newcommand\lowerwidehatsym{%
  \text{\smash{\raisebox{-1.3ex}{%
    $\widehatsym$}}}}
\newcommand\fixwidehat[1]{%
  \mathchoice
    {\accentset{\displaystyle\lowerwidehatsym}{#1}}
    {\accentset{\textstyle\lowerwidehatsym}{#1}}
    {\accentset{\scriptstyle\lowerwidehatsym}{#1}}
    {\accentset{\scriptscriptstyle\lowerwidehatsym}{#1}}
}

\rnc{\th}{\theta}
\nc{\lam}{\lambda}
\nc{\al}{\alpha}
\nc{\si}{\sigma}
\nc{\ve}{\varepsilon}

\nc{\T}{\mathcal T}
\nc{\I}{\mathcal I}
\rnc{\P}{\mathcal P}
\nc{\PT}{\P\T}
\rnc{\S}{\mathcal S}
\rnc{\O}{\mathcal O}
\nc{\PO}{\P\O}
\nc{\OI}{\O\I}
\nc{\Q}{\mathcal Q}

\nc{\TXA}{\T(X,A)}
\nc{\TXal}{\T(X,\al)}
\nc{\IXA}{\I(X,A)}

\rnc{\L}{\mathrel\mathscr{L}}
\nc{\R}{\mathrel\mathscr{R}}
\nc{\J}{\mathrel\mathscr{J}}
\nc{\D}{\mathrel\mathscr{D}}
\rnc{\H}{\mathrel\mathscr{H}}
\nc{\K}{\mathrel\mathscr{K}}
\nc{\leqL}{\leq_{\L}}
\nc{\leqR}{\leq_{\R}}
\nc{\leqH}{\leq_{\H}}
\nc{\leqJ}{\leq_{\J}}
\nc{\leqK}{\leq_{\K}}
\nc{\geqK}{\geq_{\K}}
\nc{\leqLa}{\leq_{\La}}
\nc{\leqRa}{\leq_{\Ra}}
\nc{\leqHa}{\leq_{\Ha}}
\nc{\leqJa}{\leq_{\Ja}}
\nc{\leqJP}{\leq_{\J^P}}
\nc{\leqaJa}{\leq_{\aJa}}
\nc{\leqKa}{\leq_{\Ka}}
\nc{\geqKa}{\geq_{\Ka}}
\nc{\La}{\mathrel{\mathscr L^a}}
\nc{\Ra}{\mathrel{\mathscr R^a}}
\nc{\Ha}{\mathrel{\mathscr H^a}}
\nc{\Da}{\mathrel{\mathscr D^a}}
\nc{\Ja}{\mathrel{\mathscr J^a}}
\nc{\Ka}{\mathrel{\mathscr K^a}}
\nc{\Lha}{\mathrel{\fixwidehat{\mathscr L}^a}}
\nc{\Rha}{\mathrel{\fixwidehat{\mathscr R}^a}}
\nc{\Hha}{\mathrel{\fixwidehat{\mathscr H}^a}}
\nc{\Dha}{\mathrel{\fixwidehat{\mathscr D}^a}}
\nc{\Jha}{\mathrel{\fixwidehat{\mathscr J}^a}}
\nc{\Kha}{\mathrel{\fixwidehat{\mathscr K}^a}}
\nc{\aLh}{\mathrel{{}^a\!\!\fixwidehat{\mathscr L}}}
\nc{\aRh}{\mathrel{{}^a\!\fixwidehat{\mathscr R}}}
\nc{\aHh}{\mathrel{{}^a\!\!\fixwidehat{\mathscr H}}}
\nc{\aDh}{\mathrel{{}^a\!\fixwidehat{\mathscr D}}}
\nc{\aJh}{\mathrel{{}^a\!\!\!\fixwidehat{\mathscr J}}}
\nc{\aKh}{\mathrel{{}^a\!\!\fixwidehat{\mathscr K}}}
\nc{\aLa}{\mathrel{{}^a\!\!\mathscr L^a}}
\nc{\aRa}{\mathrel{{}^a\!\mathscr R^a}}
\nc{\aHa}{\mathrel{{}^a\!\!\mathscr H^a}}
\nc{\aDa}{\mathrel{{}^a\!\mathscr D^a}}
\nc{\aJa}{\mathrel{{}^a\!\!\!\mathscr J^a}}
\nc{\aKa}{\mathrel{{}^a\!\!\mathscr K^a}}
\nc{\aL}{\mathrel{{}^a\!\!\mathscr L}}
\nc{\aR}{\mathrel{{}^a\!\mathscr R}}
\nc{\aH}{\mathrel{{}^a\!\!\mathscr H}}
\nc{\aD}{\mathrel{{}^a\!\mathscr D}}
\nc{\aJ}{\mathrel{{}^a\!\!\!\mathscr J}}
\nc{\aK}{\mathrel{{}^a\!\!\mathscr K}}
\nc{\Rh}{\widehat{R}}
\nc{\Lh}{\widehat{L}}
\nc{\Hh}{\widehat{H}}
\nc{\Dh}{\widehat{D}}
\nc{\Jh}{\widehat{J}}
\nc{\Kh}{\widehat{K}}

\nc{\bbE}{\mathbb E}

\nc{\COMMA}{,\qquad}
\nc{\COMMa}{,\quad}
\nc{\COMma}{,\ \ }
\nc{\anD}{\ \ \ \text{and}\ \ \ }
\nc{\ANd}{\quad\text{and}\quad}
\nc{\AND}{\qquad\text{and}\qquad}

\rnc{\iff}{\ \Leftrightarrow\ }
\rnc{\implies}{\ \Rightarrow\ }
\rnc{\emptyset}{\varnothing}
\nc{\la}{\langle}
\nc{\ra}{\rangle}
\nc{\sub}{\subseteq}
\nc{\sm}{\setminus}
\nc{\mt}{\mapsto}
\nc{\lmap}[1]{\mapstochar \xrightarrow {\ #1\ }}
\nc{\pre}{\preceq}

\nc{\set}[2]{\{#1:#2\}}
\nc{\bigset}[2]{\big\{#1:#2\big\}}
\nc{\pres}[2]{\la#1:#2\ra}

\nc{\dom}{\operatorname{dom}}
\nc{\im}{\operatorname{im}}
\nc{\rank}{\operatorname{rank}}
\nc{\relrank}[2]{\rank(#1\hspace{0.5truemm}{:}\hspace{0.5truemm}#2)}
\nc{\idrank}{\operatorname{idrank}}
\nc{\relidrank}[2]{\idrank(#1\hspace{0.5truemm}{:}\hspace{0.5truemm}#2)}
\nc{\Reg}{\operatorname{Reg}}
\nc{\MI}{\operatorname{MI}}
\nc{\LI}{\operatorname{LI}}
\nc{\RI}{\operatorname{RI}}
\nc{\RP}{\operatorname{RP}}
\nc{\Max}{\operatorname{Max}_{\preceq}}
\nc{\sh}{\operatorname{sh}}
\nc{\defect}{\operatorname{def}}
\nc{\col}{\operatorname{col}}
\nc{\id}{\operatorname{id}}

\nc{\pf}{\begin{proof}}
\nc{\epf}{\end{proof}}
\nc{\epfres}{\hfill\qed}
\nc{\epfreseq}{\tag*{\qed}}
\nc{\bit}{\begin{itemize}}
\nc{\eit}{\end{itemize}}
\nc{\bmc}{\begin{multicols}}
\nc{\emc}{\end{multicols}}
\nc{\itemit}[1]{\item[\emph{(#1)}]}
\nc{\itemnit}[1]{\item[(#1)]}
\nc{\pfitem}[1]{\bigskip\noindent (#1).}
\nc{\pfcase}[1]{\bigskip\noindent {\bf Case #1.}}

\nc{\DClass}[5]{
\begin{tikzpicture}[scale=#5]
\foreach \x/\y in {#3} {\fill[lightgray](\y-1,#2-\x)--(\y,#2-\x)--(\y,#2-\x+1)--(\y-1,#2-\x+1)--(\y-1,#2-\x); \draw(\y-.5,#2-\x+.5)node{$\scriptscriptstyle{#4}$};}
\foreach \x in {0,1,...,#1} {\draw(\x,0)--(\x,#2);}
\foreach \x in {0,1,...,#2} {\draw(0,\x)--(#1,\x);}
\draw[line width=1.3mm] (0,0)--(#1,0)--(#1,#2)--(0,#2)--(0,0)--(#1,0);
\end{tikzpicture}
}

\nc{\DaClass}[7]{
\begin{tikzpicture}[scale=#5]
\foreach \x/\y in {#3} {\fill[lightgray](\y-1,#2-\x)--(\y,#2-\x)--(\y,#2-\x+1)--(\y-1,#2-\x+1)--(\y-1,#2-\x); \draw(\y-.5,#2-\x+.5)node{$\scriptscriptstyle{#4}$};}
\foreach \x in {0,1,...,#1} {\draw(\x,0)--(\x,#2);}
\foreach \x in {0,1,...,#2} {\draw(0,\x)--(#1,\x);}
\foreach \x in {#6} {\draw[line width=1.3mm](\x,0)--(\x,#2);}
\foreach \x in {#7} {\draw[line width=1.3mm](0,\x)--(#1,\x);}
\draw[line width=1.3mm] (0,0)--(#1,0)--(#1,#2)--(0,#2)--(0,0)--(#1,0);
\end{tikzpicture}
}

\numberwithin{equation}{section}

\newtheorem{thm}[equation]{Theorem}
\newtheorem{lemma}[equation]{Lemma}
\newtheorem{cor}[equation]{Corollary}
\newtheorem{prop}[equation]{Proposition}

\theoremstyle{definition}

\newtheorem{rem}[equation]{Remark}

\begin{document}

\title{Structure of principal one-sided ideals\vspace{-0.3cm}}
\author{James East\footnote{The author thanks Prof Jintana Sanwong and her team for some helpful conversations during his visit to Chiang Mai University in 2017.  The financial assistence of CMU is gratefully acknowledged.}\\
{\it\small Centre for Research in Mathematics and Data Science,}\\
{\it\small Western Sydney University, Locked Bag 1797, Penrith NSW 2751, Australia.}\\
{\tt\small J.East\,@\,WesternSydney.edu.au}}
\date{}

\maketitle

\vspace{-1.0cm}

\begin{abstract}
We give a thorough structural analysis of the principal one-sided ideals of arbitrary semigroups, and then apply this to full transformation semigroups and symmetric inverse monoids.  One-sided ideals of these semigroups naturally occur as semigroups of transformations with restricted range or kernel.  


\end{abstract}

\tableofcontents

\section{Introduction}\label{sect:intro}

One- and two-sided ideals play an important role in the structure theory of semigroups.  Principal ideals in particular are directly involved in the definition of Green's relations \cite{Green1951}, and also feature in results on sandwich semigroups and variants \cite{Hickey1983,Sandwiches1,Sandwiches2,DE2015}.  Moreover, many interesting semigroups are one-sided ideals in other naturally occurring semigroups.  For example, the semigroup $T_1$ of all non-negative mappings of the real numbers is a principal left ideal in the monoid $S$ of all real functions, while the semigroup $T_2$ of even functions (which satisfy the identity $(-x)f=xf$) is a principal right ideal of $S$.  Indeed, if $a$ denotes the function $\mathbb R\to\mathbb R:x\mt x^2$, then $T_1=Sa$ and $T_2=aS$.  The semigroups $T_1$ and~$T_2$ are special cases of semigroups of transformations with restricted range or kernel.  Such semigroups have been studied extensively by many authors, particularly from the Thai school of semigroup theory; see for example \cite{Symons1975,SS2017,ZL2017,SS2016,FHQS2016, SS2015,FHQS2014,Sun2013,SS2012, NYK2005,MK2010a,MK2010b,NK2007, Sanwong2011,SS2013,SS2014,FS2014,SS2008,Sullivan2008,MGS2011,MGS2010,SSS2009}.

The main motivation of the current article is to provide a general framework within which to study semigroups such as those above.  Many of the results from the articles just cited follow from general results proved below.  The basic philosophy is to ask:  
\bit
\item[] \emph{Given a semigroup $S$, and an element $a$ of $S$, how does the structure of the principal one-sided ideals~$Sa$ and $aS$ relate to that of $S$?}
\eit
Such questions have been considered extensively for two-sided ideals, and have led to some very interesting studies.  For example, the two-sided ideals of full transformation semigroups consist of all transformations of bounded rank.  Similar characterisations hold for other semigroups of (linear) transformations, endomorphisms and diagrams; for some studies, see for example \cite{EG2017,DEG2017,HM1990,Gray2014,Gray2007,Gray2008,DE2015,DE2018a,DE2018b,Klimov1977}.  In some ways, the structure of a two-sided ideal $I$ of a semigroup $S$ is quite closely tied to that of $S$ itself; for example, if $S$ is regular, then so too is $I$, and every one of Green's relations on $I$ is simply the restriction of the corresponding relation on~$S$ \cite{ERideals}.  In general, neither of these statements hold for one-sided ideals.  As a visual demonstration of this fact, let~$S$ be the full transformation semigroup of degree $5$.  The egg-box diagram of $S$ is pictured in Figure~\ref{fig:TX} (right), and some principal left and right ideals of $S$ are pictured in Figures \ref{fig:TXA} and \ref{fig:TXal}, respectively.  Although these are clearly more complicated than for $S$ itself, certain patterns do seem to emerge.  
Some of the general results we prove can be thought of as formal explanations of such patterns.
Let us now summarise the structure and main results of the paper.  

Section \ref{sect:prelim} contains the preliminary definitions and background results we will need, including new results in Section \ref{subsect:MI} on one-sided identity elements, and properties we call RI- and LI-domination.  

Section \ref{sect:PLI} then gives a thorough treatment of an arbitrary principal left ideal~$Sa$ in a semigroup $S$.  The regular elements of $Sa$ are characterised in Section \ref{subsect:Reg_Sa}, and Green's relations in Section~\ref{subsect:Green_Sa}.  A crucial role in these sections is played by certain sets denoted $P$, $P'$, $P''$ and $P'''$; for example,~$P$ and $P'$ consist of all elements $x\in Sa$ for which $x$ is $\L$- or $\J$-related (in $S$) to $ax$, respectively; when $a$ is regular in~$S$, we have $P''=P'''=Sa$.  In a sense, the main results of Sections \ref{subsect:Reg_Sa} and \ref{subsect:Green_Sa} show that many structural questions concerning $Sa$ may be reduced to the determination of these sets, a somewhat ``lower-level'' task; see especially Theorems \ref{thm:RegSa} and \ref{thm:GreenSa}, and Corollary \ref{cor:GreenSa}.  Sections \ref{subsect:P} and \ref{subsect:rank} identify a natural assumption on the element~$a$ (called \emph{sandwich-regularity} in \cite{Sandwiches1}) under which a more detailed structural analysis may be carried out.  In this case, the set $\Reg(Sa)$ of all regular elements of $Sa$ is a subsemigroup of $Sa$, indeed a right ideal, and in fact $\Reg(Sa)$ is then precisely the set $P$ mentioned above.  When~$a$ is a sandwich-regular idempotent, the structure of $P=\Reg(Sa)$ is closely related not only to that of $S$ itself, but also to the regular monoid~$aSa$.  There is a natural surmorphism (surjective homomorphism) $P\to aSa$, which allows us to describe the idempotents and idempotent-generated subsemigroup of $Sa$ in terms of those of $aSa$ (Theorem \ref{thm:E_Sa}), and describe the Green's structure of $P$ as a kind of ``inflation'' of that of $aSa$ (Theorem~\ref{thm:D_structure_P}; cf.~Remark~\ref{rem:inflation_Sa} and Figure~\ref{fig:inflation_Sa}).  The main results of Section~\ref{subsect:rank} give lower bounds for the (idempotent) ranks of the regular and idempotent-generated subsemigroups of $Sa$, and show that these are exact values when $P$ is RI-dominated; see especially Theorems~\ref{thm:rank_P} and~\ref{thm:rank_EP}.  Finally, Section \ref{subsect:inverse} shows how the whole theory simplifies under an assumption stronger than sandwich-regularity, under which the regular monoid $P=\Reg(Sa)$ is in fact inverse, and even equal to $aSa$ itself (Theorem \ref{thm:inverse_P}).  

Section \ref{sect:PRI} gives the corresponding results for principal right ideals $aS$.  These are direct duals of those of Section \ref{sect:PLI}, so only the main results are stated, and no proofs are given.

Section \ref{sect:TX} then applies the results of Sections \ref{sect:PLI} and \ref{sect:PRI} to the principal one-sided ideals of the full transformation semigroup $\T_X$, which is the semigroup of all self-maps of the set $X$.  The flavour of the results sometimes depend on whether the set $X$ is finite of infinite.  If $a\in\T_X$ is a fixed transformation, and if we write~$A$ and $\al$ for the image and kernel of $a$, then the principal one-sided ideals $\T_Xa$ and $a\T_X$ are precisely the well-studied semigroups 
\[
\TXA = \set{f\in\T_X}{\im(f)\sub A} \AND \TXal = \set{f\in\T_X}{\ker(f)\supseteq\al}
\]
discussed above; see Proposition \ref{prop:TXa_aTX}.  In Section \ref{subsect:Green_TX}, structural information concerning Green's relations and regular elements of $\TXA$ and $\TXal$ is deduced from the general theory, recovering some old results and proving new ones; see Theorems \ref{thm:Green_TXA} and \ref{thm:Green_TXal}.  Section \ref{subsect:Reg_TXA_TXal} thoroughly analyses the regular subsemigroups ${P=\Reg(\TXA)}$ and ${Q=\Reg(\TXal)}$, describing Green's relations and the ideal structure (Theorems~\ref{thm:Green_RegTXA} and~\ref{thm:Green_RegTXal}), calculating the sizes of $P$ and $Q$ (Propositions~\ref{prop:size_P_T} and~\ref{prop:size_Q_T}, and Corollaries~\ref{cor:size_P_T} and~\ref{cor:size_Q_T}), as well as their ranks (Theorems \ref{thm:rank_P_T} and \ref{thm:rank_Q_T}).  Section \ref{subsect:IG_TXA_TXal} concerns the idempotent-generated subsemigroups $\bbE(\TXA)$ and $\bbE(\TXal)$, characterising the elements of these subsemigroups (Theorems~\ref{thm:IGTA} and~\ref{thm:IGTal}), enumerating the idempotents (Proposition \ref{prop:E_TXA_TXal}) and calculating ranks and idempotent ranks (Theorem \ref{thm:E_TXA_TXal}).  Finally, egg-box diagrams are given in Section \ref{subsect:eggbox} (Figures~\ref{fig:TX}--\ref{fig:RXal}) to illustrate many of the results proved in Sections \ref{subsect:Green_TX}--\ref{subsect:IG_TXA_TXal} in special cases.

Section \ref{sect:I} briefly discusses the situation for the principal one-sided ideals of the symmetric inverse monoid~$\I_X$.  Here the strong results of Section \ref{subsect:inverse} apply, and lead to quick proofs of old and new results concerning the semigroups
\[
\set{f\in\T_X}{\im(f)\sub A} \AND \set{f\in\I_X}{\dom(f)\sub A}.
\]
The methods employed in this paper could be applied to a great many other semigroups of mappings, such as partial transformations, linear transformations of vector spaces, or more generally endomorphisms of independence algebras.  It would also be interesting to investigate principal one-sided ideals of diagram monoids such as the partition, Brauer and Temperley-Lieb monoids.

\section{Preliminaries}\label{sect:prelim}

In this section, we fix notation and give some background on semigroups; for more, see \cite{CP1,Hig,Howie,RSbook}.  For a subset $U$ of a semigroup $S$, we write $\la U\ra$ for the subsemigroup of $S$ generated by $U$, which is the smallest subsemigroup of $S$ containing $U$, and consists of all products $u_1\cdots u_k$ for $k\geq1$ and~$u_1,\ldots,u_k\in U$.

\subsection{Green's relations and pre-orders}\label{subsect:Green}

Let $S$ be a semigroup.  As usual, $S^1$ denotes $S$ if $S$ is a monoid; otherwise, $S^1$ denotes $S\cup\{1\}$, where $1$ is an adjoined identity element.  Green's pre-orders $\leqL$, $\leqR$, $\leqJ$ and $\leqH$ are defined, for $x,y\in S$, by
\[
x\leqL y \iff x\in S^1y \COMMA
x\leqR y \iff x\in yS^1 \COMMA
x\leqJ y \iff x\in S^1yS^1 \COMMA
{\leqH} = {\leqL} \cap {\leqR}.
\]
If $\K$ denotes any of $\L$, $\R$, $\J$ or $\H$, then Green's $\K$ relation is defined to be the equivalence ${\leqK}\cap{\geqK}$.  Green's $\D$ relation is defined to be the join (in the lattice of equivalence relations on $S$) of $\L$ and $\R$: i.e., ${\D}={\L}\vee{\R}$ is the smallest equivalence relation containing both $\L$ and $\R$.  It is well known that ${\D}={\J}$ if $S$ is finite, and that ${\D}={\L}\circ{\R}={\R}\circ{\L}$ in any semigroup.  
Note that for any $x,y,z\in S$, $x\leqL y\implies xz\leqL yz$ and so also $x\L y\implies xz\L yz$; the latter says that $\L$ is a \emph{right congruence} (i.e., an equivalence that is invariant under right multiplication).  Dual statements hold for $\leqR$ and~$\R$.

If $x\in S$, and if $\K$ is any of $\L$, $\R$, $\J$, $\H$ or $\D$, we will write $K_x = \set{y\in S}{y\K x}$ for the $\K$-class of $x$ in $S$.  Since ${\D}={\L}\circ{\R}={\R}\circ{\L}$, as noted above, we have $D_x=\bigcup_{y\in L_x}R_y=\bigcup_{y\in R_x}L_y$ for any $x\in S$.  If $\K$ is any of Green's relations other than $\D$, then the set $S/{\K}=\set{K_x}{x\in S}$ of all $\K$-classes of $S$ has a natural partial order induced from the pre-order $\leqK$ on $S$, and we denote this partial order also by $\leqK$: for $x,y\in S$, $K_x\leqK K_y \iff x\leqK y$.  The ordering $\leqJ$ on $\J$-classes is often denoted simply by $\leq$.

If $T$ is a subsemigroup of $S$, then Green's relations on $T$ are not necessarily just the restrictions to $T$ of the corresponding relations on $S$; thus, we will sometimes write $\K^S$ and $\K^T$ for Green's $\K$ relation on $S$ and $T$, respectively, with similar conventions for $\K^S$- and $\K^T$-classes, $K_x^S$ and $K_x^T$.

We may picture elements of a $\D$-class of a semigroup in a so-called \emph{egg-box diagram}: $\R$-related elements are in the same row, $\L$-related elements in the same column, and $\H$-related elements in the same cell.  Group $\H$-classes are usually shaded gray.  When $S$ is finite, we may draw \emph{all} the ${\D}={\J}$-classes in this way, and indicate the $\leq$ ordering on these classes as a Hasse diagram.  For some examples, see Figures \ref{fig:TX}--\ref{fig:RXal}.

\subsection{Idempotents and regularity}\label{subsect:EReg}

An element $x$ of a semigroup $S$ is an \emph{idempotent} if $x=x^2$.  We write
\[
E(S) = \set{x\in S}{x=x^2}
\]
for the set of all idempotents of $S$, and $\bbE(S)=\la E(S)\ra$ for the subsemigroup of $S$ generated by its idempotents.  Any finite semigroup contains an idempotent \cite[Theorem 1.2.2]{Howie}, but this is not necessarily the case for infinite semigroups.

An element $x$ of a semigroup $S$ is \emph{regular} if $x=xyx$ for some $y\in S$; clearly idempotents are regular.  For $x\in S$, we denote by $V(x)=\set{y\in S}{x=xyx,\ y=yxy}$ the set of \emph{inverses} of $x$.  Note that if $y\in S$ is such that $x=xyx$, then $z=yxy$ belongs to $V(x)$, and then $x\R xz$ and $x\L zx$, with $xz,zx\in E(S)$.  We write
\[
\Reg(S) = \set{x\in S}{x=xyx\ (\exists y\in S)}
\]
for the set of all regular elements of the semigroup $S$; note that $\Reg(S)$ may be empty, but not for finite~$S$ (since any finite semigroup contains an idempotent, as noted above).  Any $\D$-class $D$ of a semigroup $S$ satisfies either $D\sub\Reg(S)$ or $D\cap\Reg(S)=\emptyset$: i.e., every element of $D$ is regular, or else no element of $D$ is regular \cite[Theorem 2.11]{CP1}.  Thus, if a $\D$-class $D$ contains an idempotent, then $D$ is a regular $\D$-class.  

A semigroup $S$ is \emph{inverse} \cite{Lawson1998,Petrich1984} if $|V(x)|=1$ for all $x\in S$.  Equivalently, $S$ is inverse if $S$ is regular and its idempotents commute.  Yet another equivalent condition is that every $\R$-class and every $\L$-class contains a unique idempotent.

\subsection{Rank and idempotent rank}\label{subsect:rk}

The \emph{rank} of a semigroup $S$ is the cardinal
\begin{align*}
\rank(S) &= \min\bigset{|U|}{U\sub S,\ S=\la U\ra}.
\intertext{The \emph{relative rank} of $S$ with respect to a subset $A\sub S$ is the cardinal}
\relrank SA &= \min\bigset{|U|}{U\sub S,\ S=\la A\cup U\ra}.
\intertext{If $S$ is an idempotent-generated semigroup, then we may speak of the \emph{idempotent rank} of $S$,}
\idrank(S) &= \min\bigset{|U|}{U\sub E(S),\ S=\la U\ra},
\intertext{and the \emph{relative idempotent rank} of $S$ with respect to a subset $A\sub S$,}
\relidrank SA &= \min\bigset{|U|}{U\sub E(S),\ S=\la A\cup U\ra}.
\end{align*}
We will need the following simple lemma concerning ideals; it is probably well known, but we give a simple proof for completeness.  Recall that a subset $I$ of a semigroup $S$ is an \emph{ideal} if $xy,yx\in I$ for all $x\in I$ and~$y\in S$.  

\begin{lemma}\label{lem:rankWT}
Let $T$ be a subsemigroup of a semigroup $S$ for which $S\sm T$ is an ideal of $S$.  Then
\begin{align*}
\rank(S)&=\relrank ST+\rank(T).
\intertext{If in addition $S$ and $T$ are idempotent-generated, then}
\idrank(S)&=\relidrank ST+\idrank(T).  
\end{align*}
\end{lemma}

\pf
We just prove the second part, as the proof of the first is similar.  Suppose first that $S=\la X\ra$, where $X\sub E(S)$ and $|X|=\idrank(S)$.  Put $Y=X\cap T$ and $Z=X\sm T$.  Because $S\sm T$ is an ideal of $S$, any factorisation over $X$ of an element of $T$ can only involve factors from $Y$, so it follows that $T=\la Y\ra$, and so $|Y|\geq\idrank(T)$.  Since also $S=\la X\ra=\la Y\cup Z\ra=\la\la Y\ra\cup X\ra=\la T\cup Z\ra$, we have $|Z|\geq\relidrank ST$.  But then $\idrank(S) = |X| = |Y|+|Z| \geq \idrank(T) + \relidrank ST$.

The converse may be quickly proved: if $U\sub E(T)$ and $V\sub E(S)$ are such that $T=\la U\ra$, ${S=\la T\cup V\ra}$, $|U|=\idrank(T)$ and $|V|=\relidrank ST$, then $S=\la T\cup V\ra = \la\la U\ra\cup V\ra=\la U\cup V\ra$, and it follows that ${\idrank(S)\leq |U\cup V|\leq|U|+|V|=\idrank(T) + \relidrank ST}$.  
\epf

\subsection{Left and right groups}\label{subsect:LG}

Recall that a \emph{left zero band} is a semigroup $U$ with product $uv=u$ for all $u,v\in U$.  Recall that a \emph{left group} is a semigroup $S$ isomorphic to a direct product $U\times G$, where $U$ is a left zero band and $G$ a group; in this case, we say that $S$ is a \emph{left group of degree $|U|$ over $G$}.  It is easy to show that a semigroup is a left group if and only if it is a union of groups and its idempotents form a left zero band.  \emph{Right zero bands} and \emph{right groups} are defined analogously.  More information on left and right groups can be found in \cite[Section 1.11]{CP1}.

Here we prove two basic results concerning left groups; there are obvious dual statements for right groups, but we will omit these.  The first follows from much stronger results of Ru\v skuc \cite{Ruskuc1994} (see also \cite[Proposition~4.11]{Sandwiches1}), but we include a simple direct proof for convenience.

\begin{lemma}\label{lem:rank_left_group}
If $S$ is a left group of degree $\rho$ over $G$, then $\rank(S)=\max(\rho,\rank(G))$.  
\end{lemma}

\pf
Without loss of generality, we may assume that $S=U\times G$, where $U$ is a left zero band of size $\rho$.  Since $uv=u$ for all $u,v\in U$, clearly $\rank(U)=|U|=\rho$.  Since $U$ and $G$ are both homomorphic images of~$S$, we have $\rank(S)\geq\rank(U)=\rho$ and $\rank(S)\geq\rank(G)$, so that $\rank(S)\geq\max(\rho,\rank(G))$.  

For the converse, write $U=\set{u_i}{i\in I}$ where $|I|=\rho$, and let $X=\set{x_j}{j\in J}$ be a generating set for~$G$ with $|J|=\rank(G)$.  For notational convenience, we will assume that $|I|\leq|J|$; the other case is treated in almost identical fashion.  Without loss of generality, we may assume that $I\sub J$.  For each $j\in J\sm I$, let $u_j$ be an arbitrary element of $U$.  So also $U=\set{u_j}{j\in J}$.  Now put $Z=\set{(u_j,x_j)}{j\in J}$.  Since $|Z|=|J|=\rho=\max(\rho,\rank(G))$, the proof will be complete if we can show that $S=\la Z\ra$.  To do so, let $u\in U$ and $g\in G$ be arbitrary.  Now, $u=u_j$ for some $j\in J$.  Since $G=\la X\ra$, we have $x_j^{-1}g = x_{j_1}\cdots x_{j_k}$ for some $j_1,\ldots,j_k\in J$.  But then $(u,g)=(u_j,x_j)(u_{j_1},x_{j_1})\cdots(u_{j_k},x_{j_k})\in\la Z\ra$, as required.
\epf

The next result is a little more general than we need, but there is no greater difficulty in proving the stronger statement.

\begin{lemma}\label{lem:LG_subs}
Let $U$ be a left zero band and $M$ a monoid with identity $e$.  Suppose $T$ is a subsemigroup of $U\times M$ such that $T$ contains $U\times\{e\}$.  Then $T=U\times W$ for some submonoid $W$ of $M$.
\end{lemma}

\pf
Put $W=\set{x\in M}{(u,x)\in T\ (\exists u\in U)}$.  Clearly $W$ is a submonoid of $M$, and clearly $T\sub U\times W$.  Conversely, let $(u,w)\in U\times W$ be arbitrary.  By definition of $W$, there exists $v\in U$ such that $(v,w)\in T$.  By assumption, $(u,e)\in T$.  But then $(u,w)=(u,e)(v,w)\in T$, showing that $U\times W\sub T$.
\epf

\subsection{One-sided identities and mid-identities}\label{subsect:MI}

In our investigations of principal one-sided ideals, a crucial role will be played by one-sided identities and mid-identities.  Here we review the definitions, and prove some results that will highlight the importance of these kinds of elements.

Recall that a \emph{right identity} of a semigroup $S$ is an element $u\in S$ such that $x=xu$ for all $x\in S$.  \emph{Left identities} are defined analogously.  We write $\RI(S)$ and $\LI(S)$ for the sets of all right and left identities of~$S$, respectively.  Note that either or both of these sets might be empty, but if they are both non-empty, then~$S$ is a monoid and $\RI(S)=\LI(S)$ consists entirely of the (unique, two-sided) identity element of $S$.  

Recall \cite{Yamada1955} that a \emph{mid-identity} of a semigroup $S$ is an element $u\in S$ such that $xy=xuy$ for all $x,y\in S$.  We write $\MI(S)$ for the set of all mid-identities of $S$.  Again, $\MI(S)$ may be empty, but we note that $\MI(S)$ always contains both $\RI(S)$ and~$\LI(S)$.  
The next lemma contains some further basic results.

\begin{lemma}\label{lem:MI}
Let $S$ be a semigroup.
\bit
\itemit{i} If $u\in\MI(S)$ and if $u=uv$ or $u=vu$ for some $v\in S$, then $u\in E(S)$.
\itemit{ii} If $S$ is regular or if $S$ has a left or right identity, then $\MI(S)\sub E(S)$.
\itemit{iii} If $\RI(S)\not=\emptyset$, then $\MI(S)=\RI(S)$.
\itemit{iv} If $\LI(S)\not=\emptyset$, then $\MI(S)=\LI(S)$.
\eit
\end{lemma}

\pf
(i).  If $u\in\MI(S)$ and $u=uv$ for some $v\in S$, then $u=uv=uuv=uu$.  The $u=vu$ case is similar.

\pfitem{ii}  This follows from (i), since if $S$ is regular or if $S$ has a left or right identity, then any mid-identity $u$ of~$S$ satisfies $u=uv$ or $u=vu$ for some $v\in S$.

\pfitem{iii) and (iv}  We just prove (iii), as (iv) is dual.  Suppose $\RI(S)\not=\emptyset$, and let $e\in\RI(S)$.  We have already noted that $\RI(S)\sub\MI(S)$.  For the converse, suppose $u\in\MI(S)$, and let $x\in S$ be arbitrary.  Then since $e$ is a right identity and $u$ a mid-identity, we have $x=xe=xue=xu$, so that $u\in\RI(S)$.
\epf

\begin{rem}\label{rem:aa=aaa}
We need not have $\MI(S)\sub E(S)$ in general.  For example, consider the semigroup $S$ given by the presentation $\pres{a}{a^3=a^2}$, so that $S=\{a,a^2\}$ with $a\not=a^2$.  Then $\MI(S)=S$, while $E(S)=\{a^2\}$.
\end{rem}

Recall \cite{Mitsch1986} that there is a natural partial order $\pre$ on a regular semigroup $S$ defined, for $x,y\in S$, by $x\pre y$ if and only if $x=ey=yf$ for some idempotents $e,f\in E(S)$.  If $e,f\in E(S)$, then it is easy to show that $e\pre f$ if and only if $e=fef$ (which is itself equivalent to $e=ef=fe$).  

Recall \cite{Sandwiches1} that a regular semigroup $S$ is \emph{MI-dominated} if each idempotent of $S$ is $\pre$-below a mid-identity.
The concept of MI-domination was used in \cite{Sandwiches1} to describe the structure of sandwich semigroups, and it will be used in the current article (in an equivalent form to be described shortly) to describe the structure of principal one-sided ideals.  

If $S$ is a semigroup and $e\in E(S)$ an idempotent of $S$, then $eSe$ is a subsemigroup of $S$ called the \emph{local monoid} of $S$ with respect to $e$; as the name suggests, $eSe$ is a monoid with identity $e$.
MI-domination is especially useful because of the next result, which is \cite[Proposition~4.3]{Sandwiches1}, and which shows (among other things) that MI-dominated semigroups are unions of local monoids corresponding to mid-identities, all of which are naturally isomorphic.

\begin{prop}\label{prop:MI}
Let $S$ be a regular semigroup, write $M=\MI(S)$, and suppose $M\not=\emptyset$.
\bit
\itemit{i} If $e\in M$, then the map $S\to eSe:x\mt exe$ is a surmorphism.
\itemit{ii} If $e,f\in M$, then the maps $eSe\to fSf:x\mt fxf$ and $fSf\to eSe:x\mt exe$ are mutually inverse isomorphisms.
\itemit{iii} The set $\bigcup_{e\in M}eSe = MSM$ is a subsemigroup of $S$.
\itemit{iv} $S$ is MI-dominated if and only if $S=\bigcup_{e\in M}eSe$. \epfres
\end{itemize}
\end{prop}

It turns out that the MI-domination property has an equivalent reformulation in terms of one-sided identity elements if the semigroup has any of these.
We say that a semigroup $S$ is \emph{RI-dominated} if every element of $S$ is $\leqR$-below a right identity of $S$.  (Note that any element of any semigroup is trivially $\leqL$-below any right identity the semigroup may contain.)  \emph{LI-dominated} semigroups are defined analogously.

\begin{lemma}\label{lem:RILI}
Let $S$ be a regular semigroup.
\bit
\itemit{i} If $\RI(S)\not=\emptyset$, then $S$ is MI-dominated if and only if it is RI-dominated.
\itemit{ii} If $\LI(S)\not=\emptyset$, then $S$ is MI-dominated if and only if it is LI-dominated.
\eit
\end{lemma}

\pf
We just prove (i), as (ii) is dual.  Suppose $\RI(S)\not=\emptyset$.  By Lemma \ref{lem:MI}(iii), we have $\MI(S)=\RI(S)$.

\pfitem{$\Rightarrow$}  Suppose first that $S$ is MI-dominated.  Let $x\in S$ be arbitrary; we must show that $x$ is $\leqR$-below some right identity.  Since $S$ is regular, $x=ex$ for some $e\in E(S)$.  Since $S$ is MI-dominated, $e\pre u$ for some $u\in\MI(S)=\RI(S)$, and so $e=ueu$.  But then $x=ex=ueux\leqR u$.

\pfitem{$\Leftarrow$}  Conversely, suppose $S$ is RI-dominated.  Let $e\in E(S)$ be arbitrary; we must show that $e$ is $\pre$-below some mid-identity.  Since $S$ is RI-dominated, $e\leqR u$ for some $u\in\RI(S)=\MI(S)$.  Since $u$ is a right identity, $e=eu$, while $e\leqR u$ gives $e=ux$ for some $x\in S^1$.  But then $e=ux=uux=ue=ueu$, so that $e\pre u$.  
\epf

\subsection{Transformation semigroups}\label{subsect:trans}

Let $X$ be an arbitrary set.  A \emph{partial transformation} of $X$ is a function from a subset of $X$ into $X$.  The set of all such partital transformations is denoted by $\PT_X$, and is a semigroup under composition, known as the \emph{partial transformation semigroup over $X$}.  For $f\in\PT_X$, we write $\dom(f)$ and $\im(f)$ for the domain and image (or range) of $f$, which are defined in the standard way; we also write
\[
\ker(f) = \set{(x,y)\in\dom(f)\times\dom(f)}{xf=yf} \AND \rank(f)=|{\im(f)}|
\]
for the \emph{kernel} and \emph{rank} of $f$.  Note that $\dom(f)$ and $\im(f)$ are subsets of $X$, $\ker(f)$ is an equivalence on $f$, and $\rank(f)$ is a cardinal between $0$ and $|X|$.  As usual, if $\si$ is an equivalence on a set $Y$, we write $Y/\si$ for the set of $\si$-classes of $Y$; for brevity, we will write $\Vert\si\Vert=|Y/\si|$ for the number of such $\si$-classes.  Note that for $f\in\PT_X$, we also have $\rank(f)=\Vert{\ker(f)}\Vert$.

The \emph{full transformation semigroup} and \emph{symmetric inverse monoid} over $X$ are, respectively, the subsemigroups $\T_X$ and $\I_X$ of $\PT_X$ defined by
\[
\T_X = \set{f\in\PT_X}{\dom(f)=X} \AND \I_X = \set{f\in\PT_X}{f\text{ is injective}}.
\]
Green's relations and pre-orders may easily be described on these monoids in terms of the parameters defined above.  The next result is easily established; see for example \cite[Section 3.1]{Sandwiches2}.  If $\si$ is an equivalence relation on a set $Y$, and if $Z\sub Y$, we write $\si|_Z=\si\cap(Z\times Z)$ for the restriction of $\si$ to $Z$.

\begin{thm}\label{thm:T}
Let $\Q_X$ be any of the semigroups $\PT_X$, $\T_X$ or $\I_X$.  Then $\Q_X$ is a regular monoid.  Further, if $f,g\in\Q_X$, then
\bit
\itemit{i} $f\leqL g \iff \im(f)\sub\im(g)$,
\itemit{ii} $f\leqR g \iff \dom(f)\sub\dom(g)$ and $\ker(f)\supseteq\ker(g)|_{\dom(f)}$,
\itemit{iii} $f\leqJ g \iff \rank(f)\leq\rank(g)$,
\itemit{iv} $f\L g \iff \im(f)=\im(g)$,
\itemit{v} $f\R g \iff \ker(f)=\ker(g)$,
\itemit{vi} $f\J g \iff f\D g \iff \rank(f)=\rank(g)$.  \epfres
\eit
\end{thm}

\begin{rem}
There are simplifications of the $\leqR$ relation in the case of $\T_X$ and $\I_X$ because of the  form of the elements of these monoids.  In $\T_X$, $f\leqR g \iff \ker(f)\supseteq\ker(g)$.  In $\I_X$, $f\leqR g \iff \dom(f)\sub\dom(g)$.
\end{rem}

We also require some combinatorial data concerning Green's classes.  For cardinals $\mu,\nu$ with $\nu\leq\mu$, we write
\bit
\item $\mu!$ for the number of permutations of a set of size $\mu$,
\item $\binom\mu\nu$ for the number of subsets of size $\nu$ of a set of size $\mu$,
\item $S(\mu,\nu)$ for the number of equivalence classes with $\nu$ classes in a set of size $\mu$.
\eit
Note that if $\mu$ is infinite, then $\mu!=2^\mu$, $\binom\mu\nu=\mu^\nu$, $S(\mu,1)=1$, and $S(\mu,\nu)=2^\mu$ if $\nu\geq2$; see \cite{Jech2003}.  If $\mu$ is finite, then $\mu!$, $\binom\mu\nu$ and $S(\mu,\nu)$ have their usual meanings, as factorials, binomial coefficients and Stirling numbers (of the second kind), respectively.  

We write $\S_X$ for the \emph{symmetric group} over $X$, which consists of all permutations of $X$, and is the group of units of $\PT_X$, $\T_X$ and $\I_X$.  If $\mu$ is a cardinal, then we may consider the semigroups $\PT_\mu$, $\S_\mu$, etc., by interpreting $\mu$ as an ordinal (and hence as a set).

If $0\leq\mu\leq|X|$ is an arbitrary cardinal, and if $\Q_X$ is any of $\PT_X$, $\T_X$ or $\I_X$, we write
\[
D_\mu(\Q_X)=\set{f\in\Q_X}{\rank(f)=\mu}.
\]
The next result is easily established; see also \cite[Corollary 2.4]{Sandwiches2}.

\begin{prop}\label{prop:combinatorics}
Let $X$ be a set, let $\Q_X$ be any of $\PT_X$, $\T_X$ or $\I_X$, and let $z=1$ if $\Q_X=\T_X$ or $z=0$ otherwise.  Then the ${\D}={\J}$-classes of $\Q_X$ are the sets 
\[
D_\mu(\Q_X)=\set{f\in\Q_X}{\rank(f)=\mu} \qquad\text{for $z\leq\mu\leq|X|$.}
\]
These form a chain under the $\J$-class ordering: $D_\mu(\Q_X)\leq D_\nu(\Q_X) \iff \mu\leq\nu$.  Further, if $z\leq\mu\leq|X|$ is a cardinal, then
\bit
\itemit{i} $|D_\mu(\Q_X) / {\L}| = \binom{|X|}\mu$, 
\itemit{ii} $|D_\mu(\PT_X) / {\R}|=S(|X|+1,\mu+1)$, \quad $|D_\mu(\T_X) / {\R}|=S(|X|,\mu)$, \quad $|D_\mu(\I_X) / {\R}|=\binom{|X|}\mu$, 
\itemit{iii} $|D_\mu(\PT_X) / {\H}|=\binom{|X|}\mu S(|X|+1,\mu+1)$, \quad $|D_\mu(\T_X) / {\H}|=\binom{|X|}\mu S(|X|,\mu)$, \quad $|D_\mu(\I_X) / {\H}|=\binom{|X|}\mu^2$, 
\itemit{iv} any $\H$-class of $\Q_X$ contained in $D_\mu(\Q_X)$ has size $\mu!$, 
\itemit{v} any group $\H$-class of $\Q_X$ contained in $D_\mu(\Q_X)$ is isomorphic to $\S_\mu$.  \epfres
\eit
\end{prop}

We also need to know about the idempotent-generated subsemigroups $\bbE(\T_X)$ and $\bbE(\PT_X)$ in the finite case.  
The next result is \cite[Theorem~I]{Howie1966}; the case of infinite $X$ is also given in \cite[Theorem III]{Howie1966}.  We write $\id_X$ for the identity mapping on $X$.

\begin{thm}\label{thm:IGT}
If $X$ is a finite set with $|X|\geq2$, then $\bbE(\T_X)=\{\id_X\}\cup(\T_X\sm\S_X)$. Further,
\[
\epfreseq
\rank(\bbE(\T_X))=\idrank(\bbE(\T_X))=
\begin{cases}
3 &\text{if $|X|=2$}\\
\binom{|X|}2+1 &\text{if $|X|\geq3$.}  
\end{cases}
\]
\end{thm}

Finally, we recall some standard notation for partial transformations.  If $f\in\PT_X$, we write $f=\binom{F_i}{f_i}_{i\in I}$ to indicate that
\[
\dom(f) = \bigcup_{i\in I}F_i \COMMA \im(f) = \set{f_i}{i\in I} \COMMA xf=f_i\ (\forall x\in F_i) \COMMA \dom(f)/\ker(f) = \set{F_i}{i\in I}.
\]
Sometimes we will write $f=\binom{F_i}{f_i}$, with the indexing set $I$ being implied, rather than explicitly stated.  If $f=\binom{F_i}{f_i}$ belongs to $\T_X$, then $X=\bigcup_{i\in I}F_i$, while if $f$ belongs to $\I_X$, then $|F_i|=1$ for all $i$.

\section{Principal left ideals}\label{sect:PLI}

A subset $I$ of a semigroup $S$ is a \emph{left ideal} if it is closed under left multiplication by arbitrary elements of~$S$: i.e., for all $x\in S$ and $y\in I$, we have $xy\in I$.  The \emph{principal left ideal} generated by an element $a$ of the semigroup $S$ is the set
\[
Sa = \set{xa}{x\in S}.
\]
\emph{(Principal) right ideals} of $S$ are defined dually.  The purpose of this paper is to develop a structure theory of principal left and right ideals; since these theories are dual, we give a detailed treatment of left ideals in the current section, and then simply state the corresponding results concerning right ideals in Section \ref{sect:PRI}.

Note that some authors would define the principal left ideal generated by $a$ to be $S^1a=Sa\cup\{a\}$.  In many cases we have $S^1a=Sa$, such as when $S$ is a monoid (or just has a left identity element) or when $a$ is regular.  In order to be as general as possible, the results that follow concern $Sa$, but results concerning $S^1a$ may be easily be obtained by simply replacing $S$ by $S^1$, and considering $S^1a$ as a principal left ideal (in our sense) of $S^1$.

This section has five subsections.
Subsection \ref{subsect:Reg_Sa} characterises the regular elements of $Sa$, and gives a sufficient condition for the set $\Reg(Sa)$ to be a subsemigroup (indeed, right ideal) of $Sa$.
Subsection \ref{subsect:Green_Sa} describes Green's relations on $Sa$, characterising these in terms of the corresponding relations on $S$ and certain subsets of $Sa$.
Subsection \ref{subsect:P} investigates the structure of the regular subsemigroup $\Reg(Sa)$ in the case that~$a$ is a so-called \emph{sandwich-regular} idempotent of $S$.  It is shown that the structure of $\Reg(Sa)$ is closely related to that of the (local) monoid $aSa$; crucial use is made of a natural surmorphism $\phi:\Reg(Sa)\to aSa$.  The idempotent-generated subsemigroup of $Sa$ is also related to that of $aSa$.
Subsection \ref{subsect:rank} explores the rank (and idempotent rank, where appropriate) of the regular and idempotent-generated subsemigroups of $Sa$, again relating these to corresponding (idempotent) ranks in $aSa$.  Lower bounds for these (idempotent) ranks are given, and shown to be exact values in the case of $\Reg(Sa)$ being RI-dominated.
Finally, Subsection~\ref{subsect:inverse} identifies a property stronger than sandwich-regularity under which the whole theory simplifies greatly, as we will show that $\Reg(Sa)=aSa$ is an inverse monoid.

\subsection[Regular elements of $Sa$]{\boldmath Regular elements of $Sa$}\label{subsect:Reg_Sa} 

For the duration of this subsection, we fix a semigroup $S$ and an element $a$ of $S$.  Our main goal here is to characterise the set
\[
\Reg(Sa) = \set{x\in Sa}{x=xyx\ (\exists y\in Sa)}
\]
of regular elements of the semigroup $Sa$.  We will see later that under some mild regularity assumptions on~$a$ and $S$ (which hold if $S$ is regular, for example), the set $\Reg(Sa)$ is in fact a subsemigroup of $Sa$.

A crucial role in all that follows is played by the set $P$ defined by
\[
P = \set{x\in Sa}{x\L ax}.
\]
Since $ax\leqL x$ for any $x\in S$, we could equivalently have defined $P$ as $\set{x\in Sa}{x\leqL ax}$.
Note that if $x\in P$, then $x=wax$ for some $w\in S^1$; in fact, we may assume that $w\in S$, since if $w=1$, then $x=ax=aax$.  

\begin{lemma}\label{lem:PRightIdeal}
The set $P$ is a right ideal of $Sa$.
\end{lemma}

\pf
Let $x\in P$ and $y\in Sa$.  Certainly $xy\in Sa$.  But also $x\L ax$ implies $xy\L axy$, since $\L$ is a right congruence, so it follows that $xy\in P$, as required.
\epf

The next result characterises the set $\Reg(Sa)$ of all regular elements of $Sa$.

\begin{thm}\label{thm:RegSa}
Let $S$ be a semigroup, let $a\in S$, and define $P=\set{x\in S}{x\L ax}$.  Then
\[
\Reg(Sa) = \Reg(S) \cap P.
\]
\end{thm}

\pf
First suppose $x\in\Reg(Sa)$.  So $x\in Sa$ and $x=xyx$ for some $y\in Sa$.  Certainly then $x\in\Reg(S)$.  Also, since $y\in Sa$, we have $y=za$ for some $z\in S$, in which case $x=xyx=x(za)x=(xz)ax$, so that $x\L ax$, which gives $x\in P$.

Conversely, suppose $x\in \Reg(S)\cap P$.  Since $x\in\Reg(S)$, we have $x=xyx$ for some $y\in S$.  Since $x\in P$, we have $x\in Sa$ and $x=zax$ for some $z\in S^1$.
But then $x=xyx=xy(zax) = x(yza)x$; since $yza\in Sa$, it follows that $x\in\Reg(Sa)$.
\epf

It follows from Theorem \ref{thm:RegSa} that $\Reg(Sa)=P$ if $S$ is regular, or even if every element of $P$ is regular (in $S$).  The next result identifies a weaker property than regularity of $S$ that ensures $\Reg(Sa)=P$.

\begin{cor}\label{cor:RegSa}
If $aSa\sub\Reg(S)$, then $P\sub\Reg(S)$.  Consequently, $\Reg(Sa)=P$ is a right ideal of $Sa$ in this case.
\end{cor}

\pf
The second assertion follows from the first, because of Theorem \ref{thm:RegSa} and Lemma \ref{lem:PRightIdeal}.  To prove the first assertion, suppose $aSa\sub\Reg(S)$, and let $x\in P$.  So $x\in Sa$ and $x=yax$ for some $y\in S^1$.  Since $ax\in aSa\sub\Reg(S)$, we have $ax=axzax$ for some $z\in S$.  But then $x=yax=y(axzax)=x(za)x$, so that $x\in\Reg(S)$, as required.
\epf

\begin{rem}
Note that the condition $aSa\sub\Reg(S)$ does not imply that $a\in\Reg(S)$ in general.  For example, if $S$ is the semigroup defined by the presentation $\pres{a}{a^3=a^2}$, as in Remark \ref{rem:aa=aaa}, then we have ${aSa=\Reg(S)=\{a^2\}}$.  In \cite{Sandwiches1}, an element $a$ of a semigroup $S$ satisfying $\{a\}\cup aSa\sub\Reg(S)$ was called \emph{sandwich-regular}; this property will play an important role in subsequent discussions.
\end{rem}

\subsection[Green's relations in $Sa$]{\boldmath Green's relations in $Sa$}\label{subsect:Green_Sa}

We now consider Green's relations on the principal left ideal $Sa$.  Theorem \ref{thm:GreenSa} characterises these in terms of Green's relations on~$S$ and certain subsets of $S$, including $P$ defined above.  Corollary \ref{cor:GreenSa} shows how these characterisations simplify in the case that $a$ is a regular element of $S$.

We will continue to write Green's relations on $S$ as $\L$, $\R$, etc., and we will continue to write $K_x$ for the $\K$-class of $x\in S$ for any of Green's relations $\K$.  However, in order to avoid confusion, we will write $\K^a$ for Green's $\K$-relation on $Sa$.  If $x\in Sa$, we will write $K_x^a=\set{y\in Sa}{x\K^ay}$ for the $\K^a$-class of $x$ in~$Sa$.  It is clear that $K_x^a\sub K_x\cap Sa$ for any $x\in Sa$ and for any $\K$.  

Our characterisation of Green's relations on $Sa$ (Theorem \ref{thm:GreenSa}) uses the set $P$ defined above, as well as three more sets:
\[
P' = \set{x\in Sa}{x\J ax} \COMMA P'' = \set{x\in S}{x\in xSa} \COMMA P''' = \set{x\in S}{x\in S^1xSa}.
\]
Note that we could have equivalently defined $P''$ as $\set{x\in S}{x\in xS^1a}$; indeed, if $x=xa$, then $x=xaa\in xSa$.  Similarly, we could have defined $P'''$ as $\set{x\in S}{x\in S^1xS^1a}$.  Also observe that clearly $P''\sub P'''\sub Sa$.  If $a$ is regular, then we may make a much stronger statement:

\begin{lemma}\label{lem:P''P'''}
If $a$ is a regular element of $S$, then $P''=P'''=Sa$.
\end{lemma}

\pf
In light of the above observation, it suffices to show that $Sa\sub P''$.  Let $b\in S$ be such that $a=aba$, and suppose $x\in S$ is arbitrary.  Then $xa=xaba\in (xa)Sa$, so that $xa\in P''$, showing that $Sa\sub P''$.
\epf

\begin{rem}
If $a$ is not regular, then it is possible for $P''$ and $P'''$ to be proper subsets of $Sa$.  For example, let $S$ be defined by the presentation $\pres{a}{a^4=a^3}$.  Then $Sa=\{a^2,a^3\}$, while $P''=P'''=\{a^3\}$.  We also clearly have $P\sub P'$.  Although $P$ and $P'$ are not always equal, they are if $S$ is \emph{left-stable} (which is the case, for example, if $S$ is finite); cf.~\cite{EH2019}.
\end{rem}

The next technical lemma will be used on a number of occasions in the proof of the theorem that follows.

\begin{lemma}\label{lem:PP''}
Let $x\in S$.
\bit
\itemit{i} If $x\in P$, and if $y\in Sa$ satisfies $y\leqR x$, then $y\in P$.  In particular, $x\in P \implies R_x\cap Sa\sub P$.
\itemit{ii} If $x\in P''$, and if $y\in Sa$ satisfies $y\leqL x$, then $y\in P''$.  In particular, $x\in P'' \implies L_x\cap Sa\sub P''$.
\eit
\end{lemma}

\pf
We just prove (i), as the proof of (ii) is almost identical.  It clearly suffices to prove the first assertion, so suppose $x\in P$, and let $y\in Sa$ with $y\leqR x$.  Then we have $x=uax$ and $y=xv$ for some $u,v\in S^1$, and so $y=xv=uaxv=uay$, which gives $y\L ay$ and $y\in P$.
\epf

\begin{thm}\label{thm:GreenSa}
Let $S$ be a semigroup, let $a\in S$, and define the sets
\[
P = \set{x\in Sa}{x\L ax} \COMma P' = \set{x\in Sa}{x\J ax} \COMma P'' = \set{x\in S}{x\in xSa} \COMma P''' = \set{x\in S}{x\in S^1xSa}.
\]
Then for any $x\in Sa$,
\bit\bmc2
\itemit{i} $L_x^a = \begin{cases} L_x\cap P &\hspace{2.4mm}\text{if $x\in P$} \\ \{x\} &\hspace{2.4mm}\text{if $x\not\in P$,} \end{cases}$
\itemit{ii} $R_x^a = \begin{cases} R_x\cap P'' &\text{if $x\in P''$} \\ \{x\} &\text{if $x\not\in P''$,} \end{cases}$
\itemit{iii} $H_x^a = \begin{cases} H_x &\hspace{7.3mm}\text{if $x\in P\cap P''$} \\ \{x\} &\hspace{7.4mm}\text{if $x\not\in P\cap P''$,} \end{cases}$
\itemit{iv} $D_x^a = \begin{cases} D_x\cap P\cap P'' &\text{if $x\in P\cap P''$} \\ R_x^a &\text{if $x\not\in P$} \\ L_x^a &\text{if $x\not\in P''$,} \end{cases}$
\itemit{v} $J_x^a = \begin{cases} J_x\cap P'\cap P''' &\hspace{0.3mm}\text{if $x\in P'\cap P'''$} \\ D_x^a &\hspace{0.3mm}\text{if $x\not\in P'\cap P'''$.} \end{cases}$
\item[] ~
\emc\eit
\end{thm}

\pf
(i).  Suppose $|L_x^a|\geq2$.  Let $y\in L_x^a\sm\{x\}$.  Then $x=uy$ and $y=vx$ for some $u,v\in Sa$.  Since $v\in Sa$, we may write $v=wa$ for some $w\in S$, and so $x=uy=uvx=uwax$, which gives $x\L ax$.  Since also $x\in Sa$, we have $x\in P$.  We have shown that $|L_x^a|\geq2 \implies x\in P$.  The contrapositive of this says that $x\not\in P\implies L_x^a=\{x\}$.  

Now suppose $x\in P$, so that $x=wax$ for some $w\in S$.  To complete the proof of (i), we must show that $L_x^a=L_x\cap P$.  To show the forwards inclusion, suppose $y\in L_x^a$.  Certainly $y\in L_x\cap P$ if $y=x$, so suppose $y\not=x$.  Then certainly $y\in L_x$, and also $|L_y^a|=|L_x^a|\geq2$, so the previous paragraph gives $y\in P$; thus, $y\in L_x\cap P$.
Conversely, suppose $y\in L_x\cap P$, so that $x=uy$, $y=vx$ and $y=zay$ for some $u,v\in S^1$ and $z\in S$.  Then $x=uy=uzay$ and $y=vx=vwax$; since $uza,vwa\in Sa$, it follows that $x\L^ay$, and $y\in L_x^a$ as required.

\pfitem{ii}  Suppose $|R_x^a|\geq2$.  Let $y\in R_x^a\sm\{x\}$.  Then $x=yu$ and $y=xv$ for some $u,v\in Sa$, and so $x=yu=xvu\in xSa$ (since $u\in Sa$), so that $x\in P''$.  We have shown that $|R_x^a|\geq2 \implies x\in P''$.  The contrapositive of this says that $x\not\in P'' \implies R_x^a=\{x\}$.  

Now suppose $x\in P''$, so that $x=xwa$ for some $w\in S$.  To complete the proof of (ii), we must show that $R_x^a=R_x\cap P''$.  To show the forwards inclusion, suppose $y\in R_x^a$.  Certainly $y\in R_x\cap P''$ if $y=x$, so suppose $y\not=x$.  Then $|R_y^a|=|R_x^a|\geq2$, so the previous paragraph gives $y\in P''$; thus, $y\in R_x\cap P''$.
Conversely, suppose $y\in R_x\cap P''$, so that $x=yu$, $y=xv$ and $y=yza$ for some $u,v\in S^1$ and $z\in S$.  Then $x=xwa=yuwa$ and $y=yza=xvza$; since $uwa,vza\in Sa$, it follows that $x\R^ay$, and $y\in R_x^a$ as required.

\pfitem{iii} If $x\not\in P$, then $H_x^a\sub L_x^a=\{x\}$ by (i), and so $H_x^a=\{x\}$.  Similarly, (ii) shows that $H_x^a=\{x\}$ if $x\not\in P''$.  
Finally suppose $x\in P\cap P''$.  Then by (i) and (ii), $H_x^a=L_x^a\cap R_x^a=(L_x\cap P)\cap(R_x\cap P'') = H_x\cap(P\cap P'')$, so it remains to show that $H_x\sub P\cap P''$.  With this in mind, let $y\in H_x$.  Since $y\leqL x$ and $x\in Sa$, it follows that $y\in Sa$.  But then $y\in H_x\cap Sa\sub R_x\cap Sa\sub P$, by Lemma \ref{lem:PP''}(i).  A similar calculation using Lemma~\ref{lem:PP''}(ii) gives $y\in P''$.

\pfitem{iv}  If $x\not\in P$, then $D_x^a=\bigcup_{y\in L_x^a}R_y^a=R_x^a$, since $L_x^a=\{x\}$ by (i).  A similar argument works for $x\not\in P''$.

Finally, suppose $x\in P\cap P''$.  We must show that $D_x^a=D_x\cap P\cap P''$.  We begin with the forwards inclusion.  Clearly $D_x^a\sub D_x$.  Next, note that $R_x^a = R_x\cap P'' \sub R_x\cap Sa \sub P$, by part (ii) above and Lemma~\ref{lem:PP''}(i).  Together with part (i) above, it follows that
\[
D_x^a = \bigcup_{y\in R_x^a}L_y^a = \bigcup_{y\in R_x^a}(L_x\cap P) \sub P. 
\]
Similarly, $L_x^a\sub P''$ and so $D_x^a = \bigcup_{y\in L_x^a}R_y^a = \bigcup_{y\in L_x^a}(R_x\cap P'') \sub P''$.  Thus, $D_x^a\sub D_x\cap P\cap P''$.

To prove the backwards inclusion, suppose $y\in D_x\cap P\cap P''$.  So $y\L z\R x$ for some $z\in S$.  First note that $z\L y$ and $y\in Sa$ together imply that $z\in Sa$.  Since $x\in P$, Lemma \ref{lem:PP''}(i) gives $z\in R_x\cap Sa\sub P$.  But then $z\in L_y\cap P=L_y^a$ by part (i) above, since $y\in P$, and so $z\L^ay$.  Since $y\in P''$, Lemma \ref{lem:PP''}(ii) gives $z\in L_y\cap Sa\sub P''$.  But then $z\in R_x\cap P''=R_x^a$ by part (ii) above, since $x\in P''$, and so $z\R^ax$.  Thus, $y\L^az\R^ax$, which gives $y\D^ax$, and $y\in D_x^a$ as required.

\pfitem{v}  We begin with the backwards inclusion.  Since $D_x^a\sub J_x^a$ for any $x\in Sa$, it suffices to show that $J_x\cap P'\cap P'''\sub J_x^a$ if $x\in P'\cap P'''$.  To do so, suppose $x\in P'\cap P'''$, and let $y\in J_x\cap P'\cap P'''$.  Since $x,y\in P'$, we have
\begin{align*}
x&=uaxv &&\text{and} & y&=u'ayv' &&\text{for some $u,u',v,v'\in S^1$.}
\intertext{In fact, we may assume that $u,u'\in S$; for example, $x=uaxv=ua(uaxv)v=uau(ax)v^2$ with $uau\in S$.  
Since $x,y\in P'''$, we have}
x&=pxqa &&\text{and}& y&=p'yq'a &&\text{for some $p,p'\in S^1$ and $q,q'\in S$.}
\intertext{Since $x\J y$, we have}
x&=syt &&\text{and}& y&=s'xt' &&\text{for some $s,s',t,t'\in S^1$.}
\end{align*}
But then $x=pxqa=p(syt)qa=ps(u'ayv')tqa=(psu'a)y(v'tqa)$.  Since $u',q\in S$, it follows that $psu'a,v'tqa\in Sa$.  Similarly, $y=(p's'ua)x(vt'q'a)$, with $p's'ua,vt'q'a\in Sa$.  It follows that $y\J^ax$, and so $y\in J_x^a$ as required.

To prove the forwards inclusion, let $y\in J_x^a$.  We must show that $y$ belongs to $J_x\cap P'\cap P'''$ if $x\in P'\cap P'''$, or to $D_x^a$ otherwise.  Since this is clearly true if $y=x$, we will assume that $y\not=x$.  Since also $y\J^a x$, it follows that one of (a)--(c) must hold, and also one of (d)--(f):
\bit\bmc2
\itemnit{a} $x=yv$ for some $v\in Sa$,
\itemnit{b} $x=uy$ for some $u\in Sa$,
\itemnit{c} $x=uyv$ for some $u,v\in Sa$,
\itemnit{d} $y=xt$ for some $t\in Sa$,
\itemnit{e} $y=sx$ for some $s\in Sa$,
\itemnit{f} $y=sxt$ for some $s,t\in Sa$.
\emc\eit
It may appear that we need to consider all nine combinations separately.  However, we may reduce to just three.  Indeed, in cases (a), (b), (d) and (e), we respectively define $u=1$, $v=1$, $s=1$ and $t=1$.  Then in all combinations, we have $x=uyv$ and $y=sxt$, with $u,v,s,t\in Sa\cup\{1\}$, and with $\{u,v\}\not=\{1\}$ and $\{s,t\}\not=\{1\}$.  Note that (a) and (d) both hold if and only if $\{u,s\}=\{1\}$, while (b) and (e) both hold if and only if $\{v,t\}=\{1\}$.  For any other combination, we have $x=(usu)y(vtv)$ and $y=(sus)x(tvt)$, with $usu,vtv,sus,tvt\in Sa$, so that (c) and (f) both hold.
Thus, the only combinations we need to consider are:
\[
\text{(a) and (d)} \COMMA \text{(b) and (e)} \COMMA \text{(c) and (f)}.
\]

Suppose first that (a) and (d) both hold, noting then that $x\R^ay$: i.e., $y\in R_x^a$.  If $x\not\in P'\cap P'''$, then we are done, since $y\in R_x^a\sub D_x^a$.  Now suppose $x\in P'\cap P'''$.  Since $y\in J_x^a\sub J_x$, we just need to show that $y\in P'\cap P'''$ as well.  Since $x\in P'$, we have $x=waxz$ for some $w,z\in S^1$, and then $y=xt=(waxz)t=wa(yv)zt=w(ay)vzt$, so that $y\J ay$, and $y\in P'$.  
Since $t\in Sa$, we also have $y=xt=yvt\in S^1ySa$, so that $y\in P'''$.

Next suppose (b) and (e) both hold, noting then that $x\L^ay$: i.e., $y\in L_x^a$.  If $x\not\in P'\cap P'''$, then we are done, since $y\in L_x^a\sub D_x^a$.  Now suppose $x\in P'\cap P'''$.  
Again, we just need to show that $y\in P'\cap P'''$.  Write $u=pa$ where $p\in S$.  Then $y=sx=suy=sp(ay)$, so that $y\in P\sub P'$.  Also, since $x\in P'''$, we have $x=wxza$ for some $w\in S^1$ and $z\in S$.  But then $y=sx=s(wxza)=sw(uy)za\in S^1ySa$, so that $y\in P'''$.

Finally, suppose (c) and (f) both hold.  Since $s,v\in Sa$, we have $s=pa$ and $v=qa$ for some $p,q\in S$.  Now, $x=uyv=u(sxt)v=u(pa)xtv=up(ax)tv$, so that $x\in P'$.  Also, ${x=usxtv=usxt(qa)=(us)x(tq)a\in S^1xSa}$, so that $x\in P'''$.  This shows that $x\in P'\cap P'''$.  A similar argument shows that $y\in P'\cap P'''$.  Since also $y\in J_x^a\sub J_x$, it follows that $y\in J_x\cap P'\cap P'''$, completing the proof in this case.
\epf

By Lemma \ref{lem:P''P'''}, $P''=P'''=Sa$ if $a$ is regular.  Thus, several parts of Theorem \ref{thm:GreenSa} simplify in the case of $a$ being regular.  Since all of our applications involve $a$ (indeed, $S$) being regular, it will be convenient to state this simplification explicitly:

\begin{cor}\label{cor:GreenSa}
Let $S$ be a semigroup, let $a\in \Reg(S)$, and define the sets
\[
P = \set{x\in Sa}{x\L ax} \AND P' = \set{x\in Sa}{x\J ax}.
\]
Then for any $x\in Sa$,
\bit\bmc2
\itemit{i} $L_x^a = \begin{cases} L_x\cap P &\text{if $x\in P$} \\ \{x\} &\text{if $x\not\in P$,} \end{cases}$
\itemit{ii} $R_x^a = R_x\cap Sa$,
\itemit{iii} $H_x^a = \begin{cases} H_x &\hspace{4.9mm}\text{if $x\in P$} \\ \{x\} &\hspace{4.9mm}\text{if $x\not\in P$,} \end{cases}$
\itemit{iv} $D_x^a = \begin{cases} D_x\cap P &\text{if $x\in P$} \\ R_x^a &\text{if $x\not\in P$,} \end{cases}$
\itemit{v} $J_x^a = \begin{cases} J_x\cap P' &\hspace{0.8mm}\text{if $x\in P'$} \\ R_x^a &\hspace{0.8mm}\text{if $x\not\in P'$.} \end{cases}$
\item[] ~
\emc\eit
\end{cor}

\pf
Given the comments before the statement, the only part that is slightly non-obvious is the $x\not\in P'$ case of (v).  Here we have $x\not\in P'=P'\cap Sa=P'\cap P'''$, so Theorem \ref{thm:GreenSa}(v) gives $J_x^a=D_x^a$.  Since $x\not\in P'$, certainly $x\not\in P$, so Theorem \ref{thm:GreenSa}(iv) gives $D_x^a=R_x^a$.
\epf

\subsection[Sandwich-regularity and the structure of $\Reg(Sa)$]{\boldmath Sandwich-regularity and the structure of $\Reg(Sa)$}\label{subsect:P}

We have already seen that the structure of a principal left ideal $Sa$ is easier to describe in the case that the element $a\in S$ is regular; cf.~Theorem \ref{thm:GreenSa} and Corollary \ref{cor:GreenSa}.  In the remaining subsections, we will concentrate exclusively on the case in which $a$ is regular.  In fact, we will identify a natural property, called \emph{sandwich-regularity} in \cite{Sandwiches1}, that allows for an even more detailed analysis.  In all of our motivating examples,~$S$ is itself regular, in which case every element of $S$ is sandwich regular.

We begin with a simple lemma; it shows that if we wish to study $Sa$ with $a$ a regular element of $S$, then we may assume without loss of generality that $a$ is in fact an idempotent.

\begin{lemma}
If $a$ is a regular element of $S$, then $Sa=Se$ for some idempotent $e$ of $S$.
\end{lemma}

\pf
Let $b\in S$ be such that $a=aba$, and define the idempotent $e=ba$.  Then $Sa = Saba \sub Sba \sub Sa$, so that $Sa=Sba=Se$.
\epf

If $a$ is an idempotent of $S$, then we may also consider the \emph{local monoid} $aSa=\set{axa}{x\in S}$, which is the largest monoid contained in $S$ that has $a$ as its (two-sided) identity element.  This monoid $aSa$ will play an important role in all that follows.   
The next result gathers some basic properties that we will need.
We will keep the notation of the previous section, in particular $P=\set{x\in Sa}{x\L ax}$.

\newpage

\begin{lemma}\label{lem:Sa_equiv}
If $a\in E(S)$, then 
\bit
\itemit{i} $aSa=aP\sub P=Pa$,
\itemit{ii} the following are equivalent:
\item[] 
\emph{(a)} $aSa\sub\Reg(S)$, \qquad
\emph{(b)} $P\sub\Reg(S)$, \qquad
\emph{(c)} $\Reg(Sa)=P$, \qquad
\emph{(d)} $aSa$ is a regular monoid.
\eit
\end{lemma}

\pf
(i).  Since $P\sub Sa$, and since $a$ is an idempotent, we clearly have $P=Pa$, and also $aP\sub aSa$.  
To show that $aSa\sub P$ and $aSa\sub aP$, let $x\in aSa$.
Then $x\in Sa$ and also $x=ax$ (as $a$ is an idempotent) so certainly $x\L ax$, which gives $x\in P$.
But then also $x=ax\in aP$ as well.

\pfitem{ii}  Corollary \ref{cor:RegSa} gives (a)$\implies$(b), Theorem \ref{thm:RegSa} gives (b)$\implies$(c), and (d)$\implies$(a) is clear, so it remains only to show that (c)$\implies$(d).  
To do so, suppose (c) holds.  Let $x\in aSa$.  The proof will be complete if we can show that $x$ is regular in~$aSa$.  By part (i), just proved, $x=ay$ for some $y\in P$.  Since $P=\Reg(Sa)$ by assumption, there exists $z\in Sa$ such that $y=yzy$.  Since $y,z\in Sa$, we have $y=ya$ and $z=za=za^2$, and so $x = ay = ayzy = a(ya)(za^2)y = (ay)(aza)(ay) = x(aza)x$, so that $x$ is indeed regular in $aSa$.
\epf

The remainder of this section is devoted to the study of the structure of $\Reg(Sa)$ in the case that $a\in E(S)$ satisfies the conditions of Lemma \ref{lem:Sa_equiv}(ii).  In \cite{Sandwiches1}, a regular element $a\in S$ (not necessarily an idempotent) for which $aSa\sub\Reg(S)$ was called \emph{sandwich-regular}, and we will continue to use that terminology here.
\bit
\item[] {\bf \boldmath For the remainder of this subsection, we fix a sandwich-regular idempotent $a\in E(S)$.}
\eit
Thus, by Corollary \ref{cor:RegSa}, $P=\Reg(Sa)$ is a (regular) subsemigroup of $Sa$, indeed a right ideal.  Thus, we may study $P$ as a semigroup in its own right.
In what follows, we will see that the structure of $P=\Reg(Sa)$ is closely related to that of the regular monoid $aSa$.
In later sections, we will see that when $S$ belongs to a natural family of semigroups, such as full or partial transformation semigroups, the local monoid $aSa$ will be another member of this family.

Lemma \ref{lem:Sa_equiv}(i) says that $aSa$ is a subsemigroup of $P$.  It turns out that $aSa$ is also a natural homomorphic image of $P$, as we will demonstrate in the next lemma.  We will see later that $P$ contains a number of subsemigroups isomorphic to $aSa$; see Remark \ref{rem:MI_P}.

\begin{lemma}\label{lem:phi}
If $a$ is a sandwich-regular idempotent of $S$, then the map $\phi:P\to aSa:x\mt ax$ is a surmorphism.
\end{lemma}

\pf
Since $aSa=aP$, by Lemma \ref{lem:Sa_equiv}(i), $\phi$ is surjective.  To show that $\phi$ is a homomorphism, suppose $x,y\in P$.  Since $a$ is a right identity for $Pa=P$, $x=xa$, and so $(xy)\phi = a(xy) = a(xa)y = (x\phi)(y\phi)$.
\epf

The map
\[
\phi:P\to aSa:x\mt ax
\]
from Lemma \ref{lem:phi} will play a crucial role in all that follows; in particular, we will use $\phi$ to relate many structural properties of $P=\Reg(Sa)$ to corresponding properties of $aSa$.
As a first such application, we show how (products of) idempotents in $Sa$ are related to (products of) idempotents in $aSa$.  Recall that for any semigroup $T$, we write $\bbE(T)=\la E(T)\ra$ for the idempotent-generated subsemigroup of $T$.  Since all idempotents of $Sa$ are regular, and since $P=\Reg(Sa)$, it is clear that $E(Sa)=E(P)$ and $\bbE(Sa)=\bbE(P)$.

\begin{thm}\label{thm:E_Sa}
If $a$ is a sandwich-regular idempotent of the semigroup $S$, then 
\bit\bmc2
\itemit{i} $E(Sa)=E(aSa)\phi^{-1}$,
\itemit{ii} $\bbE(Sa)=\bbE(aSa)\phi^{-1}$.
\emc\eit
\end{thm}

\pf
Since any homomorphism maps (products of) idempotents to (products of) idempotents, it is enough to prove the backwards containments in both parts.  To do so, let $x\in P$; since $P$ is regular, there exists $e\in E(P)=E(Sa)$ such that $x=ex$.

\pfitem{i}  If $ax=x\phi\in E(aSa)$, then $ax=axax$ and so $x=ex=eax=eaxax=exx=xx$, so that $x\in E(Sa)$.

\pfitem{ii}  If $ax=x\phi\in \bbE(aSa)$, then $x=ex=eax\in\bbE(Sa)$, since $e\in E(Sa)$ and $ax\in\bbE(aSa)\sub\bbE(Sa)$.
\epf

In the remainder of the current subsection, we investigate the connection, via $\phi$, between Green's relations on $P$ and $aSa$, leading to a detailed description of $P$ as a kind of ``inflation'' of $aSa$; see Theorem \ref{thm:D_structure_P} and Remark \ref{rem:inflation_Sa}.

Since $P$ is a regular subsemigroup of $Sa$, the $\R$-, $\L$- and $\H$-relations on $P$ are simply the restrictions to~$P$ of the corresponding relations on $Sa$; see for example, \cite[Proposition A.1.16]{RSbook}.  Since $P$ consists of \emph{all} regular elements of $Sa$, \cite[Lemma 2.8]{Sandwiches1} says that this is also the case for the $\D$-relation.  Thus, if $\K$ is any of Green's relations other than $\J$, we will continue to write $\Ka$ for Green's $\K$ relation on $P$; we will also continue to write $K_x^a$ for the $\Ka$-class of $x\in P$ for any such $\K$.  We will write $\J^P$ for Green's $\J$-relation on $P$, and denote $\J^P$-classes by $J_x^P$.  

Together with Corollary \ref{cor:GreenSa}, the previous paragraph may be summarised as follows:

\begin{lemma}\label{lem:Green_P}
If $a$ is a sandwich-regular idempotent of $S$, and if $x\in P$, then
\bit\bmc4
\itemit{i} $L_x^a=L_x\cap P$,
\itemit{ii} $R_x^a=R_x\cap P$,
\itemit{iii} $D_x^a=D_x\cap P$,
\itemit{iv} $H_x^a=H_x$. \epfres
\emc\eit
\end{lemma}

Green's $\J$-relation on $P$ is not as easy to describe.  However, if Green's $\J$ and $\D$ relations on $S$ coincide, then the same is true in $P$ (though it need not be true in $Sa$ itself; see for example Theorem \ref{thm:Green_TXA}):

\begin{cor}\label{cor:J=D_P}
If ${\J}={\D}$ in $S$, then ${\J^P}={\D^a}$ in $P$.
\end{cor}

\pf
Since the $\D$ relation is contained in the $\J$ relation in any semigroup, it suffices to show that ${\J^P}\sub{\D^a}$.  So suppose $x,y\in P$ are such that $(x,y)\in{\J^P}$.  Since $P$ is a subsemigroup of $S$, it follows that $(x,y)\in{\J}={\D}$, and so $y\in D_x$.  But also $x,y\in P$, and so $y\in D_x\cap P=D_x^a$, by Lemma \ref{lem:Green_P}(iii), whence $(x,y)\in{\D^a}$, as required.
\epf

We will also need to refer to Green's relations on the monoid $aSa$.  Again, to avoid confusion, we will use superscripts to identify these relations: the $\K$ relation on $aSa$ will be denoted by $\aKa$, and $\aKa$-classes in~$aSa$ will be denoted by~${}^a\!K_x^a$.  Clearly ${}^a\!K_x^a\sub K_x\cap aSa$ for any $x\in aSa$ and for any $\K$.

\begin{lemma}\label{lem:Green_aSa}
If $a$ is a sandwich-regular idempotent of $S$, and if $x\in aSa$, then
\bit\bmc3
\itemit{i} ${}^a\!L_x^a=L_x\cap aSa$,
\itemit{ii} ${}^a\!R_x^a=R_x\cap aSa$,
\itemit{iii} ${}^a\!D_x^a=D_x\cap aSa$,
\itemit{iv} ${}^a\!J_x^a=J_x\cap aSa$,
\itemit{v} ${}^a\!H_x^a=H_x$.
\emc\eit
\end{lemma}

\pf
(i) and (ii).  These also follow from \cite[Proposition A.1.16]{RSbook} since $aSa$ is a regular subsemigroup of~$S$.

\pfitem{iii}  We noted before the lemma that ${}^a\!D_x^a\sub D_x\cap aSa$.  
To demonstrate the reverse inclusion, let $y\in D_x\cap aSa$.  So $x\L z\R y$ for some $z\in S$.  Then $z=ux=yv$ for some $u,v\in S^1$.  From $z=ux$ and $x\in Sa$, we obtain $z=za$, and similarly $z=az$.  It follows that $z=aza\in aSa$.  But then $z\in L_x\cap aSa={}^a\!L_x^a$ by (i), and similarly $z\in{}^a\!R_y^a$.  Thus, $x\aLa z\aRa y$, so that $x\aDa y$, and $y\in{}^a\!D_x^a$ as required.

\pfitem{iv}  
To show the backwards inclusion (which is again all that is required), let $y\in J_x\cap aSa$.  Since $y\J x$, we have $x=syt$ and $y=uxv$ for some $s,t,u,v\in S^1$.  Since $x,y\in aSa$, we have $x=axa$ and $y=aya$.  It then follows that $x = axa = asyta = as(aya)ta = (asa) y (ata)$, and similarly $y = (aua) x (ava)$.  Since $asa,ata,aua,ava\in aS^1a=aSa$, it follows that $x\aJa y$, and $y\in{}^a\!J_x^a$.

\pfitem{v}  From (i) and (ii), we obtain ${}^a\!H_x^a={}^a\!L_x^a\cap {}^a\!R_x^a = (L_x\cap aSa)\cap(R_x\cap aSa)=H_x\cap aSa$, so it remains to show that $H_x\sub aSa$.  To do so, let $y\in H_x$.  Since $y\L x$ and $x\in Sa$, it follows that $y\in Sa$, and so $y=ya$.  Similarly, $y\R x$ and $x\in aS$ give $y=ay$.  It follows that $y=aya\in aSa$.  As noted above, this completes the proof.
\epf

\begin{rem}\label{rem:H_classes}
Even though the last parts of Lemmas \ref{lem:Green_P} and \ref{lem:Green_aSa} say that $\Ha$-classes of $P$ and $\aHa$-classes of $aSa$ are simply $\H$-classes of $S$, we will continue to use superscripts to indicate whether a certain set of elements is to be thought of as an $\H$-class of $S$, an $\Ha$-class of $P$, or an $\aHa$-class of $aSa$.
\end{rem}

\begin{cor}
If $a$ is a sandwich-regular idempotent of $S$, then the group of units of $aSa$ is ${}^a\!H_a^a = H_a$.
\end{cor}

\pf
The group of units of any monoid is the $\H$-class of the identity element.  Thus, the group of units of $aSa$ is the $\aHa$-class of $a$; by Lemma \ref{lem:Green_aSa}(v), this is ${}^a\!H_a^a=H_a$.
\epf

We now wish to show how the internal structure of a $\Da$-class $D_x^a$ of $P$ is related to that of the corresponding $\aDa$-class ${}^a\!D_{ax}^a=D_{ax}\cap aSa$ of $aSa$.  To do so, we introduce a number of new relations on $P$.
Associated to each of Green's relations $\K$, we define a relation $\Kha$ on $P$ by
\[
{\Kha} = \bigset{(x,y)\in P\times P}{(ax,ay)\in{\aKa}}.
\]
So $\Kha$ is the pre-image under the map $\phi:P\to aSa:x\mt ax$ of the $\aKa$-relation on $aSa$.  Clearly ${\Ka}\sub{\Kha}$ for any $\K$.
Theorem \ref{thm:D_structure_P} (and Remark \ref{rem:inflation_Sa}) gives the promised description of the $\Da$-classes of $P$.  We begin with two technical lemmas.

\begin{lemma}\label{lem:leqJa_leqaJa}
If $x,y\in P$, then 
\bit\bmc2
\itemit{i} $x\leqLa y \iff ax \leq_{\aLa} ay$.
\itemit{ii} $x\leqJP y \iff ax \leqaJa ay$.
\emc\eit
\end{lemma}

\pf
We just prove (ii), as the proof of (i) is similar, but slightly easier.

\pfitem{$\Rightarrow$}  Suppose $x\leqJP y$.  Then, since $P$ is regular, $x=uyv$ for some $u,v\in P$ (not just in $P^1$).  But then $ax=auyv=a(ua)(ya)v = (au) ay (av)$, with $au,av\in aP=aSa$, and so $ax \leqaJa ay$.

\pfitem{$\Leftarrow$}  Suppose $ax \leqaJa ay$.  Then $ax=u(ay)v$ for some $u,v\in aSa$.  Since $P$ is regular, there exists an idempotent $e\in E(P)$ such that $x=ex$.  But then $x=ex=eax=e(uayv)=eua\cdot y\cdot v$; since $e,a\in P$ and since $u,v\in aSa\sub P$ (by Lemma \ref{lem:Sa_equiv}(i)), we have $eua,v\in P$, so that $x\leqJP y$.
\epf

\begin{lemma}\label{lem:Khat_Sa}
We have
\bit\bmc2
\itemit{i} ${\Lha}={\La}$,
\itemit{ii} ${\Ra}\sub{\Rha}\sub{\Da}$,
\itemit{iii} ${\Ha}\sub{\Hha}\sub{\Da}$,
\itemit{iv} ${\Dha}={\Da}\sub{\Jha}={\J^P}$.
\emc\eit
\end{lemma}

\pf
(i).  This follows quickly from Lemma \ref{lem:leqJa_leqaJa}(i).

\pfitem{ii} Clearly ${\Ra}\sub{\Rha}$.  To show that ${\Rha}\sub{\Da}$, let $(x,y)\in{\Rha}$.  Since $P$ is regular, we have $x \Ra e$ and $y\Ra f$ for some $e,f\in E(P)$.  We claim that $e\Da f$, and since then $x\Da e\Da f\Da y$, this will complete the proof of (ii).  To show that $e\Da f$, we will show that $e\Ra ef\La f$.  Since $ef\leqRa e$ and $ef\leqLa f$, it remains to show the reverse inequalities.  
Since ${\Ra}\sub{\Rha}$, we have $e\Rha x\Rha y\Rha f$, so that $ae \aRa af$ (in~$aSa$).  Since $ae,af\in E(aSa)$, it follows that $ae=(af)(ae)$ and $af=(ae)(af)$.  But then $e=ee=(ea)e=e(afae)=efe\leqRa ef$.  Similarly, $f=fef\leqLa ef$.

\pfitem{iii} We have ${\Ha}={\La}\cap{\Ra}\sub{\Lha}\cap{\Rha}={\Hha}$ and ${\Hha}={\Lha}\cap{\Rha}={\La}\cap{\Rha}\sub{\Da}\cap{\Da}={\Da}$.

\pfitem{iv} It is clear that ${\Da}\sub{\Dha}\sub{\Jha}$, and we obtain ${\J^P}={\Jha}$ from Lemma \ref{lem:leqJa_leqaJa}(ii).  It remains only to observe that ${\Dha}={\Lha}\vee{\Rha}\sub{\Da}\vee{\Da}={\Da}$.  
\epf

The next result describes the structure of the $\Hha$-classes of $P$ in terms of left groups (as defined in Section~\ref{subsect:LG}); see also Remark \ref{rem:inflation_Sa}.
Since ${\Ra}\sub{\Rha}$, any $\Rha$-class of $P$ is a union of $\Ra$-classes; thus, if $x\in P$, we may consider the set~$\Rh_x^a/{\Ra}$ of all $\Ra$-classes of $P$ contained in $\Rh_x^a$.  
Recall that if $x\in P$, then the $\Ha$-class of $x$ in $P$ is $H_x^a=H_x$ (by Lemma \ref{lem:Green_P}(iv)), and that the $\aHa$-class of $ax$ in $aSa$ is ${}^a\!H_{ax}^a=H_{ax}$ (by Lemma \ref{lem:Green_aSa}(v)).  However, as in Remark \ref{rem:H_classes}, we will continue to refer to these classes as $H_x^a$ and ${}^a\!H_{ax}^a$, so that it is clear that we are thinking of them as $\Ha$- or $\aHa$-classes of $P$ or $aSa$, respectively.

\begin{thm}\label{thm:D_structure_P}
Let $x\in P$, and let $r=|\Rh_x^a/{\Ra}|$ be the number of $\Ra$-classes contained in $\Rh_x^a$.  Then
\bit
\itemit{i} the restriction to $H_x^a$ of the map $\phi:P\to aSa$ is a bijection $\phi|_{H_x^a}:H_x^a\to {}^a\!H_{ax}^a$,
\itemit{ii} $H_x^a$ is a group if and only if ${}^a\!H_{ax}^a$ is a group, in which case these groups are isomorphic,
\itemit{iii} if $H_x^a$ is a group, then $\Hh_x^a$ is a left group of degree $r$ over $H_x^a$,
\itemit{iv} if $H_x^a$ is a group, then $E(\Hh_x^a)$ is a left zero band of size $r$.
\eit
\end{thm}

\pf
(i).  Since $x\in P$, we have $x\L ax$, and so $x=uax$ for some $u\in S$.  By Green's Lemma \cite[Lemma~2.2.1]{Howie} in the semigroup $S$, it follows that the maps
\[
\th_1:H_x\to H_{ax}:z\mt az \AND \th_2:H_{ax}\to H_x:z\mt uz
\]
are mutually inverse bijections.  But $H_x=H_x^a$ and $H_{ax}={}^a\!H_{ax}^a$, by Lemmas \ref{lem:Green_P}(iv) and \ref{lem:Green_aSa}(v).  Since $\th_1$ has the same action as $\phi$ on $H_x=H_x^a$, it follows that $\phi|_{H_x^a}=\th_1$ is a bijection.

\pfitem{ii}  If $H_x^a$ is a group, then without loss of generality, we may assume that $x$ is an idempotent; but then so too is $ax=x\phi$, and so ${}^a\!H_{ax}^a$ is a group.  Conversely, if ${}^a\!H_{ax}^a$ is a group, then we may assume $ax$ is an idempotent; but then so too is $x$, by Theorem \ref{thm:E_Sa}(i), and so $H_x^a$ is a group.  

By (i), $\phi|_{H_x^a}:H_x^a\to {}^a\!H_{ax}^a$ is a bijection.  If $H_x^a$ is a group, then $\phi|_{H_x^a}$ is also a homomorphism---as it is a restriction of a homomorphism to a sub(semi)group---and hence an isomorphism.

\pfitem{iii) and (iv}  Suppose $H_x^a$ is a group.  Since $\Hh_x^a$ is a union of $\H^a$-classes, we may write $\Hh_x^a=\bigsqcup_{y\in Y}H_y^a$ for some subset $Y\sub P$.  (Here, ``$\sqcup$'' means \emph{disjoint} union.)  By Lemma \ref{lem:Khat_Sa}(i), we have $\Hh_x^a\sub\Lh_x^a=L_x^a$, and so all of the elements of $Y$ are $\La$-related.
For any $y\in Y$, $H_y^a\phi\sub\Hh_x^a\phi={}^a\!H_{ax}^a$, and since ${}^a\!H_{ax}^a$ is a group (and since $\Hh_x^a\phi=\Hh_y^a\phi={}^a\!H_{ay}^a$), part~(ii) says that $H_y^a$ is a group.  Thus, $\Hh_x^a=\bigsqcup_{y\in Y}H_y^a$ is a union of groups, each isomorphic to~$H_x^a$.  Thus, if we can prove (iv), then (iii) will also follow.
Since each $H_y^a$ is a group, we may assume without loss of generality, that each element of $Y$ is an idempotent, so that $Y=E(\Hh_x^a)$.  Now, if $y,z\in Y$, then since~$y\La z$, we have $yz=y$, from which it follows that~$Y=E(\Hh_x^a)$ is a left zero band.  

It remains only to show that $|Y|=r$.  To do so, it suffices to show that $\Rh_x^a = \bigsqcup_{y\in Y}R_y^a$.  First note that since the elements of $Y$ are all $\La$-related but are mutually $\Ha$-unrelated (as they are all idempotents), it follows that they are mutually $\Ra$-unrelated, and so the $\Ra$-classes $R_y^a$ ($y\in Y$) are indeed pairwise disjoint.
Next, consider some $y\in Y$.  Since $y\in \Hh_x^a \sub \Rh_x^a$, and since ${\Ra}\sub{\Rha}$ by Lemma \ref{lem:Khat_Sa}(ii), we have~$R_y^a\sub\Rh_x^a$.  Since this is true for any $y\in Y$, it follows that $\bigsqcup_{y\in Y}R_y^a\sub\Rh_x^a$.  
To prove the reverse containment, suppose~$z\in\Rh_x^a$.  Since ${\Rha}\sub{\Da}$ by Lemma \ref{lem:Khat_Sa}(ii), we have $z\Da x$, and so $R_z^a\cap L_x^a$ is non-empty.  Let $w\in R_z^a\cap L_x^a$ be arbitrary.  Since $z\in\Rh_x^a$ and since ${\Ra}\sub{\Rha}$, we have $w\in R_z^a\sub\Rh_x^a$.  Since also $w\in L_x^a$, it follows that $w \in \Rh_x^a\cap L_x^a = \Rh_x^a\cap\Lh_x^a = \Hh_x^a = \bigsqcup_{y\in Y}H_y^a$, and so $w\in H_y^a \sub R_y^a$ for some $y\in Y$.  Since $w\in R_z^a$, it follows that $z\R^aw\R^ay$, whence $z\in R_y^a\sub\bigsqcup_{y\in Y}R_y^a$, as required.
\epf

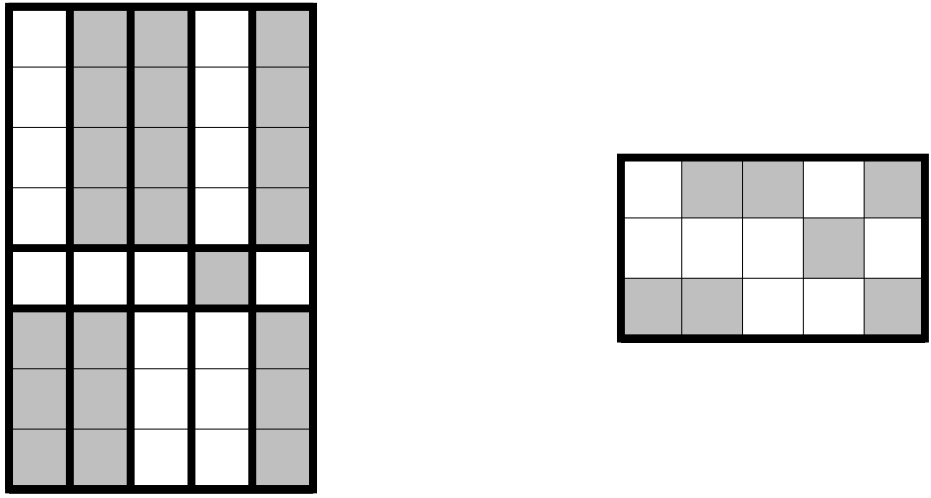
\begin{figure}[ht]
\begin{center}
\scalebox{.8}{
\begin{tikzpicture}[scale=1]
\node (D1) at (0,0) {\DaClass{5}{8}{
1/2,1/3,1/5,
2/2,2/3,2/5,
3/2,3/3,3/5,
4/2,4/3,4/5,
5/4,
6/1,6/2,6/5,
7/1,7/2,7/5,
8/1,8/2,8/5
}{}{1}
{0,1,2,3,4,5}
{0,3,4,8}
};
\node (D2) at (10,0) {\DClass{5}{3}{1/2,1/3,1/5,2/4,3/1,3/2,3/5}{}{1}};
\end{tikzpicture}
}
\end{center}
\vspace{-5mm}
\caption{A ${\Da}$-class of $P=\Reg(Sa)$ (left) and its corresponding $\aDa$-class of $aSa$ (right).  See Remark \ref{rem:inflation_Sa} for more information.}
\label{fig:inflation_Sa}
\end{figure}

\begin{rem}\label{rem:inflation_Sa}
By the preceding series of results, the structure of $P=\Reg(Sa)$, in terms of Green's relations, is a kind of ``inflation'' of the corresponding structure of the regular monoid $aSa$:
\bit
\itemnit{i} The partially ordered sets $(P/{\J^P},\leqJP)$ and $(aSa/{\aJa},\leqaJa)$ are order-isomorphic, via $J_x^P\mt {}^a\!J_{ax}^a$.
\itemnit{ii} The sets $P/{\D^a}$ and $aSa/{\aDa}$ are in one-one correspondence, via $D_x^a\mt {}^a\!D_{ax}^a$.
\itemnit{iii} Each $\Kha$-class in $P$ is a union of $\K^a$-classes.
\itemnit{iv} The $\aRa$-, $\aLa$- and $\aHa$-classes contained within a single $\aDa$-class ${}^a\!D_{ax}^a$ of $aSa$ ($x\in P$) are in one-one correspondence with the $\Rha$-, ${\Lha}={\La}$- and $\Hha$-classes in the ${\Dha}={\Da}$-class $D_x^a$ of $P$. 
\itemnit{v} An $\Hha$-class $\Hh_x^a$ in $P$ is a union of $\H^a$-classes, and these are either all non-groups (if $H_{ax}={}^a\!H_{ax}^a$ is a non-group $\aHa$-class of $aSa$) or else all groups (if $H_{ax}$ is a group); in the latter case, $\Hh_x^a$ is a left group.
\eit
Figure \ref{fig:inflation_Sa} illustrates the last two points in an \emph{egg-box diagram} (as described in Section \ref{subsect:Green}).  The left egg-box displays a ${\Dha}={\D^a}$-class in $P$, and the right egg-box displays the corresponding $\aDa$-class in $aSa$.  Group $\Ha$- and $\aHa$-classes are shaded gray, and solid lines in the left egg-box denote boundaries between ${\Rha}$-classes and ${\Lha}={\La}$-classes.   See also Figures \ref{fig:TXA}--\ref{fig:RXal}.
\end{rem}

\subsection{Rank and idempotent rank}\label{subsect:rank}

This subsection mainly concerns the rank (and idempotent rank, where appropriate) of the regular and idempotent-generated subsemigroups $P=\Reg(Sa)$ and $\bbE(Sa)$ in the case that $a$ is a sandwich-regular idempotent of the semigroup $S$.  (The concepts of (relative) rank and (relative) idempotent rank were defined in Section \ref{subsect:rk}.)  The main results are Theorems \ref{thm:rank_P} and \ref{thm:rank_EP}, which give lower bounds for these (idempotent) ranks, and show that these bounds are exact values in the case that $P$ is RI-dominated.
\bit
\item[] {\bf \boldmath For the duration of this subsection, we fix a sandwich-regular idempotent $a\in E(S)$.}
\eit
We begin by giving numerous characterisations of the mid-identities of the regular semigroup $P=\Reg(Sa)$.  For $x\in P$, we write
\[
V_P(x)=\set{y\in P}{x=xyx,\ y=yxy}
\]
for the set of all inverses of $x$ in $P$.  (The notation is chosen in order to distinguish $V_P(x)$ from the set ${V(x)=\set{y\in S}{x=xyx,\ y=yxy}}$ of all inverses of $x$ in $S$.)

\begin{prop}\label{prop:MI_P}
If $a$ is a sandwich-regular idempotent of a semigroup $S$, then
\[
\MI(Sa)=\RI(Sa)=\MI(P)=\RI(P)=V_P(a)= V(a)\cap P=V(a)\cap Sa =V(a)a=E(\Hh_a^a)=a\phi^{-1} . 
\]
\end{prop}

\pf 
As $a$ is a right identity of both $Sa$ and $P$, Lemma \ref{lem:MI}(iii) gives $\MI(Sa)=\RI(Sa)$ and ${\MI(P)=\RI(P)}$.  We complete the proof by demonstrating a series of set containments. 
\bit
\item
Suppose $u\in\MI(Sa)$.  Since $Sa$ has a right identity, Lemma \ref{lem:MI}(ii) gives $u\in E(Sa)\sub P$.  Clearly $xy=xuy$ for all $x,y\in P$ (since the same is true of all $x,y\in Sa$, as $u\in\MI(Sa)$), so $u\in\MI(P)$.  This shows that $\MI(Sa)\sub\MI(P)$.

\item Next suppose $u\in\RI(P)$.  Since $a$ and $u$ are both right identities, $a=au=aua$ and $u=ua=uau$.  This shows that $\RI(P)\sub V_P(a)$.

\item Since $V_P(a)\sub V(a)$ and $V_P(a)\sub P\sub Sa$, we have $V_P(a)\sub V(a)\cap P\sub V(a)\cap Sa$.

\item Next suppose $u\in V(a)\cap Sa$.  Then $u=ua\in V(a)a$.  This shows that $V(a)\cap Sa\sub V(a)a$.

\item Next suppose $u\in V(a)a$, so $u=va$ for some $v\in V(a)$.  Then $u=va=(vav)a=(va)(va)=u^2$, so $u$ is an idempotent.  We also have
\[
u\phi=au=a(va)=a=a\phi \implies a\phi \aHa u\phi \implies a\Hha u \implies u\in\Hh_a^a.
\]
Thus, $u\in E(\Hh_a^a)$.  This shows that $V(a)a\sub E(\Hh_a^a)$.

\item
Next suppose $u\in E(\Hh_a^a)$.  Then $u\phi\in{}^a\!H_a^a$.  Since $u$ is an idempotent, so too is $u\phi$, and so $u\phi=a$ (as~$a$ is the unique idempotent of the group ${}^a\!H_a^a$), whence $u\in a\phi^{-1}$.  This shows that $E(\Hh_a^a)\sub a\phi^{-1}$.

\item
Finally, suppose $u\in a\phi^{-1}$, so that $a=u\phi=au$.  Then for any $x\in Sa$, $x=xa=xau=xu$, so that $u\in\RI(Sa)$.  This shows that $a\phi^{-1}\sub\RI(Sa)$, and completes the proof. \qedhere
\eit
\epf

\begin{rem}\label{rem:MI_P}
Consider Proposition \ref{prop:MI}, as applied to the (regular) semigroup $P$.  It refers to the local monoids $ePe$, where $e\in\MI(P)$.  Since $\MI(P)=\RI(P)$, by Proposition \ref{prop:MI_P}, each such local monoid is in fact a principal right ideal: $ePe=eP$.  Proposition \ref{prop:MI}(ii) says that each of these local monoids are isomorphic to $aPa=aP$, and Lemma \ref{lem:Sa_equiv}(i) says that $aP=aSa$.
Thus, $P$ generally contains several (local) monoids isomorphic to $aSa$.
Moreover, by Proposition \ref{prop:MI}(iv) and Lemma \ref{lem:RILI}(i), we have $P=\bigcup_{e\in\RP(P)}eP$ if and only if $P$ is RI-dominated.
\end{rem}

Recall that we wish to prove results about the (idempotent) ranks of $P$ and $\bbE(Sa)$; see Theorems \ref{thm:rank_P} and~\ref{thm:rank_EP}.  To prove these theorems, it will be convenient to first prove a more general result; see Proposition~\ref{prop:UW}.  This result concerns submonoids of $aSa$ satisfying certain conditions; these are automatically satisfied by $\bbE(aSa)$, but not always by $aSa$ itself.  In the latter case, the group of units, ${}^a\!H_a^a=H_a$ of $aSa$ plays a crucial role.  For a monoid $U$, we write $G_U$ for the group of units of $U$.  If $U$ is a submonoid of a monoid $M$, then $G_U\sub U\cap G_M$, but we need not have equality (consider the non-negative integers in the additive monoid of all integers).  The next two results concern submonoids $U$ of $aSa$ for which $G_U=U\cap G_{aSa} = U\cap {}^a\!H_a^a$.

\begin{lemma}\label{lem:UW}
Let $a$ be a sandwich-regular idempotent of the semigroup $S$.  Suppose $U$ is a submonoid of~$aSa$ for which $G_U=U\cap{}^a\!H_a^a$, and 
 $U\sm G_U$ is an ideal of $U$.  
Write $\rho = |\Rh_a^a/{\Ra}|$, $W=U\phi^{-1}$ and $T=G_U\phi^{-1}$.  Then 
\bit\bmc2
\itemit{i} $T$ is a left group of degree $\rho$ over $G_U$, 
\itemit{ii} $W\sm T$ is an ideal of $W$.
\emc\eit
\end{lemma}

\pf
(i).  Note first that
$
T=G_U\phi^{-1}=(U\cap{}^a\!H_a^a)\phi^{-1}=W\cap\Hh_a^a
$.
Since $U$ is a submonoid of $aSa$, we have $a\in U$, and so $W$ contains $a\phi^{-1}$; recall that $a\phi^{-1}=E(\Hh_a^a)$ by Proposition \ref{prop:MI_P}.  
For convenience, we will write $F=E(\Hh_a^a)$ for the rest of the proof.  We have just shown that $W$ (and hence $T$) contains $F$.  By Lemma \ref{lem:LG_subs}, since $\Hh_a^a \cong F\times H_a^a$, it follows that $T=FK$ for some submonoid $K$ of $H_a^a$.  Since $K\sub H_a^a$, we have $K=aK$, and so $K=a(aK)=(F\phi)(K\phi)=(FK)\phi=T\phi=G_U$.

\pfitem{ii}  Since $T=W\cap\Hh_a^a$, we may prove this part by showing that for all $x,y\in W$, $xy\in\Hh_a^a \implies x,y\in\Hh_a^a$.  With this in mind, suppose $x,y\in W$ are such that $xy\in\Hh_a^a$.  Then
\begin{align*}
xy\in\Hh_a^a &\implies xy \Hha a 
\implies (ax)(ay) = (x\phi)(y\phi) = (xy)\phi \aHa a\phi = a
\implies (ax)(ay) \in {}^a\!H_a^a\cap U = G_U.
\end{align*}
Since $U\sm G_U$ is an ideal of $U$, it follows that $x\phi=ax$ and $y\phi=ay$ both belong to $G_U\sub {}^a\!H_a^a$.  Thus, $x\phi,y\phi \aHa a=a\phi$, and so $x,y\Hha a$: i.e., $x,y\in\Hh_a^a$.
\epf

\begin{prop}\label{prop:UW}
Let $a$ be a sandwich-regular idempotent of the semigroup $S$.  Suppose $U$ is a submonoid of $aSa$ for which $G_U=U\cap{}^a\!H_a^a$, and $U\sm G_U$ is an ideal of $U$.  Write $\rho = |\Rh_a^a/{\Ra}|$ and $W=U\phi^{-1}$.  Then
\[
\rank(W) \geq \relrank U{G_U} + \max(\rho,\rank(G_U)),
\]
with equality if $P$ is RI-dominated.
\end{prop}

\pf
For convenience, write $T=G_U\phi^{-1}=W\cap\Hh_a^a$.  By Lemma \ref{lem:UW}(ii), $W\sm T$ is an ideal of $W$.  Thus, by Lemma~\ref{lem:rankWT},
\[
\rank(W) = \relrank WT + \rank(T).
\]
By Lemmas \ref{lem:rank_left_group} and \ref{lem:UW}(i), we have $\rank(T)=\max(\rho,\rank(G_U))$.  Thus, it remains to show that
\bit
\itemnit{i} $\relrank WT\geq\relrank U{G_U}$, and
\itemnit{ii} $\relrank WT=\relrank U{G_U}$ if $P$ is RI-dominated.
\eit
(i).  If $X\sub W$ is such that $W=\la T\cup X\ra$ and $|X|=\relrank WT$, then $U=W\phi=\la T\phi\cup X\phi\ra=\la G_U\cup X\phi\ra$, and so $\relrank U{G_U} \leq|X\phi|\leq|X|=\relrank WT$.

\pfitem{ii}  Suppose now that $P$ is RI-dominated.  By (i), it remains to show that $\relrank WT\leq\relrank U{G_U}$.
To do so, let $Y\sub U$ be such that $U=\la G_U\cup Y\ra$ and $|Y|=\relrank U{G_U}$.  Let $Z\sub W$ be such that $Z\phi= Y$ and $|Z|=|Y|$.  
For each $y\in G_U\cup Y$, let $z_y\in T\cup Z$ be such that $y=z_y\phi=az_y$.  Now let $w\in W$ be arbitrary.  Then $aw=w\phi\in U$, and so $aw=y_1\cdots y_k=(az_{y_1})\cdots(az_{y_k})=a(z_{y_1}\cdots z_{y_k})$ for some $y_1,\ldots,y_k\in G_U\cup Y$.  Since $P$ is RI-dominated, $w\leqRa e$ for some $e\in\RI(P)$.  Since $e\in E(P)$, it follows that $w=ew$, and so
\[
w=ew=eaw=ea(z_{y_1}\cdots z_{y_k})=e(z_{y_1}\cdots z_{y_k}).
\]
But the $z_{y_i}$ all belong to $T\cup Z$, and by Proposition \ref{prop:MI_P}, $e\in\RI(P)=E(\Hh_a^a)\sub W\cap\Hh_a^a=T$, so it follows that $w\in\la T\cup Z\ra$.  Thus, $W=\la T\cup Z\ra$, and so $\relrank WT\leq|Z|=|Y|=\relrank U{G_U}$, as required.
\epf

The hypotheses of Proposition \ref{prop:UW} are clearly satisfied by $U=aSa$ as long as $aSa\sm {}^a\!H_a^a$ is an ideal of~$aSa$, so we immediately obtain the following.

\begin{thm}\label{thm:rank_P}
Let $a$ be a sandwich-regular idempotent of the semigroup $S$, write $\rho = |\Rh_a^a/{\Ra}|$, and suppose $aSa\sm {}^a\!H_a^a$ is an ideal of $aSa$.  Then
\[
\rank(P) \geq \relrank {aSa}{{}^a\!H_a^a} + \max(\rho,\rank({}^a\!H_a^a)),
\]
with equality if $P$ is RI-dominated.  \epfres
\end{thm}

Next, we wish to apply Proposition \ref{prop:UW} to $U=\bbE(aSa)$, and also prove a corresponding statement concerning \emph{idempotent} ranks.  To do so, we require the following two lemmas; the first is \cite[Lemma~3.9]{Sandwiches1}, and the second is part of \cite[Lemma 2.1(iv)]{IBM}.

\begin{lemma}\label{lem:IGU}
If $U$ is an idempotent-generated monoid with identity $e$, then
\bit\bmc2
\itemit{i} $G_U=\{e\}$,
\itemit{ii} $U\sm G_U$ is an ideal of $U$,
\itemit{iii} $\rank(U)=1+\relrank U{G_U}$,
\itemit{iv} $\idrank(U)=1+\relidrank U{G_U}$.  \epfres
\emc\eit
\end{lemma}

\begin{lemma}\label{lem:IGU2}
If $M$ is a monoid with identity $e$, then $\bbE(M)\cap G_M=\{e\}$.  \epfres
\end{lemma}

\begin{thm}\label{thm:rank_EP}
Let $a$ be a sandwich-regular idempotent of the semigroup $S$, and write $\rho = |\Rh_a^a/{\Ra}|$.  Then
\[
\rank(\bbE(Sa)) \geq \rank(\bbE(aSa))+ \rho - 1
\AND
\idrank(\bbE(Sa)) \geq \idrank(\bbE(aSa))+ \rho - 1,
\]
with equality in both if $P$ is RI-dominated.  
\end{thm}

\pf
Put $U=\bbE(aSa)$ and $W=U\phi^{-1}$.  Then $W=\bbE(Sa)$, by Theorem \ref{thm:E_Sa}(ii).  By Lemma \ref{lem:IGU}(ii), $U\sm G_U$ is an ideal of $U$.  By Lemmas \ref{lem:IGU}(i) and \ref{lem:IGU2}, we also have $G_U=\{a\}=\bbE(aSa)\cap G_{aSa}=U\cap {}^a\!H_a^a$.  Obviously $U$ is a submonoid of $aSa$.  So by Proposition \ref{prop:UW}, Lemma \ref{lem:IGU}(iii), and the fact that $\rank(G_U)=\rank(\{a\})=1$, it follows that
\[
\rank(W) \geq \relrank{U}{G_U} + \max(\rho,\rank(G_U)) = \rank(U) - 1 +\rho,
\]
with equality throughout if $P$ is RI-dominated.  

For the statement concerning idempotent ranks, consider the proof of Proposition \ref{prop:UW} in the case that $U=\bbE(aSa)$.  First, since $G_U=\{a\}$ by Lemma \ref{lem:IGU}(i), we have $T=a\phi^{-1}=E(\Hh_a^a)$ by Proposition \ref{prop:MI_P}.  By Lemma \ref{lem:UW}(ii), $W\sm T$ is an ideal of $W$.  Lemma \ref{lem:rankWT} then gives 
\[
\idrank(W) = \relidrank WT + \idrank(T).
\]
Since $T$ is a left zero band of size $\rho$, we have $\idrank(T)=\rho$.  As in the proof of Proposition~\ref{prop:UW}, we may show that:
\bit
\itemnit{i} If $X\sub E(W)$ is such that $W=\la T\cup X\ra$ and $|X|=\relidrank WT$, then $U=\la G_U\cup X\phi\ra$.
\itemnit{ii} If $P$ is RI-dominated, and if $Y\sub E(U)$ is such that $U=\la G_U\cup Y\ra$ and $|Y|=\relidrank U{G_U}$, then there exists $Z\sub W$ with $|Z|=|Y|$, $Z\phi=Y$ and $W=\la T\cup Z\ra$; since $Y\sub E(U)$, Theorem \ref{thm:E_Sa}(i) gives $Z\sub E(W)$.
\eit
From (i), and using Lemma \ref{lem:IGU}(iv), it follows that
\[
\relidrank WT=|X|\geq|X\phi|\geq \relidrank U{G_U}=\idrank(U)-1.
\]
Similarly,~(ii) and Lemma \ref{lem:IGU}(iv) give $\relidrank WT\leq|Z|=|Y|=\relidrank U{G_U}=\idrank(U)-1$ if $P$ is RI-dominated.
\epf

Now that we have explored the structure of $P=\Reg(Sa)$ in more detail, we can prove a result concerning the idempotent-generated subsemigroup $\bbE(Sa)$ of $Sa$ in a particular special case that arises in all our motivating examples.
By Lemma \ref{lem:IGU2}, if $M$ is a monoid with identity $e$, then $\bbE(M)\sub\{e\}\cup(M\sm G_M)$.  In particular, $\bbE(aSa)\sub\{a\}\cup(aSa\sm {}^a\!H_a^a)$; the next result describes the situation in which $\bbE(aSa)=\{a\}\cup(aSa\sm {}^a\!H_a^a)$.

\begin{prop}\label{prop:singular_ESa}
Suppose $a$ is a sandwich-regular idempotent of the semigroup $S$, and that $\bbE(aSa)=\{a\}\cup(aSa\sm {}^a\!H_a^a)$.  Then $\bbE(Sa)=a\phi^{-1}\cup(P\sm\Hh_a^a)=E(\Hh_a^a)\cup(P\sm\Hh_a^a)$.
\end{prop}

\pf
By Theorem \ref{thm:E_Sa}(ii), We have $\bbE(Sa) = \bbE(aSa)\phi^{-1} = a\phi^{-1} \cup (aSa\sm {}^a\!H_a^a)\phi^{-1} = a\phi^{-1} \cup (P\sm\Hh_a^a)$.
\epf

\begin{rem}
Note that the set $a\phi^{-1}=E(\Hh_a^a)$ has many equivalent formulations; cf.~Proposition \ref{prop:MI_P}.
\end{rem}

\subsection{Inverse monoids}\label{subsect:inverse}

We continue to assume that $a$ is a sandwich-regular idempotent of $S$.  Recall that for $x\in S$, we write $V(x)=\set{y\in S}{x=xyx,\ y=yxy}$ for the set of all inverses of $x$ in $S$.  Recall also that if $x\in P=\Reg(S)$, we write $V_P(x)=\set{y\in P}{x=xyx,\ y=yxy}$ for the set of all inverses of $x$ in $P$; of course $V_P(x)=V(x)\cap P\sub V(x)$ for any such $x$.  

In \cite{Sandwiches1}, an element $x\in S$ was called \emph{uniquely regular} if $|V(x)|=1$.  Thus, a semigroup is \emph{inverse} if every element is uniquely regular.  In \cite{Sandwiches1}, an element $a\in S$ was called \emph{uniquely sandwich-regular} if each element of $\{a\}\cup aSa$ is uniquely regular.  Every element of an inverse semigroup is uniquely sandwich regular.

\begin{thm}\label{thm:inverse_P}
If $a$ is a uniquely sandwich-regular idempotent of the semigroup $S$, then $\Reg(Sa)=P=aSa$ is an inverse monoid.
\end{thm}

\pf
By Lemma \ref{lem:Sa_equiv}(ii), $P=\Reg(Sa)$.
Let $x\in P=\Reg(Sa)$ be arbitrary, and let $y\in V_P(x)$.  It is easy to check that $x$ and $ax$ both belong to $V(ay)$.  But $ay\in aP=aSa$ is uniquely regular, so it follows that $x=ax$.  Since $x\in P$ was arbitrary, it follows that $P=aP=aSa$.  

To show that $P=aSa$ is inverse, let $x\in P$ be arbitrary.  We must show that $|V_P(x)|=1$.  Since $P$ is regular, certainly $|V_P(x)|\geq1$.  Since $V_P(x)\sub V(x)$ and $|V(x)|=1$ (by the uniquely sandwich-regularity assumption), the proof is complete.
\epf

\begin{rem}\label{rem:inverse_P}
In the case that $a$ is uniquely sandwich-regular, many of the results in the preceeding subsections become trivial or even vacuous, as $\phi:P\to aSa=P:x\mt ax$ is just the identity map.  For example, the $\Kha$ relations are precisely the $\Ka$ relations, and these are the same as the $\aKa$ relations.  Also, in Theorems \ref{thm:rank_P} and \ref{thm:rank_EP}, we have $\rho=1$.  Theorem \ref{thm:rank_P} reduces to the statement
\[
\rank(aSa)=\relrank{aSa}{{}^a\!H_a^a}+\rank({}^a\!H_a^a) \qquad\text{if $aSa\sm{}^a\!H_a^a$ is an ideal of $aSa$,}
\]
which is just a special case of Lemma \ref{lem:rankWT}.  Since $\bbE(Sa)=\bbE(P)=\bbE(aSa)$, Theorem \ref{thm:rank_EP} becomes completely vacuous.
\end{rem}

\section{Principal right ideals}\label{sect:PRI}

In this section, we describe the corresponding results for a principal \emph{right} ideal $aS$ generated by an element~$a$ of the semigroup $S$.  These results a direct duals of those in Section \ref{sect:PLI}, so we will not provide any proofs.  We will also only state the main results.

We begin with a description of the regular elements of $aS$.  The next result is the dual of Theorem \ref{thm:RegSa} and Corollary \ref{cor:RegSa}.

\begin{thm}\label{thm:RegaS}
Let $S$ be a semigroup, let $a\in S$, and define $Q = \set{x\in aS}{x\R xa}$.  Then
\[
\Reg(aS) = \Reg(S) \cap Q.
\]
If $aSa\sub\Reg(S)$, then $Q\sub\Reg(S)$.  Consequently, $\Reg(aS)=Q$ is a left ideal of $aS$ in this case.  \epfres
\end{thm}

We may also describe Green's relations on $aS$ (cf.~Theorem \ref{thm:GreenSa}).  We denote the $\K$ relation on $aS$ by~$\aK$, write ${}^a\!K_x$ for the $\aK$-class of $x$ in $aS$, and so on.  

\newpage

\begin{thm}\label{thm:GreenaS}
Let $S$ be a semigroup, let $a\in S$, and define the sets
\[
Q = \set{x\in aS}{x\R xa} \COMma Q' = \set{x\in aS}{x\J xa} \COMma Q'' = \set{x\in S}{x\in aSx} \COMma Q''' = \set{x\in S}{x\in aSxS^1}.
\]
Then for any $x\in aS$,
\bit\bmc2
\itemit{i} ${}^a\!R_x = \begin{cases} R_x\cap Q &\hspace{1.2mm}\text{if $x\in Q$} \\ \{x\} &\hspace{1.2mm}\text{if $x\not\in Q$,} \end{cases}$
\itemit{ii} ${}^a\!L_x = \begin{cases} L_x\cap Q'' &\text{if $x\in Q''$} \\ \{x\} &\text{if $x\not\in Q''$,} \end{cases}$
\itemit{iii} ${}^a\!H_x = \begin{cases} H_x &\hspace{7.0mm}\text{if $x\in Q\cap Q''$} \\ \{x\} &\hspace{7.0mm}\text{if $x\not\in Q\cap Q''$,} \end{cases}$
\itemit{iv} ${}^a\!D_x = \begin{cases} D_x\cap Q\cap Q'' &\text{if $x\in Q\cap Q''$} \\ {}^a\!L_x &\text{if $x\not\in Q$} \\ {}^a\!R_x &\text{if $x\not\in Q''$,} \end{cases}$
\itemit{v} ${}^a\!J_x = \begin{cases} J_x\cap Q'\cap Q''' &\hspace{0.6mm}\text{if $x\in Q'\cap Q'''$} \\ {}^a\!D_x &\hspace{0.6mm}\text{if $x\not\in Q'\cap Q'''$.} \end{cases}$
\item[] ~ \epfres
\emc\eit
\end{thm}

As in Corollary \ref{cor:GreenSa}, the situation is simpler if the element $a$ is regular, as then $Q''=Q'''=aS$ (cf.~Lemma~\ref{lem:P''P'''}).

\begin{cor}\label{cor:GreenaS}
Let $S$ be a semigroup, let $a\in \Reg(S)$, and define the sets
\[
Q = \set{x\in aS}{x\R xa} \AND Q' = \set{x\in aS}{x\J xa}.
\]
Then for any $x\in aS$,
\bit\bmc2
\itemit{i} ${}^a\!R_x = \begin{cases} R_x\cap Q &\text{if $x\in Q$} \\ \{x\} &\text{if $x\not\in Q$,} \end{cases}$
\itemit{ii} ${}^a\!L_x = L_x\cap aS$,
\itemit{iii} ${}^a\!H_x = \begin{cases} H_x &\hspace{5.8mm}\text{if $x\in Q$} \\ \{x\} &\hspace{5.8mm}\text{if $x\not\in Q$,} \end{cases}$
\itemit{iv} ${}^a\!D_x = \begin{cases} D_x\cap Q &\text{if $x\in Q$} \\ {}^a\!L_x &\text{if $x\not\in Q$,} \end{cases}$
\itemit{v} ${}^a\!J_x = \begin{cases} J_x\cap Q' &\hspace{1.1mm}\text{if $x\in Q'$} \\ {}^a\!L_x &\hspace{1.1mm}\text{if $x\not\in Q'$.} \end{cases}$
\item[] \epfres
\emc\eit
\end{cor}

Again, if $a$ is a regular element of $S$, then $aS=eS$ for some idempotent $e$ of $S$; thus, when studying $aS$ with $a$ regular, we may assume that $a$ is in fact an idempotent.  As with Lemma \ref{lem:Sa_equiv}, we have
\[
aSa=Qa\sub Q=aQ \AND \text{$a$ is sandwich-regular} \iff Q\sub\Reg(S) \iff \Reg(aS)=Q.
\]
\bit
\item[] {\bf \boldmath For the remainder of this subsection, we fix a sandwich-regular idempotent $a\in E(S)$.}
\eit
We again have a surmorphism
\[
\psi:Q\to aSa:x\mt xa,
\]
which allows us to link the structure of $Q=\Reg(aS)$ with that of the regular monoid $aSa$.  The idempotents $E(aS)$ and the idempotent-generated subsemigroup $\bbE(aS)$ of $aS$ may quickly be described; cf.~Theorem \ref{thm:E_Sa}.

\begin{thm}\label{thm:E_aS}
If $a$ is a sandwich-regular idempotent of the semigroup $S$, then 
\bit\bmc2
\itemit{i} $E(aS)=E(aSa)\psi^{-1}$,
\itemit{ii} $\bbE(aS)=\bbE(aSa)\psi^{-1}$. \epfres
\emc\eit
\end{thm}

Green's non-$\J$ relations on $Q$ are also easily characterised.  These are simply the restrictions to $Q$ of the corresponding relations on $aS$, and will also be denoted by $\aK$, with the $\J$-relation denoted by~$\J^Q$; cf.~Lemma \ref{lem:Green_P}.

\begin{lemma}\label{lem:Green_Q}
If $a$ is a sandwich-regular idempotent of $S$, and if $x\in Q$, then
\bit\bmc4
\itemit{i} ${}^a\!L_x=L_x\cap Q$,
\itemit{ii} ${}^a\!R_x=R_x\cap Q$,
\itemit{iii} ${}^a\!D_x=D_x\cap Q$,
\itemit{iv} ${}^a\!H_x=H_x$. \epfres
\emc\eit
\end{lemma}

\begin{cor}\label{cor:J=D_Q}
If ${\J}={\D}$ in $S$, then ${\J^Q}={\aD}$ in $Q$. \epfres
\end{cor}

To describe the internal structure of a $\aD$-class of $Q=\Reg(aS)$, we use the $\aKh$ relations for each of Green's relations $\K$, defined by
\[
\aKh = \bigset{(x,y)\in Q\times Q}{(xa,ya)\in{\aKa}}.
\]
(Recall that $\aKa$ is the $\K$-relation on the monoid $aSa$.)  As in Lemma \ref{lem:Khat_Sa}, we have the following:

\newpage

\begin{lemma}\label{lem:Khat_aS}
We have
\bit\bmc2
\itemit{i} ${\aL}\sub{\aLh}\sub{\aD}$,
\itemit{ii} ${\aRh}={\aR}$,
\itemit{iii} ${\aH}\sub{\aHh}\sub{\aD}$,
\itemit{iv} ${\aDh}={\aD}\sub{\aJh}={\J^Q}$.  \epfres
\emc\eit
\end{lemma}

We then obtain the following analogue of Theorem \ref{thm:D_structure_P}.

\begin{thm}\label{thm:D_structure_Q}
Let $x\in Q$, and let $l=|{}^a\!\Lh_x/{\aL}|$ be the number of $\aL$-classes contained in ${}^a\!\Lh_x$.  Then
\bit
\itemit{i} the restriction to ${}^a\!H_x$ of the map $\psi:Q\to aSa$ is a bijection $\psi|_{{}^a\!H_x}:{}^a\!H_x\to {}^a\!H_{xa}^a$,
\itemit{ii} ${}^a\!H_x$ is a group if and only if ${}^a\!H_{xa}^a$ is a group, in which case these groups are isomorphic,
\itemit{iii} if ${}^a\!H_x$ is a group, then ${}^a\!\Hh_x$ is a right group of degree $l$ over ${}^a\!H_x$,
\itemit{iv} if ${}^a\!H_x$ is a group, then $E({}^a\!\Hh_x)$ is a right zero band of size $l$.  \epfres
\eit
\end{thm}

As in Remark \ref{rem:inflation_Sa}, the Green's structure of $Q=\Reg(aS)$ may be thought of as a kind of ``inflation'' of that of $aSa$.  We leave the reader to supply the details, and to draw a diagram akin to Figure \ref{fig:inflation_Sa} (in $Q$, the ``stretching'' happens in the horizontal direction, rather than the vertical, as in $P=\Reg(Sa)$); compare Figures \ref{fig:TXal} and \ref{fig:RXal} to Figures \ref{fig:TXA} and \ref{fig:RXA}.

RI-domination played an important role in the further study of $P=\Reg(Sa)$ in Section \ref{sect:PLI}, but when studying $aS$ and $Q=\Reg(aS)$, it is \emph{LI}-domination that plays the corresponding role.  The next result, analogous to Proposition \ref{prop:MI_P}, gives several characterisations of the left identities (equivalently, mid-identities) in $aS$ and $Q=\Reg(aS)$.

\begin{prop}\label{prop:MI_Q}
If $a$ is a sandwich-regular idempotent of a semigroup $S$, then
\[
\epfreseq
\MI(aS)=\LI(aS)=\MI(Q)=\LI(Q)=V_Q(a)=V(a)\cap Q=V(a)\cap aS=aV(a)=E({}^a\!\Hh_a)=a\psi^{-1}. 
\]
\end{prop}

After proving an intermediate result analogous to Proposition \ref{prop:UW}, we obtain the following two results concerning the rank (and idempotent rank if appropriate) of $Q=\Reg(aS)$ and $\bbE(aS)$.

\begin{thm}\label{thm:rank_Q}
Let $a$ be a sandwich-regular idempotent of the semigroup $S$, write $\lam = |{}^a\!\Lh_a/{\aL}|$, and suppose $aSa\sm {}^a\!H_a^a$ is an ideal of $aSa$.  Then
\[
\rank(Q) \geq \relrank {aSa}{{}^a\!H_a^a} + \max(\lam,\rank({}^a\!H_a^a)),
\]
with equality if $Q$ is LI-dominated.  \epfres
\end{thm}

\begin{thm}\label{thm:rank_EQ}
Let $a$ be a sandwich-regular idempotent of the semigroup $S$, and write $\lam = |{}^a\!\Lh_a/{\aL}|$.  Then
\[
\rank(\bbE(aS)) \geq \rank(\bbE(aSa))+ \lam - 1
\AND
\idrank(\bbE(aS)) \geq \idrank(\bbE(aSa))+ \lam - 1,
\]
with equality in both if $Q$ is LI-dominated.  \epfres
\end{thm}

We also have the following; cf.~Proposition \ref{prop:singular_ESa}.

\begin{prop}\label{prop:singular_EaS}
Suppose $a$ is a sandwich-regular idempotent of the semigroup $S$, and that $\bbE(aSa)=\{a\}\cup(aSa\sm {}^a\!H_a^a)$.  Then $\bbE(aS)=a\psi^{-1}\cup(Q\sm{}^a\!\Hh_a)=E({}^a\!\Hh_a)\cup(Q\sm{}^a\!\Hh_a)$. \epfres
\end{prop}

As in Theorem~\ref{thm:inverse_P}, the whole theory simplifies significantly if $a$ is uniquely sandwich-regular.

\begin{thm}\label{thm:inverse_Q}
If $a$ is a uniquely sandwich-regular idempotent of the semigroup $S$, then $\Reg(aS)=aSa$ is an inverse monoid.  \epfres
\end{thm}

\begin{rem}
Theorems \ref{thm:inverse_P} and \ref{thm:inverse_Q} together say that when $a$ is uniquely sandwich-regular, we have $\Reg(Sa)=\Reg(aS)=aSa$.
\end{rem}

\section{Full transformation semigroups}\label{sect:TX}

In this section, we apply the general theory developed above to the principal one-sided ideals of the full transformation semigroups.  We will see in Proposition \ref{prop:TXa_aTX} that these one-sided ideals are certain well-known semigroups of restricted range or kernel.  These semigroups of restricted transformations have been studied by several authors \cite{SS2008,MGS2010,Sanwong2011,SS2013,FS2014}.  For example, Green's relations and the regular elements have been described in \cite{SS2008,MGS2010}; these descriptions may be quickly deduced from the general results of Sections~\ref{sect:PLI} and~\ref{sect:PRI}.
Some results concerning ranks of various semigroups we consider may be found in the literature; where possible, these have been acknowledged in the text.  Many other results presented in this section are new.

For the duration of this section, we fix a non-empty set $X$ (which may be finite or infinite), and denote by $\T_X$ the full transformation semigroup over $X$, as defined in Section \ref{subsect:trans}.  We also fix a transformation $a\in\T_X$, with the intention of studying the principal one-sided ideals
\[
\T_Xa=\set{fa}{f\in\T_X} \AND a\T_X = \set{af}{f\in\T_X}.
\]
Since $\T_X$ is regular, by Theorem \ref{thm:T}, we may assume without loss of generality that $a$ is an idempotent.  Using the notation described at the end of Section \ref{subsect:trans}, we will write
\[
a=\tbinom{A_i}{a_i}_{i\in I} \COMMA A=\im(a) \COMMA \al=\ker(a).
\]
So $A=\set{a_i}{i\in I}$, and $\al$ has equivalence classes $\set{A_i}{i\in I}$; since $a$ is an idempotent, we have $a_i\in A_i$ for each $i$.  Since $\T_X$ is regular, $a$ is sandwich-regular, meaning that the theory developed in Sections \ref{sect:PLI} and~\ref{sect:PRI} apply to the principal one-sided ideals $\T_Xa$ and $a\T_X$.  The next result follows quickly from parts (i) and (ii) of Theorem \ref{thm:T}.

\begin{prop}\label{prop:TXa_aTX}
Let $X$ be a non-empty set, let $a\in\T_X$, and write $A=\im(a)$ and $\al=\ker(a)$.  Then
\[
\epfreseq
\T_Xa = \set{f\in\T_X}{\im(f)\sub A} \AND a\T_X = \set{f\in\T_X}{\ker(f)\supseteq\al}.
\]
\end{prop}

The semigroups in Proposition \ref{prop:TXa_aTX} are commonly denoted in the literature by
\[
\TXA = \set{f\in\T_X}{\im(f)\sub A} \AND \TXal = \set{f\in\T_X}{\ker(f)\supseteq\al},
\]
and we will continue to use this notation here.  It is easy to see that
\[
|\TXA|=|A|^{|X|} \AND |\TXal| = |X|^{\Vert\al\Vert}.
\]
Most results of this section will be stated in terms of $A$ or $\al$ without reference to the other, but we will always have the transformation $a$ (which links $A$ and $\al$) in mind.

\subsection[Green's relations and regular elements in $\TXA$ and $\TXal$]{\boldmath Green's relations and regular elements in $\TXA$ and $\TXal$}\label{subsect:Green_TX}

Since $\T_X$ is regular, Green's relations and regular elements in its principal one-sided ideals are governed by the sets $P$, $P'$, $Q$ and $Q'$, as defined in Sections \ref{sect:PLI} and \ref{sect:PRI}.
(Regularity of $\T_X$ ensures that we do not need to explicitly refer to the sets $P''$, $P'''$, $Q''$ and $Q'''$; see Lemma \ref{lem:P''P'''} and its dual.)  To describe these sets, we first recall some terminlogy.  Let $B$ be a subset of $X$, and $\si$ an equivalence on $X$.  We say that
\bit
\item $B$ \emph{saturates} $\si$ if every $\si$-class contains at least one element of $B$,
\item $\si$ \emph{separates} $B$ if every $\si$-class contains at most one element of $B$,
\item $B$ is a \emph{cross section} of $\si$ if every $\si$-class contains exactly one element of $B$.
\eit
Recall that $\Vert\si\Vert$ denotes the number of $\si$-classes of $X$.  If $f\in\T_X$, we write
\[
\si f^{-1}=\set{(x,y)\in X\times X}{(xf,yf)\in\si}.
\]
If $f,g\in\T_X$, then $\ker(fg)=\ker(g)f^{-1}$.

\newpage

\begin{prop}\label{prop:PQP'Q'_T}
Let $X$ be a non-empty set, let $A$ be a non-empty subset of $X$, and let $\al$ be an equivalence relation on $X$.  Then 
\bit
\itemit{i} $\Reg(\TXA)=P=\set{f\in\TXA}{A\text{ saturates }\ker(f)}$ is a right ideal of $\TXA$, 
\itemit{ii} $\Reg(\TXal)=Q=\set{f\in\TXal}{\al\text{ separates }\im(f)}$ is a left ideal of $\TXal$,
\itemit{iii} $P'=\set{f\in\TXA}{|Af|=\rank(f)}$,
\itemit{iv} $Q'=\set{f\in\TXal}{\Vert\al f^{-1}\Vert=\rank(f)}$.
\eit
\end{prop}

\pf
(i) and (ii).  We have $\Reg(\TXA)=\Reg(\T_Xa)=P$ and $\Reg(\TXal)=\Reg(a\T_X)=Q$ from Corollary \ref{cor:RegSa} and Theorem \ref{thm:RegaS}, since $\T_X$ is regular.  Now consider some $f\in\T_X$, and write $f=\binom{F_j}{f_j}$.  Then by Theorem \ref{thm:T}(iv),
\begin{align*}
f\L af \iff \im(f)=\im(af) \iff f_j\in\im(af) \ (\forall j) \iff F_j\cap A\not=\emptyset \ (\forall j) \iff A\text{ saturates }\ker(f).
\end{align*}
Similarly, one may show that $f\R fa \iff \al$ separates $\im(f)$.  

\pfitem{iii}  If $f\in\TXA$, then
\[
f\in P' \iff f\J af \iff \rank(f)=\rank(af) = |{\im(af)}| = |{\im(a)f}| = |Af|.
\]
(iv).  If $f\in\TXal$, then
\[
f\in Q' \iff f\J fa \iff \rank(f)=\rank(fa) = \Vert{\ker(fa)}\Vert = \Vert{\ker(a)f^{-1}}\Vert = \Vert\al f^{-1}\Vert. \qedhere
\]
\epf

\begin{rem}
Since $Af=\set{f_j}{F_j\cap A\not=\emptyset}$, it is clear that $A$ saturates $\ker(f)$ if and only if $Af=\im(f)$.  Thus, we have the alternative characterisation $\Reg(\TXA)=\set{f\in\TXA}{\im(f)=Af}$.  With this in mind, we see that Proposition \ref{prop:PQP'Q'_T}(i) is \cite[Lemma 2.2 and Theorem 2.4]{SS2008}.  
Proposition \ref{prop:PQP'Q'_T}(ii) is \cite[Theorem~2.3]{MGS2010}; in \cite{MGS2010}, the term ``partial cross-section'' was used to describe a set separated by an equivalence relation.
\end{rem}

We now use Corollary \ref{cor:GreenSa}, Proposition \ref{prop:PQP'Q'_T} and Theorem~\ref{thm:T} to give descriptions of Green's relations on~$\TXA=\T_Xa$.

\begin{thm}\label{thm:Green_TXA}
Let $X$ be a non-empty set, let $A$ be a non-empty subset of $X$, and let $f,g\in\TXA$.  Then in the semigroup~$\TXA$,
\bit
\itemit{i} $f\L g \iff f=g$ or $[\im(f)=\im(g)$ and $A$ saturates both $\ker(f)$ and $\ker(g)]$,
\itemit{ii} $f\R g \iff \ker(f)=\ker(g)$,
\itemit{iii} $f\H g \iff f=g$ or $[\im(f)=\im(g)$ and $A$ saturates $\ker(f)=\ker(g)]$,
\itemit{iv} $f\D g \iff \ker(f)=\ker(g)$ or $[\rank(f)=\rank(g)$ and $A$ saturates both $\ker(f)$ and $\ker(g)]$,
\itemit{v} $f\J g \iff \ker(f)=\ker(g)$ or $|Af|=\rank(f)=\rank(g)=|Ag|$.
\eit
Further, ${\D}={\J}$ in $\TXA$ if and only if $A$ is finite or $A=X$.
\end{thm}

\pf
Green's $\K$ relation in $\TXA$ is the $\Ka$ relation in the principal one-sided ideal $\T_Xa$ of $\T_X$.  

\pfitem{i}  Using Corollary \ref{cor:GreenSa}(i), we have
\[
f \L g \text{ in } \TXA \iff f\La g\text{ in } \T_Xa \iff [f=g\not\in P] \text{ or } [f\L g\text{ in $\T_X$ and } f,g\in P].
\]
Using Theorem~\ref{thm:T}(i) and Proposition \ref{prop:PQP'Q'_T}(i), this is clearly equivalent to the stated conditions.

\pfitem{ii)--(v}  These are treated in similar fashion, using the relevant parts of Theorem~\ref{thm:T}, Corollary \ref{cor:GreenSa} and Proposition \ref{prop:PQP'Q'_T}.

\bigskip\noindent For the final statement, we begin with the backwards implication.  First, if $A=X$, then $\TXA=\T_X$, and so ${\D}={\J}$ in $\TXA$, by Theorem \ref{thm:T}(vi).  Next, suppose $A$ is finite.  Since ${\D}\sub{\J}$ in any semigroup, we just need to prove that ${\J}\sub{\D}$.  To do so, let $(f,g)\in{\J}$.  By part (v), we have $\ker(f)=\ker(g)$ or else $|Af|=\rank(f)=\rank(g)=|Ag|$.  If the former holds, then $(f,g)\in{\D}$, by part (iv), so suppose the latter holds.  Since $f\in\TXA$ with $A$ finite, it follows that $\rank(f)$ is finite; but then it is easy to see that $|Af|=\rank(f)$ is equivalent to $A$ saturating $\ker(f)$.  A similar statement holds for $g$, and it then quickly follows that $(f,g)\in{\D}$, again by (iv).

For the converse, we prove the contrapositive.  Suppose $A$ is infinite and $A\not=X$.  Write $B=X\sm A$, and fix some $x\in A$.  Let $f,g\in\TXA$ be such that: $f$ maps $A$ identically, and all of $B$ onto $x$; $g$ maps $A$ bijectively onto $A\sm\{x\}$, and all of $B$ onto $x$.  Then $\rank(f)=\rank(g)=|A|$, and also $Af=A$ and $Ag=A\sm\{x\}$, so that $|Af|=|Ag|=|A|$; it follows that $(f,g)\in{\J}$ (in $\TXA$), by part (v).  However, since $\ker(f)\not=\ker(g)$, and since $A$ does not saturate~$\ker(g)$, part~(iv) says that $(f,g)\not\in{\D}$ (in $\TXA$).
\epf

\begin{rem}
Parts (i)--(v)  of Theorem \ref{thm:Green_TXA} may be found in \cite[Theorems 3.2, 3.3, 3.6, 3.7 and 3.9]{SS2008}.  The implication $[A$ is finite$]\implies[{\D}={\J}$ in $\TXA]$ is \cite[Theorem 3.12]{SS2008}, but our full characterisation of when ${\D}={\J}$ holds in $\TXA$ appears to be new.
\end{rem}

Here is the corresponding result concerning $\TXal=a\T_X$.  We write $\Delta$ for the trivial relation on $X$: i.e., $\Delta = \set{(x,x)}{x\in X}$.

\begin{thm}\label{thm:Green_TXal}
Let $X$ be a non-empty set, let $\al$ be an equivalence on $X$, and let $f,g\in\TXal$.  Then in the semigroup~$\TXal$,
\bit
\itemit{i} $f\L g \iff \im(f)=\im(g)$,
\itemit{ii} $f\R g \iff f=g$ or $[\ker(f)=\ker(g)$ and $\al$ separates both $\im(f)$ and $\im(g)]$,
\itemit{iii} $f\H g \iff f=g$ or $[\ker(f)=\ker(g)$ and $\al$ separates $\im(f)=\im(g)]$,
\itemit{iv} $f\D g \iff \im(f)=\im(g)$ or $[\rank(f)=\rank(g)$ and $\al$ separates both $\im(f)$ and $\im(g)]$,
\itemit{v} $f\J g \iff \im(f)=\im(g)$ or $\Vert\al f^{-1}\Vert=\rank(f)=\rank(g)=\Vert\al g^{-1}\Vert$.
\eit
Further, ${\D}={\J}$ in $\TXal$ if and only if $\Vert\al\Vert$ is finite or $\al=\Delta$. 
\end{thm}

\pf
Parts (i)--(v) are treated in similar fashion to Theorem \ref{thm:Green_TXA}, as is the backwards implication of the final statement; the details are omitted.  If $\Vert\al\Vert$ is infinite and $\al\not=\Delta$, then we construct a pair $(f,g)\in{\J}\sm{\D}$ as follows.  Let $j\in I$ be such that $|A_j|\geq2$, let $x\in A_j\sm\{a_j\}$ be arbitrary, and let $k\in I\sm\{j\}$.  Then we define $f=a=\binom{A_i}{a_i}_{i\in I}$ and $g=\big( \begin{smallmatrix}A_i&A_k\\a_i&x\end{smallmatrix}\big)_{i\in I\sm\{k\}}$.
\epf

\begin{rem}
Parts (i)--(v) of Theorem \ref{thm:Green_TXal} may be found in \cite[Theorems 2.5, 2.6, 2.7 and 2.10]{MGS2010}.  The implication $[\;\! \Vert\al\Vert$ is finite$]\implies[{\D}={\J}$ in $\TXal]$ is \cite[Corollary 2.13]{MGS2010}, but our full characterisation of when ${\D}={\J}$ holds in $\TXal$ appears to be new.
\end{rem}

\subsection[The regular subsemigroups $\Reg(\TXA)$ and $\Reg(\TXal)$]{\boldmath The regular subsemigroups $\Reg(\TXA)$ and $\Reg(\TXal)$}\label{subsect:Reg_TXA_TXal}

We now concentrate on the regular subsemigroups $P=\Reg(\TXA)$ and $Q=\Reg(\TXal)$; as in Sections~\ref{sect:PLI} and~\ref{sect:PRI}, the results on these involve the local monoid $a\T_Xa=\set{afa}{f\in\T_X}$.  It is well known that $a\T_Xa$ is isomorphic to~$\T_A$.  More specifically, we have the following (see for example \cite[Section~3.3]{Sandwiches2}):

\begin{lemma}\label{lem:TA}
The map $\xi:a\T_Xa\to\T_A:f\mt f|_A$ is an isomorphism. \epfres
\end{lemma}

As a result of Lemma \ref{lem:TA}, instead of utilising the maps
\[
\phi:\Reg(\TXA)=\Reg(\T_Xa)\to a\T_Xa:f\mt af \ANd \psi:\Reg(\TXal)=\Reg(a\T_X)\to a\T_Xa:f\mt fa,
\]
we may compose these with $\xi$, and work with the equivalent surmorphisms
\[
\Phi:\Reg(\TXA)\to\T_A:f\mt (af)|_A=f|_A \AND \Psi:\Reg(\TXal)\to\T_A:f\mt (fa)|_A.
\]
(Note that $(af)|_A=f|_A$ for any $f\in\Reg(\TXA)$ follows from Proposition \ref{prop:PQP'Q'_T}(i).)

Green's relations on $P=\Reg(\TXA)$ and $Q=\Reg(\TXal)$ may easily be described, using Lemmas~\ref{lem:Green_P} and~\ref{lem:Green_Q}, and Corollaries \ref{cor:J=D_P} and \ref{cor:J=D_Q} (and Theorem \ref{thm:T}).  The $\J$-class ordering follows from Lemma \ref{lem:leqJa_leqaJa}(ii) and its dual.

\newpage

\begin{thm}\label{thm:Green_RegTXA}
Let $X$ be a non-empty set, let $A$ be a non-empty subset of $X$, and let $f,g\in P=\Reg(\TXA)$.  Then in the semigroup $P$, 
\bit\bmc2
\itemit{i} $f\L g \iff \im(f)=\im(g)$,
\itemit{ii} $f\R g \iff \ker(f)=\ker(g)$,
\itemit{iii} $f\H g \iff \im(f)=\im(g)$ and $\ker(f)=\ker(g)$,
\itemit{iv} $f\D g \iff f\J g \iff \rank(f)=\rank(g)$.  
\emc\eit
The ${\D}={\J}$-classes of $P$ are the sets 
\[
D_\mu(P) = \set{f\in P}{\rank(f)=\mu} \qquad\text{for each cardinal $1\leq \mu\leq|A|$,}
\]
and they form a chain: $D_\mu(P)\leq D_\nu(P) \iff \mu\leq\nu$.  \epfres
\end{thm}

\begin{thm}\label{thm:Green_RegTXal}
Let $X$ be a non-empty set, let $\al$ be an equivalence relation on $X$, and let $f,g\in Q=\Reg(\TXal)$.  Then in the semigroup $Q$, 
\bit\bmc2
\itemit{i} $f\L g \iff \im(f)=\im(g)$,
\itemit{ii} $f\R g \iff \ker(f)=\ker(g)$,
\itemit{iii} $f\H g \iff \im(f)=\im(g)$ and $\ker(f)=\ker(g)$,
\itemit{iv} $f\D g \iff f\J g \iff \rank(f)=\rank(g)$.  
\emc\eit
The ${\D}={\J}$-classes of $Q$ are the sets 
\[
D_\mu(Q) = \set{f\in Q}{\rank(f)=\mu} \qquad\text{for each cardinal $1\leq \mu\leq\Vert\al\Vert$,}
\]
and they form a chain: $D_\mu(Q)\leq D_\nu(Q) \iff \mu\leq\nu$.  \epfres
\end{thm}

\begin{rem}
Theorem \ref{thm:Green_RegTXA} was originally proved in \cite[Lemma 3]{Sanwong2011}.  The fact that ${\D}={\J}$ on $\Reg(\TXal)$, which is part of Theorem \ref{thm:Green_RegTXal}, was proved in \cite[Theorem 2.9]{MGS2010}; Green's $\R$, $\L$ and $\H$ relations on $\Reg(\TXal)$ were not described in \cite{MGS2010}.
\end{rem}

The Green's structure of $P=\Reg(\TXA)$ and $Q=\Reg(\TXal)$ may be thought of as an ``inflation'' of $a\T_Xa\cong\T_{A}$, in the sense of Remark \ref{rem:inflation_Sa} and its dual; the only additional information required to fully understand the nature of this expansion is the values of the parameters $r$ and $l$ from Theorems \ref{thm:D_structure_P} and \ref{thm:D_structure_Q}, defined in terms of the relations $\Rha$ and $\Lha$, respectively.  To keep notation the same as that of Sections \ref{sect:PRI} and \ref{sect:PLI}, we denote Green's relations on $P$ and~$Q$ by $\Ka$ and $\aK$, respectively.

\begin{lemma}\label{lem:r_T}
Let $f\in P=\Reg(\TXA)$, and write $\mu=\rank(f)$.  Then
\[
|\Rh_f^a/{\Ra}|=\mu^{|X\sm A|}.
\]
\end{lemma}

\pf
An $\Ra$-class $R_g^a$ ($g\in P$) contained in $\Rh_f^a$ is completely determined by the common kernel of each of its elements: i.e., by $\ker(g)$.  If $g\in P$, then
\[
g\in\Rh_f^a \iff \big(g|_A,f|_A\big) = (g\Phi,f\Phi)\in{\R} \text{ in }\T_A \iff {\ker}\big(g|_A\big)={\ker}\big(f|_A\big).
\]
Thus, it suffices to calculate the number of equivalence relations $\ve$ on $X$ such that $\ve=\ker(g)$ for some $g\in P$ and $\ve|_A={\ker}\big(f|_A\big)$.  Now, $\ve|_A={\ker}\big(f|_A\big)$ is a fixed equivalence on $A$ with $\mu$ classes; if we denote these classes by $\set{B_j}{j\in J}$, then the definition of $\ve$ may be completed by assigning each element of $X\sm A$ arbitrarily to any of the $B_j$ (each $\ve=\ker(g)$-class must contain at least one element of $A$, by Proposition~\ref{prop:PQP'Q'_T}(i)).  Since $|J|=\mu$, the result quickly follows.
\epf

\begin{lemma}\label{lem:l_T}
Let $f\in Q=\Reg(\TXal)$, and put $J=\set{i\in I}{\im(f)\cap A_i\not=\emptyset}$.  Then
\[
|{}^a\!\Lh_f/{\aL}|=\prod_{j\in J}|A_j|.
\]
\end{lemma}

\pf
An $\aL$-class ${}^a\!L_g$ ($g\in Q$) contained in ${}^a\!\Lh_f$ is completely determined by the common image of each of its elements: i.e., by $\im(g)$.  
If $g\in Q$, then
\[
g\in{}^a\!\Lh_f \iff \big((ga)|_A,(fa)|_A\big) = (g\Psi,f\Psi)\in{\L} \text{ in }\T_A \iff {\im}\big((ga)|_A\big)={\im}\big((fa)|_A\big).
\]
Thus, it suffices to calculate the number of subsets $B$ of $X$ such that $B=\im(g)$ for some $g\in Q$ and $\set{i\in I}{A_i\cap B\not=\emptyset}=J$; by Proposition~\ref{prop:PQP'Q'_T}(ii), the condition $g\in Q$ forces $|A_j\cap B|=1$ for all $j\in J$.  Such a set~$B$ is determined by choosing an arbitrary element of $A_j$ for each $j\in J$; since these choices can be made in~$\prod_{j\in J}|A_j|$ ways, the result follows.
\epf

\begin{rem}
Lemmas \ref{lem:r_T} and \ref{lem:l_T} respectively give the values of the parameters $r$ and $l$ from Theorems~\ref{thm:D_structure_P} and~\ref{thm:D_structure_Q}.  Thus, the parameter $r$ depends only on $\rank(f)$, meaning that the (vertical) ``stretching'' described in Remark \ref{rem:inflation_Sa} is uniform within $\D$-classes; this can be seen in Figure \ref{fig:RXA}.  In contrast to this, the parameter~$l$ depends not only on $\rank(f)$, but also on the set $J=\set{i\in I}{\im(f)\cap A_i\not=\emptyset}$; as a result, the (horizontal) stetching is not uniform in general, as can be seen in Figure \ref{fig:RXal}.
\end{rem}

We may use Lemmas \ref{lem:r_T} and \ref{lem:l_T} to calculate the sizes of the regular semigroups $P=\Reg(\TXA)$ and $Q=\Reg(\TXal)$.  For an explanation of the notation, see Section \ref{subsect:trans}.

\begin{prop}\label{prop:size_P_T}
Let $X$ be a non-empty set and $A$ a non-empty subset of $X$.  Then the size of the semigroup $P=\Reg(\TXA)$ is given by
\[
|P| = \sum_{\mu=1}^{|A|} \mu!\mu^{|X\sm A|}S(|A|,\mu)\binom{|A|}\mu.
\]
\end{prop}

\pf
Since the ${\J}={\D}$-classes of $P$ are $D_\mu(P)$ for $1\leq\mu\leq|A|$, by Theorem \ref{thm:Green_RegTXA}, we have
\[
|P| = \sum_{\mu=1}^{|A|} |D_\mu(P)|.
\]
Now fix some $1\leq\mu\leq|A|$.  Then $|D_\mu(P)|=\lam\cdot\rho\cdot\eta$, where $\lam=|D_\mu(P)/{\L}|$, $\rho=|D_\mu(P)/{\R}|$, and $\eta$ is the size of any $\H$-class contained in $D_\mu(P)$.  So the proof will be complete if we can show that
\[
\lam=\binom{|A|}\mu \COMMA
\rho=S(|A|,\mu)\mu^{|X\sm A|} \COMMA
\eta=\mu!.
\]
By Remark \ref{rem:inflation_Sa}(iv) and Proposition \ref{prop:combinatorics}(i), we have $\lam=|D_\mu(P)/{\L}|=|D_\mu(\T_A)/{\L}|=\binom{|A|}\mu$.  By Remark~\ref{rem:inflation_Sa}(iv) and Proposition \ref{prop:combinatorics}(ii), $D_\mu(P)$ contains $S(|A|,\mu)$ $\Rha$-classes; by Lemma \ref{lem:r_T}, each of these $\Rha$-classes contains $\mu^{|X\sm A|}$ $\R$-classes; together, these imply that $\rho=S(|A|,\mu)\mu^{|X\sm A|}$.  Now let $f\in D_\mu(P)$ be arbitrary.  By Lemma~\ref{lem:Green_P}(iv), the $\H$-class of $f$ in $P$ is precisely the $\H$-class of $f$ in $\T_X$, which has size $\mu!$ by Proposition \ref{prop:combinatorics}(iv): i.e., $\eta=\mu!$.
\epf

\begin{prop}\label{prop:size_Q_T}
Let $X$ be a non-empty set and $\al$ an equivalence relation on $X$ with equivalence classes $\set{A_i}{i\in I}$.  Then the size of the semigroup $Q=\Reg(\TXal)$ is given by
\[
|Q| = \sum_{\mu=1}^{\Vert\al\Vert}  \mu! S(\Vert\al\Vert,\mu)\sum_{J\sub I\atop |J|=\mu}\prod_{j\in J}|A_j|.
\]
\end{prop}

\pf
This is proved in similar fashion to Proposition \ref{prop:size_P_T}.  We have $|Q| = \sum_{\mu=1}^{\Vert\al\Vert} |D_\mu(Q)|$, and for fixed $1\leq\mu\leq\Vert\al\Vert$, $|D_\mu(Q)|=\lam\cdot\rho\cdot\eta$, where $\lam=|D_\mu(Q)/{\L}|$, $\rho=|D_\mu(Q)/{\R}|$, and $\eta$ is the size of any $\H$-class contained in $D_\mu(Q)$.  This time, we use Remark \ref{rem:inflation_Sa}(iv), Proposition \ref{prop:combinatorics}, and Lemmas \ref{lem:Green_Q} and \ref{lem:l_T} to show that
\[
\lam = \sum_{J\sub I\atop |J|=\mu}\prod_{j\in J}|A_j| \COMMA
\rho = S(\Vert\al\Vert,\mu) \COMMA
\eta = \mu!.
\]
(For the value of $\lam$, note that the $\aLh$-classes in $D_\mu(P)$ are in one-one correspondence with the $\L$-classes in $D_\mu(\T_A)$, which are indexed by the subsets of $A=\im(a)$ of size $\mu$, and hence by the subsets of $I$ of size $\mu$; the number of $\aL$-classes contained in an $\aLh$-class induced by a given subset $J\sub I$ is given in Lemma \ref{lem:l_T}.)
\epf

\begin{rem}
If $A=X$ or $\al=\Delta$, then $P=\TXA=\T_X$ and $Q=\TXal=\T_X$, and Propositions~\ref{prop:size_P_T} and \ref{prop:size_Q_T} both reduce to the well-known formulae $|\T_X|=\sum_{\mu=1}^{|X|}\mu!S(|X|,\mu)\binom{|X|}\mu$.  (Of course we also have~${|\T_X|=|X|^{|X|}}$.)
\end{rem}

In the case of infinite $X$, the expressions for $|P|$ and $|Q|$ in Propositions \ref{prop:size_P_T} and \ref{prop:size_Q_T} simplify significantly:

\begin{cor}\label{cor:size_P_T}
Let $X$ be an infinite set and $A$ a non-empty subset of $X$.  Then the size of the semigroup $P=\Reg(\TXA)$ is given by
\[
|P|=\begin{cases}
1 &\text{if $|A|=1$}\\
2^{|X|} &\text{if $|A|\geq2$.}
\end{cases}
\]
\end{cor}

\pf
The statement for $|A|=1$ being clear, suppose $|A|\geq2$.  Since $|P|\leq|\T_X|=2^{|X|}$, it suffices to show that $|P|\geq2^{|X|}$.  To do so, we show that the $\mu=2$ term of the sum in Proposition \ref{prop:size_P_T} is at least $2^{|X|}$.  We denote this term by $\xi$.  
First, if $|A|<|X|$, then $|X\sm A|=|X|$, and we have $\xi\geq2^{|X\sm A|}=2^{|X|}$.  
On the other hand, if $|A|=|X|$, then $\xi\geq S(|A|,2)=S(|X|,2)=2^{|X|}$.
\epf

\begin{cor}\label{cor:size_Q_T}
Let $X$ be an infinite set and $\al$ an equivalence relation on $X$ with equivalence classes $\set{A_i}{i\in I}$.  Then the size of the semigroup $Q=\Reg(\TXal)$ is given by
\[
|Q| = 2^{\Vert\al\Vert}\prod_{i\in I}|A_i|.
\]
\end{cor}

\pf
For simplicity, we will write $\pi=\prod_{i\in I}|A_i|$ throughout the proof.  
If $\Vert\al\Vert=1$, then $P=\TXal$ consists of all constant mappings, of which there are $|X|$; but we also note that $2^{\Vert\al\Vert}\pi$ simplifies to $|X|$ in this case (here, $X$ is the only equivalence class).  
For the rest of the proof, we assume that $\Vert\al\Vert\geq2$.  For $1\leq\mu\leq\Vert\al\Vert$, we denote by $\xi_\mu$ the $\mu$th term of the sum in Proposition \ref{prop:size_Q_T}.  We now consider separate cases according to whether $\Vert\al\Vert$ is finite or infinite.

Suppose first that $\Vert\al\Vert$ is finite.  Since $X$ is infinite and $|X|=\sum_{i\in I}|A_i|$, at least one of the $A_i$ is infinite, and hence $\pi=\prod_{i\in I}|A_i|$ is infinite.  For any $1\leq\mu\leq\Vert\al\Vert$, 
\[
\xi_\mu = \mu! S(\Vert\al\Vert,\mu)\sum_{J\sub I\atop |J|=\mu}\prod_{j\in J}|A_j|
\leq \mu! S(\Vert\al\Vert,\mu)\sum_{J\sub I}\prod_{i\in I}|A_i|
= \mu! S(\Vert\al\Vert,\mu)2^{|I|}\pi = \pi,
\]
with the last equality holding because $\mu! S(\Vert\al\Vert,\mu)2^{|I|}$ is finite and $\pi$ infinite.  Since $\Vert\al\Vert$ is finite, it follows that $|Q|=\sum_{\mu=1}^{\Vert\al\Vert}\xi_\mu\leq\Vert\al\Vert\pi=\pi=2^{\Vert\al\Vert}\pi$.  For the reverse inequality, we have
\[
|Q|\geq\xi_{\Vert\al\Vert} = \Vert\al\Vert! S(\Vert\al\Vert,\Vert\al\Vert)\sum_{J\sub I\atop |J|=\Vert\al\Vert}\prod_{j\in J}|A_j| = \Vert\al\Vert! \prod_{i\in I}|A_i| = \Vert\al\Vert!\pi = \pi = 2^{\Vert\al\Vert}\pi,
\]
again because $\Vert\al\Vert!$ and $2^{\Vert\al\Vert}$ are finite, and $\pi$ infinite.

Now suppose $\Vert\al\Vert$ is infinite.  For any $1\leq\mu\leq\Vert\al\Vert$, 
\[
\xi_\mu = \mu! S(\Vert\al\Vert,\mu)\sum_{J\sub I\atop |J|=\mu}\prod_{j\in J}|A_j|
\leq \Vert\al\Vert! 2^{\Vert\al\Vert}\sum_{J\sub I}\prod_{i\in I}|A_i|
= 2^{\Vert\al\Vert}\cdot2^{\Vert\al\Vert}\cdot2^{|I|}\pi = 2^{\Vert\al\Vert}\pi.
\]
Since there are fewer than $2^{\Vert\al\Vert}$ terms in the sum in Proposition \ref{prop:size_Q_T}, it follows that $|Q|\leq2^{\Vert\al\Vert}\cdot2^{\Vert\al\Vert}\pi=2^{\Vert\al\Vert}\pi$.
But also
\[
|Q|\geq\xi_{\Vert\al\Vert} = \Vert\al\Vert! S(\Vert\al\Vert,\Vert\al\Vert)\sum_{J\sub I\atop|J|=|I|}\prod_{j\in J}|A_j| 
\geq \Vert\al\Vert! \prod_{i\in I}|A_i| 
= 2^{\Vert\al\Vert}\pi,
\]
completing the proof.
\epf

\begin{rem}
As observed in the above proof, we have $|Q|=|X|$ if $\Vert\al\Vert=1$ and, more generally, ${|Q|=\prod_{i\in I}|A_i|}$ if $\Vert\al\Vert$ is finite.  In fact, it then follows from $|X|=\sum_{i\in I}|A_i|=\max_{i\in I}|A_i|=\prod_{i\in I}|A_i|$ that $|Q|=|X|$ for finite $\Vert\al\Vert$.  On the other hand, if $\Vert\al\Vert$ is infinite, then $|Q|\geq2^{\Vert\al\Vert}$ is always uncountable, and can be as large as $2^{|X|}$.
\end{rem}

\begin{rem}
If $A=X$ or $\al=\Delta$, then Propositions~\ref{prop:size_P_T} and \ref{prop:size_Q_T} reduce to $|\T_X|=2^{|X|}$ (for infinite $X$).
\end{rem}

We may also calculate the ranks of $P=\Reg(\TXA)$ and $Q=\Reg(\TXal)$.  For this, we first show that the semigroups $P$ and $Q$ are RI- and LI-dominated, respectively, regardless of the values of $|A|$ and~$\Vert\al\Vert$.

\begin{prop}\label{prop:RI_P_T}
Let $X$ be a non-empty set and $A$ a non-empty subset of $X$.  Then the semigroup ${P=\Reg(\TXA)}$ is RI-dominated.
\end{prop}

\pf
Let $f=\binom{F_j}{f_j}\in P$ be arbitrary.  Since $f\in P$, Proposition \ref{prop:PQP'Q'_T}(i) says that $A\cap F_j\not=\emptyset$ for all $j\in J$.  For each $j$, let $I_j=\set{i\in I}{a_i\in F_j}$, and fix a partition $F_j=\bigsqcup_{i\in I_j}F_{j,i}$ so that $a_i\in F_{j,i}$ for each $i\in I_j$.  Put $b=\binom{F_{j,i}}{a_i}_{j\in J,\ i\in I_j}$.  Proposition \ref{prop:PQP'Q'_T}(i) immediately gives $b\in P$, as $A$ is a cross-section of $\ker(b)$.  Since $b$ maps~$A$ identically, we have $a=ab$, and it follows that $b$ is a right identity for $P$ (since $a$ is).  Finally, it is clear that $f=bf$, so that $f\leqR b$.
\epf

\begin{prop}\label{prop:LI_Q_T}
Let $X$ be a non-empty set and $\al$ an equivalence relation on $X$.  Then the semigroup $Q=\Reg(\TXal)$ is LI-dominated.
\end{prop}

\pf
Let $f=\binom{F_j}{f_j}\in Q$ be arbitrary.  For each $j\in J$, we have $f_j\in A_{i_j}$ for some $i_j\in I$.  Since $f\in Q$, Proposition \ref{prop:PQP'Q'_T}(ii) says that the map $j\mt i_j$ is injective.  Write $K=\set{i_j}{j\in J}$, and define $b=\big(\begin{smallmatrix}A_{i_j}&A_l\\f_j&a_l\end{smallmatrix}\big)_{j\in J,\ l\in I\sm K}$.  Then again one may show that $b\in Q$ is a left identity for $Q$ and that $f=fb\leqL b$.
\epf

\begin{thm}\label{thm:rank_P_T}
Let $X$ be a non-empty set and $A$ a non-empty subset of $X$.  Then the rank of the semigroup $P=\Reg(\TXA)$ is given by
\[
\rank(P) = 
\begin{cases}
1 &\text{if $|A|=1$}\\
2^{|X|} &\text{if $|A|\geq2$ and $X$ is infinite}\\
3 &\text{if $3\leq|A|=|X|$ is finite}\\
1+|A|^{|X\sm A|} &\text{otherwise.}
\end{cases}
\]
\end{thm}

\pf
If $|A|=1$ then $|P|=1$ and the result is clear, so we assume $|A|\geq2$ for the rest of the proof.  If $X$ is infinite, then by Corollary \ref{cor:size_P_T}, $|P|=2^{|X|}$ is uncountable, and so $\rank(P)=|P|$, completing the proof in this case.  

For the rest of the proof we assume $X$ is finite (and $|A|\geq2$).  It follows that $A$ is finite as well, and so $\T_A\sm\S_A$ is an ideal of $\T_A$.  Given the isomorphism $\xi:a\T_Xa\to\T_A$ from Lemma~\ref{lem:TA}, it follows that $a\T_Xa\sm G_{a\T_Xa}$ is an ideal of $a\T_Xa$.  Combining this with Proposition \ref{prop:RI_P_T}, it follows that Theorem \ref{thm:rank_P} applies, and it gives
\begin{equation}\label{eq:rank_P_T}
\rank(P) = \relrank{\T_A}{\S_A} + \max\left(|A|^{|X\sm A|},\rank(\S_A)\right).
\end{equation}
If $A=X$, then $P=\TXA=\T_X$, and so $\rank(P)=\rank(\T_X)$ in this case.  It is well known that $\rank(\T_X)=2$ if $|X|=2$ and $\rank(\T_X)=3$ for finite $|X|\geq3$, agreeing with the claimed values for $\rank(P)$.

Finally, suppose $2\leq|A|<|X|$.  Then $\relrank{\T_A}{\S_A}=1$; see for example \cite[Proposition~1.2]{HRH1998}.  Also, $\rank(\S_A)\leq2$ (it can only be $1$ if $|A|=2$).  Since $2\leq|A|<|X|$, we have $|A|^{|X\sm A|}\geq2$, and so $\max\left(|A|^{|X\sm A|},\rank(\S_A)\right)=|A|^{|X\sm A|}$.  By \eqref{eq:rank_P_T}, this completes the proof.
\epf

\begin{rem}
The finite case of Theorem \ref{thm:rank_P_T} was proved in \cite[Theorem 3.6]{SS2013}.  Alternative proofs of Theorems \ref{thm:rank_P_T} and \ref{thm:rank_Q_T} may be found in \cite{Sandwiches2}.
\end{rem}

Recall that $\Delta$ denotes the trivial relation on $X$; we also write $\nabla=X\times X$ for the universal relation.

\begin{thm}\label{thm:rank_Q_T}
Let $X$ be a non-empty set and $\al$ an equivalence relation on $X$ with equivalence classes $\set{A_i}{i\in I}$.  Then the rank of the semigroup $Q=\Reg(\TXal)$ is given by
\[
\rank(Q) = 
\begin{cases}
|X| &\text{if $\al=\nabla$}\\
2^{\Vert\al\Vert}\prod_{i\in I}|A_i| &\text{if $\Vert\al\Vert$ is infinite}\\
3 &\text{if $\al=\Delta$ and $|X|\geq3$ is finite}\\
1+\prod_{i\in I}|A_i| &\text{otherwise.}
\end{cases}
\]
\end{thm}

\pf
If $\al=\nabla$ then $Q$ is the right-zero band of all constant mappings, and hence $\rank(Q)=|Q|=|X|$.  
If $\Vert\al\Vert$ is infinite, then by Corollary \ref{cor:size_Q_T}, $|Q|=2^{\Vert\al\Vert}\prod_{i\in I}|A_i|$ is uncountable, so again $\rank(Q)=|Q|$.

For the rest of the proof we assume that $\Vert\al\Vert$ is finite, and that $\al\not=\nabla$.  It follows that $\rank(a)=\Vert\al\Vert$ is finite, so as in the proof of Theorem \ref{thm:rank_P_T}, it follows from Theorem \ref{thm:rank_Q}, Lemma \ref{lem:l_T} and Proposition \ref{prop:LI_Q_T} that
\begin{equation}\label{eq:rank_Q_T}
\rank(Q) = \relrank{\T_A}{\S_A} + \max\left(\pi,\rank(\S_A)\right),
\end{equation}
where again we have written $\pi=\prod_{i\in I}|A_i|$.  If $\al=\Delta$ then $\pi=1$, so it follows from \eqref{eq:rank_Q_T} and Lemma \ref{lem:rankWT} that
\[
\rank(Q) = \relrank{\T_A}{\S_A} + \rank(\S_A) = \rank(\T_A).
\]
Consulting Theorem \ref{thm:IGT}, this agrees with the claimed value(s).  If $\al\not=\Delta$, then $\pi\geq2\geq\rank(\S_A)$.  Since $\al\not=\nabla$, $|A|=\Vert\al\Vert\geq2$, so $\relrank{\T_A}{\S_A}=1$, and it follows from~\eqref{eq:rank_Q_T} that $\rank(Q)=1+\pi$.
\epf

\subsection[The idempotent-generated subsemigroups $\bbE(\TXA)$ and $\bbE(\TXal)$]{\boldmath The idempotent-generated subsemigroups $\bbE(\TXA)$ and $\bbE(\TXal)$}\label{subsect:IG_TXA_TXal}

In this section, we study the idempotent-generated subsemigroups of the principal one-sided ideals $\T_Xa=\TXA$ and $a\T_X=\TXal$.  In the literature on the semigroups $\TXA$ and $\TXal$, these subsemigroups seem not to have been explicitly investigated.  
Theorems \ref{thm:E_Sa} and \ref{thm:E_aS} (and the isomorphism ${\xi:a\T_Xa\to\T_A}$ from Lemma \ref{lem:TA}) yield immediate descriptions of these subsemigroups in terms of the corresponding idempotent-generated subsemigroup of $\T_A$, which itself was described in \cite{Howie1966}.

\begin{thm}\label{thm:IGTA}
Let $X$ be a non-empty set, let $A$ be a non-empty subset of $X$, and let $\al$ be an equivalence relation on $X$.  Then
\bit
\itemit{i} $\bbE(\TXA)=\set{f\in\TXA}{f|_A\in\bbE(\T_A)}$,
\itemit{ii} $\bbE(\TXal)=\set{f\in\TXal}{(fa)|_A\in\bbE(\T_A)}$.  \epfres
\eit
\end{thm}

In the case that $|A|$ or $\Vert\al\Vert$ is finite, Theorem \ref{thm:IGTA} takes on a particularly elegant form (regardless of whether $X$ is itself finite or infinite).  Before we state it, it will be convenient to describe the one-sided identities of $\TXA$ and $\TXal$.

\begin{lemma}\label{lem:RI_LI_T}
Let $X$ be a non-empty set, let $A$ be a non-empty subset of $X$, and let $\al$ be an equivalence relation on $X$.  Then
\bit
\itemit{i} $\RI(\TXA)=\set{f\in\TXA}{xf=x\ (\forall x\in A)}$,
\itemit{ii} $\LI(\TXal)=\set{f\in\TXal}{(xf,x)\in\al\ (\forall x\in X)}$.
\eit
\end{lemma}

\pf
We just prove (i), as (ii) is similar.  An element $f\in\TXA$ is a right identity for $\TXA$ if and only if $a=af$ (since $a$ is a right identity); it is easy to see that this is equivalent to the stated condition.
\epf

\begin{thm}\label{thm:IGTal}
Let $X$ be a non-empty set, let $A$ be a non-empty finite subset of $X$, and let $\al$ be an equivalence relation on $X$ with finitely many equivalence classes.  Then
\bit
\itemit{i} $\bbE(\TXA)=\bigset{f\in\TXA}{xf=x\ (\forall x\in A)}\cup\bigset{f\in\TXA}{\rank(f)<|A|}$,
\itemit{ii} $\bbE(\TXal)=\bigset{f\in\TXal}{(x,xf)\in\al \ (\forall x\in X)}\cup\bigset{f\in\TXal}{\rank(f)<\Vert\al\Vert}$.
\eit
\end{thm}

\pf
These follow quickly from Propositions \ref{prop:singular_ESa} and \ref{prop:singular_EaS}, together with Theorem \ref{thm:IGT} and Lemma~\ref{lem:RI_LI_T}, and the $a\phi^{-1}=\RI(Sa)$ and $a\psi^{-1}=\LI(aS)$ parts of Propositions \ref{prop:MI_P} and \ref{prop:MI_Q}.
\epf

Now that we have described the elements of $\bbE(\TXA)$ and $\bbE(\TXal)$, we wish to calculate the ranks and idempotent ranks of these semigroups.  First, we count the idempotents.

\begin{prop}\label{prop:E_TXA_TXal}
Let $X$ be a non-empty set, let $A$ be a non-empty subset of $X$, and let $\al$ be an equivalence relation on $X$ with equivalence classes $\set{A_i}{i\in I}$.  Then
\bit
\itemit{i} $|E(\TXA)|=\begin{cases}
1 &\hspace{10.8mm}\text{if $|A|=1$}\\[2mm]
2^{|X|} &\hspace{10.8mm}\text{if $X$ is infinite and $|A|\geq2$}\\[2mm]
\displaystyle\sum_{\mu=1}^{|A|}\mu^{|X|-\mu}\binom{|A|}\mu &\hspace{10.8mm}\text{otherwise,}
\end{cases}$
\itemit{ii} $|E(\TXal)|=\begin{cases}
\displaystyle2^{\Vert\al\Vert}\prod_{i\in I}|A_i| &\text{if $X$ is infinite}\\[2mm]
\displaystyle\sum_{\mu=1}^{\Vert\al\Vert}\mu^{\Vert\al\Vert-\mu}\sum_{J\sub I\atop |J|=\mu}\prod_{j\in J}|A_j| &\text{if $X$ is finite.}
\end{cases}$
\eit
\end{prop}

\pf
(i).  Again the $|A|=1$ case is trivial, so we assume $|A|\geq2$.  

Suppose first that $X$ is infinite.  Since $|E(\TXA)|=|E(P)|\leq|P|=2^{|X|}$, by Corollary \ref{cor:size_P_T}, it suffices to show that $|E(\TXA)|\geq2^{|X|}$.  Since $|A|\geq2$, we may fix distinct $x,y\in A$.  Then for any partition $X\sm\{x,y\}=B\sqcup C$, the map $\big(\begin{smallmatrix}B\cup\{x\}&C\cup\{y\}\\x&y\end{smallmatrix}\big)$ belongs to $E(\TXA)$.  Since there are $2^{|X|}$ such partitions, the result follows.

Now suppose $X$ is finite.  An idempotent $f\in E(\TXA)$ may be specified by:
\bit
\item choosing $\mu=\rank(f)$, which can be anything from $1$ to $|A|$,
\item choosing $\im(f)$, which must be a subset of $A$ of size $\mu$, 
\item choosing $xf$ for each $x\in X\sm \im(f)$ (note that $f$ must map the elements of $\im(f)$ identically).
\eit
Since there are $\binom{|A|}\mu$ choices for $\im(f)$, and $\mu^{|X\sm\im(f)|}=\mu^{|X|-\mu}$ choices for the $xf$ ($x\in X\sm\im(f)$), the stated formula follows.

\pfitem{ii} Again, for simplicity, we will write $\pi=\prod_{i\in I}|A_i|$.  Suppose first that $X$ is infinite.  As in the previous case, by Corollary \ref{cor:size_Q_T}, it suffices to show that $|E(\TXal)|\geq2^{\Vert\al\Vert}\pi$.  Since $X$ is infinite, at least one of~$\Vert\al\Vert$ or $\pi$ must be infinite.  It follows that $2^{\Vert\al\Vert}\pi=\max(2^{\Vert\al\Vert},\pi)$, so it suffices to show that 
\bit\bmc2
\itemnit{a} $|E(\TXal)|\geq\pi$, and
\itemnit{b} $|E(\TXal)|\geq2^{\Vert\al\Vert}$.
\emc\eit
First, note that for any choice function $I\to X:i\mt b_i$ with $b_i\in A_i$ for each $i$, the map $\binom{A_i}{b_i}$ is an idempotent of $\TXal$; since there are $\pi$ such choice functions, this gives (a).  To prove (b), note first that if $\Vert\al\Vert$ is finite, then $\pi$ must be infinite (as noted above), and so (a) gives $|E(\TXal)|\geq\pi\geq2^{\Vert\al\Vert}$.  Now suppose $|I|=\Vert\al\Vert$ is infinite.  Fix some distinct $j,k\in I$.  Then for any partition $I\sm\{j,k\}=M\sqcup N$, the map
\[
\left(\begin{matrix}A_j\cup\bigcup_{m\in M}A_m & A_k\cup\bigcup_{n\in N}A_n\\a_j&a_k\end{matrix}\right)
\]
is an idempotent of $\TXal$.  Since there are $2^{|I|}=2^{\Vert\al\Vert}$ such partitions, this completes the proof of (b).  As noted above, this completes the proof of (ii) in the case of infinite $X$.

Now suppose $X$ is finite.  An idempotent $f\in E(\TXal)$ may be specified by:
\bit
\item choosing $\mu=\rank(f)$, which can be anything from $1$ to $\Vert\al\Vert$,
\item choosing $\im(f)$, which must be of the form $\set{b_j}{j\in J}$ for some subset $J\sub I$ of size $\mu$, and where $b_j\in A_j$ for each $j$,
\item choosing $A_kf$ for each $k\in I\sm J$ (note that $A_jf=b_j$ for each $j$).
\eit
There are $\sum_{J\sub I, |J|=\mu}\prod_{j\in J}|A_j|$ ways to perform the second task, and $\mu^{\Vert\al\Vert-\mu}$ to do the third.
\epf

\begin{thm}\label{thm:E_TXA_TXal}
Let $X$ be a non-empty set, let $A$ be a non-empty subset of $X$, and let $\al$ be an equivalence relation on $X$ with equivalence classes $\set{A_i}{i\in I}$.  Then
\bit
\itemit{i} $\rank(\bbE(\TXA))=\idrank(\bbE(\TXA))=\begin{cases}
1 &\text{if $|A|=1$}\\[2mm]
2^{|X|} &\text{if $X$ is infinite and $|A|\geq2$}\\[2mm]
2+2^{|X|-2} &\text{if $|A|=2$ and $X$ is finite}\\[2mm]
\binom{|A|}2+|A|^{|X|-|A|} &\text{otherwise,}
\end{cases}$
\itemit{ii} $\rank(\bbE(\TXal))=\idrank(\bbE(\TXal))=\begin{cases}
|X| &\text{if $\Vert\al\Vert=1$}\\[2mm]
2^{\Vert\al\Vert}\prod_{i\in I}|A_i| &\hspace{0.4mm}\text{if $\Vert\al\Vert$ is infinite}\\[2mm]
2+\prod_{i\in I}|A_i| &\hspace{0.4mm}\text{if $\Vert\al\Vert=2$ and $X$ is finite}\\[2mm]
\binom{\Vert\al\Vert}2+\prod_{i\in I}|A_i| &\hspace{0.4mm}\text{otherwise.}
\end{cases}$
\eit
\end{thm}

\pf
(i).  Again, the $|A|=1$ case is trivial, so we assume that $|A|\geq2$.  

Next suppose $A$ is infinite.  Then so too is $X$, so Proposition \ref{prop:E_TXA_TXal}(i) gives
\[
2^{|X|} = |E(\TXA)| \leq |\bbE(\TXA)| \leq |\T_X| = 2^{|X|}.
\]
It follows that $|\bbE(\TXA)|=2^{|X|}$ is uncountable, and so
\[
\rank(\bbE(\TXA))=\idrank(\bbE(\TXA))=|\bbE(\TXA)| = 2^{|X|}.
\]

Now suppose $A$ is finite.  By Proposition \ref{prop:RI_P_T}, $P=\Reg(\TXA)$ is RI-dominated; Theorem \ref{thm:rank_EP} and Lemma \ref{lem:r_T} then give
\[
\rank(\bbE(\TXA)) = \rank(\bbE(\T_A)) + |A|^{|X|-|A|} - 1
\anD
\idrank(\bbE(\TXA)) = \idrank(\bbE(\T_A)) + |A|^{|X|-|A|} - 1.
\]
Theorem \ref{thm:IGT} completes the proof.

\pfitem{ii}  If $\Vert\al\Vert=1$, then $\bbE(\TXal)=\TXal$ consists entirely of all the constant mappings, and the stated formula again follows quickly.  
All other cases are treated in almost identical fashion to part, treating the cases of $\Vert\al\Vert$ finite and infinite separately.
\epf

\subsection{Egg-box diagrams}\label{subsect:eggbox}

\begin{figure}[ht]
\begin{center}
\includegraphics[height=2.73cm]{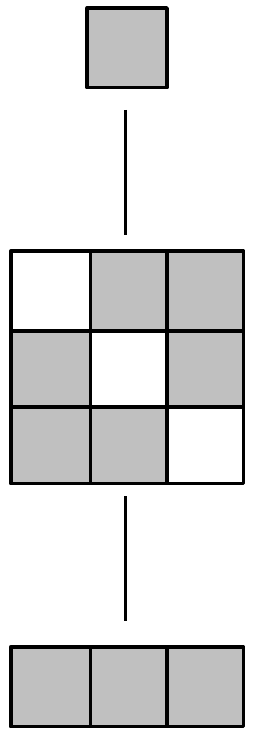} 
\qquad\qquad\qquad
\includegraphics[height=6cm]{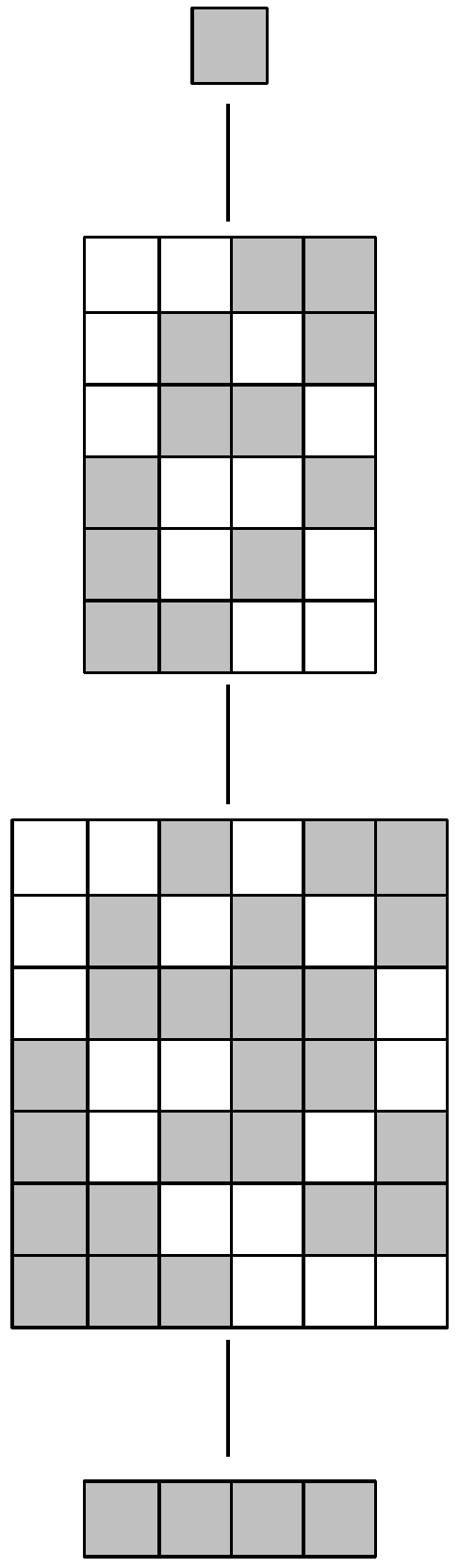} 
\qquad\qquad\qquad
\includegraphics[height=10cm]{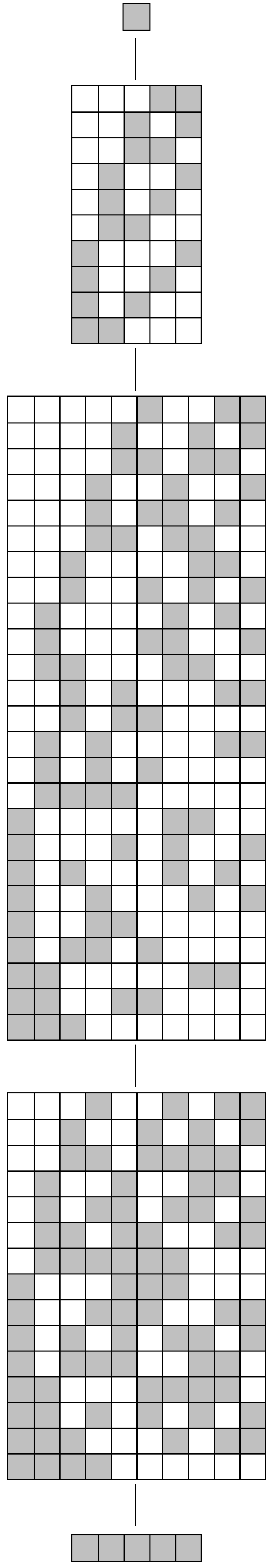} 
\caption[blah]{Left to right: egg-box diagrams of $\T_X$, where $|X|=3$, $4$ and $5$.}
\label{fig:TX}
\end{center}
\end{figure}

\begin{figure}[ht]
\begin{center}
\includegraphics[width=\textwidth]{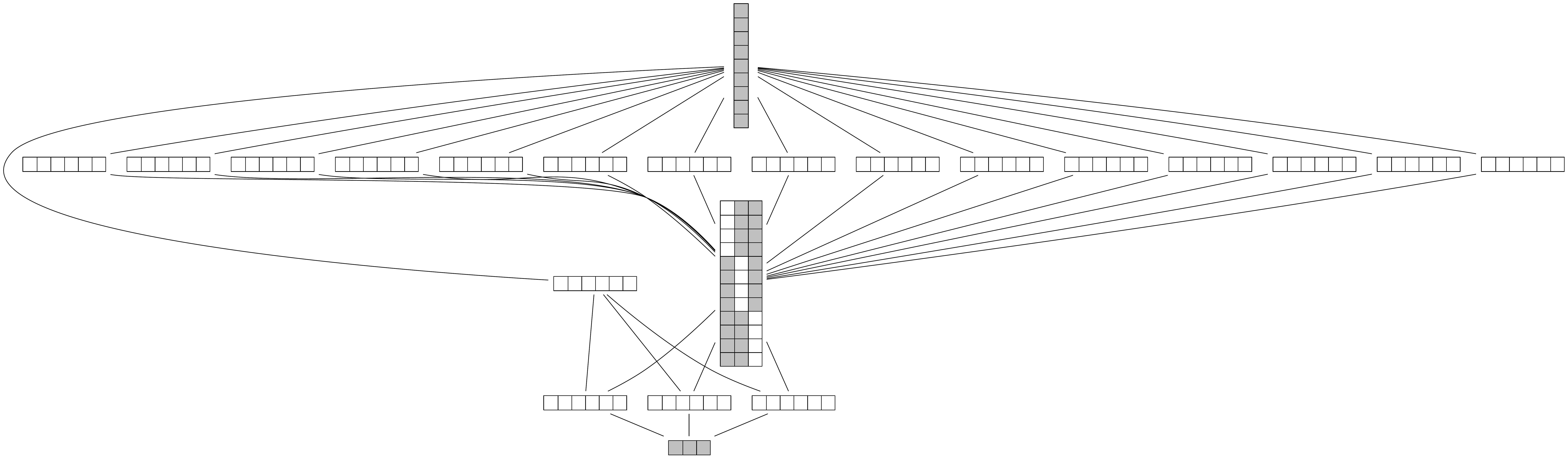} 
\\[5mm]
\includegraphics[width=\textwidth]{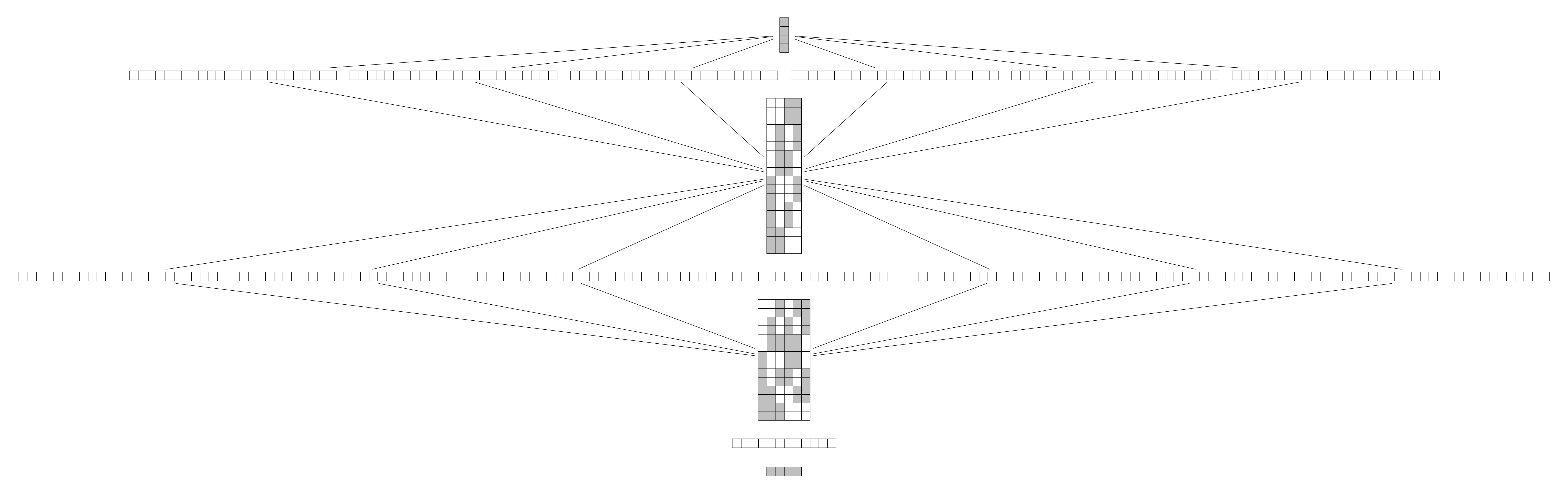} 
\caption[blah]{Egg-box diagrams of $\TXA$, where $|X|=5$ and $|A|=3$ (top) or $|A|=4$ (bottom).}
\label{fig:TXA}
\end{center}
\end{figure}

\begin{figure}[ht]
\begin{center}
\includegraphics[height=6cm]{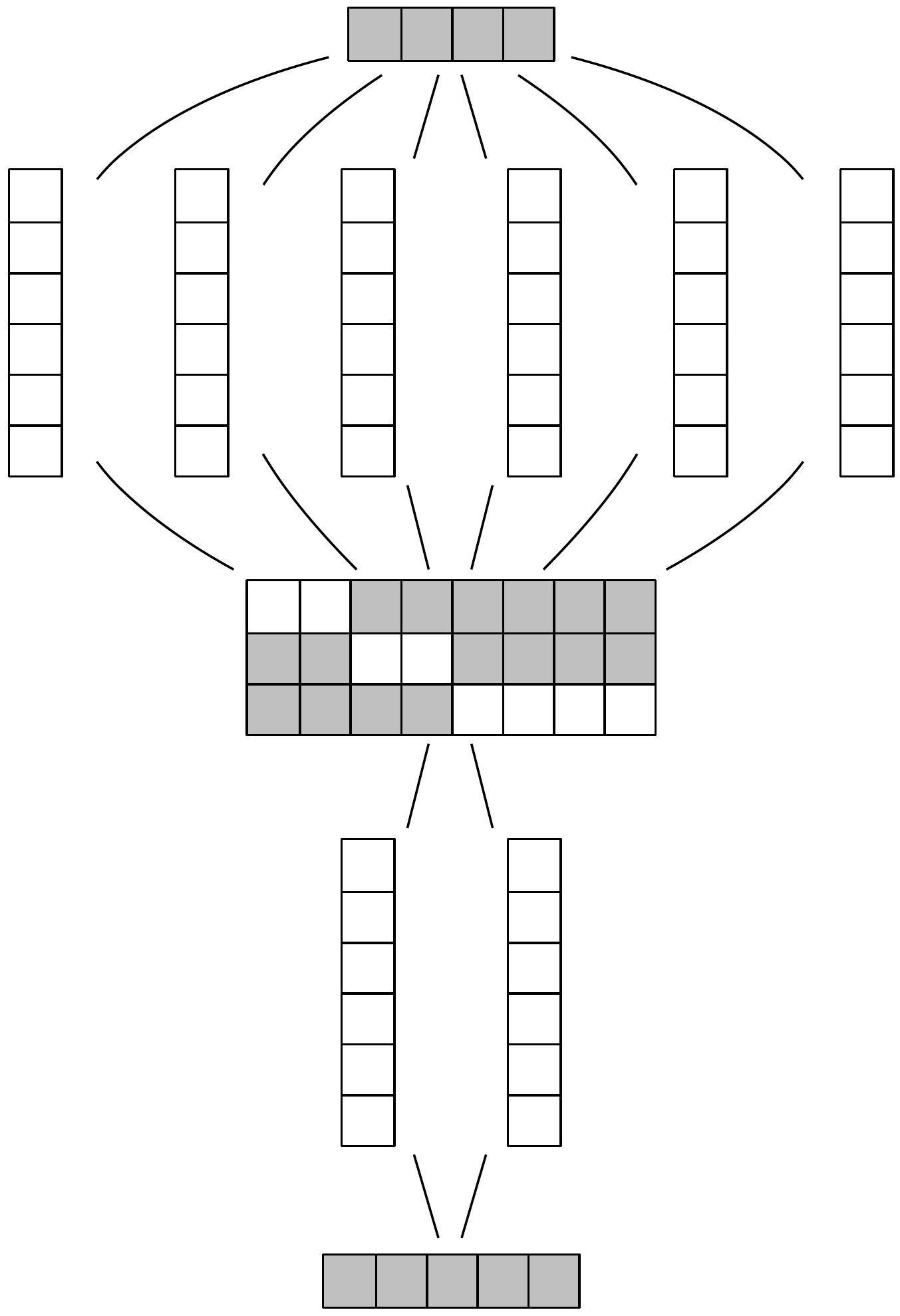} 
\qquad\qquad\qquad
\includegraphics[height=6cm]{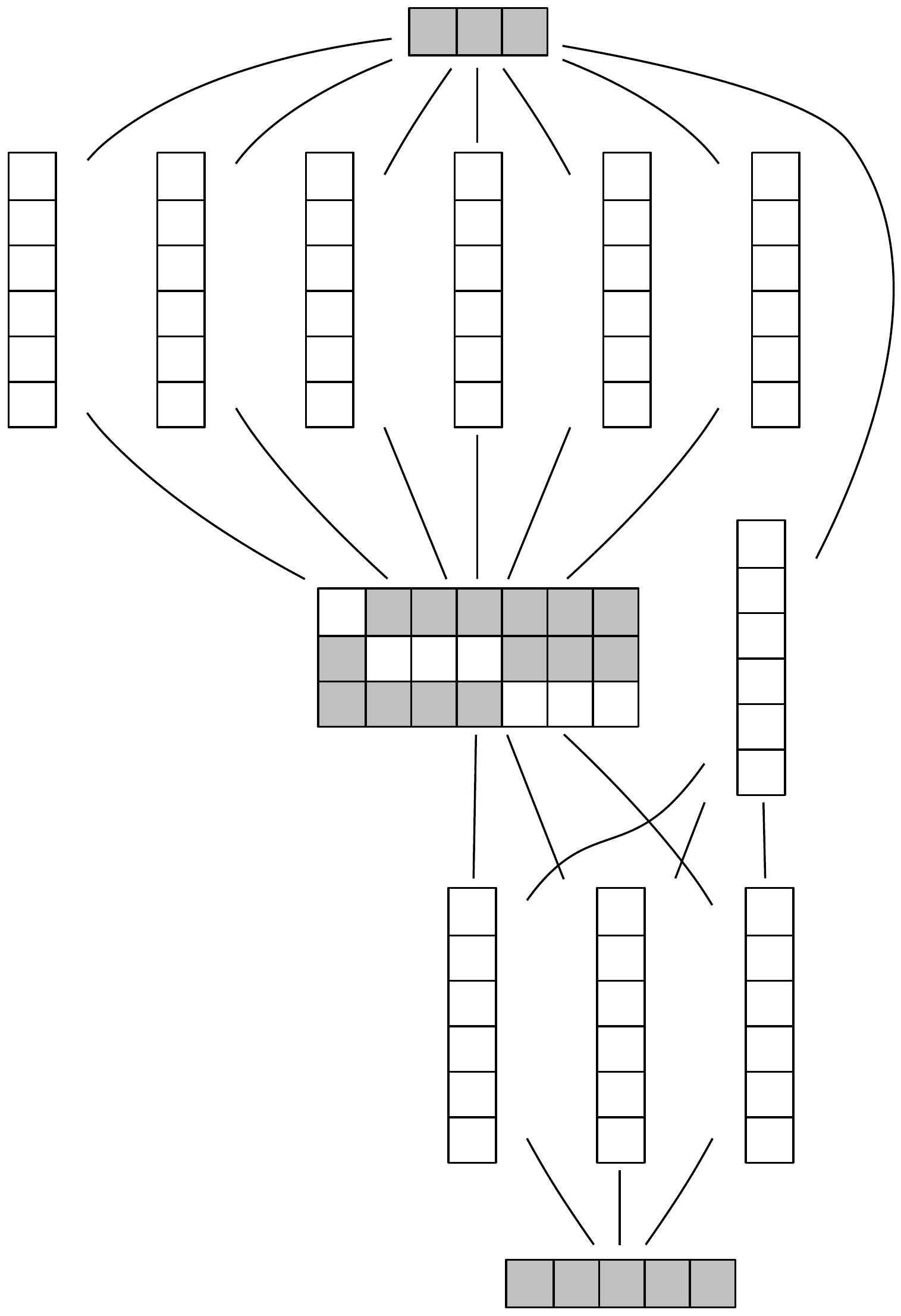} 
\caption[blah]{Egg-box diagrams of $\TXal$, where $|X|=\{1,2,3,4,5\}$ and $X/\al=\{\{1\},\{2,3\},\{4,5\}\}$ (left) or $X/\al=\{\{1\},\{2\},\{3,4,5\}\}$ (right).}
\label{fig:TXal}
\end{center}
\end{figure}

\begin{figure}[ht]
\begin{center}
\includegraphics[height=5.6cm]{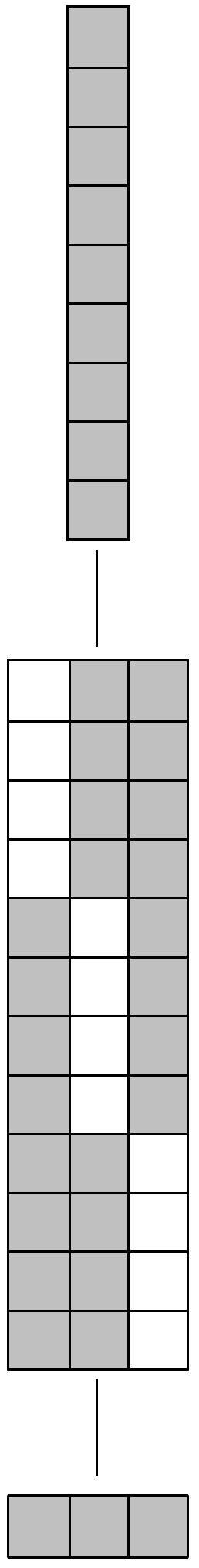} 
\qquad\qquad\qquad
\includegraphics[height=8cm]{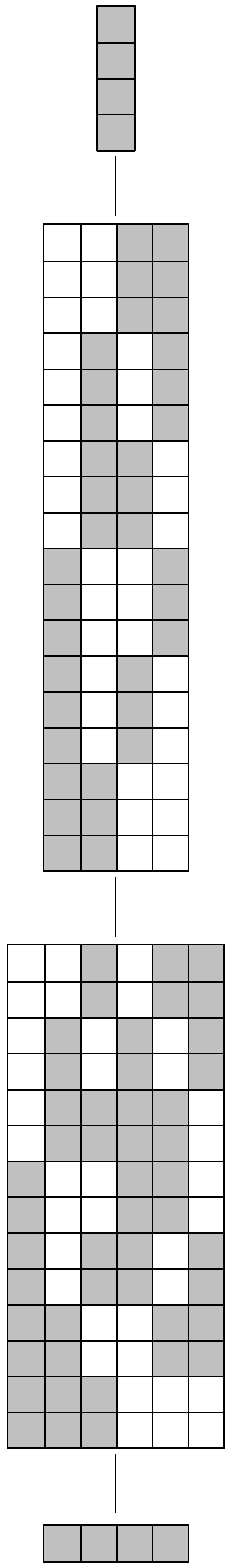} 
\caption[blah]{Egg-box diagrams of $\Reg(\TXA)$, where $|X|=5$ and $|A|=3$ (left) or $|A|=4$ (right).}
\label{fig:RXA}
\end{center}
\end{figure}

\begin{figure}[ht]
\begin{center}
\includegraphics[height=4cm]{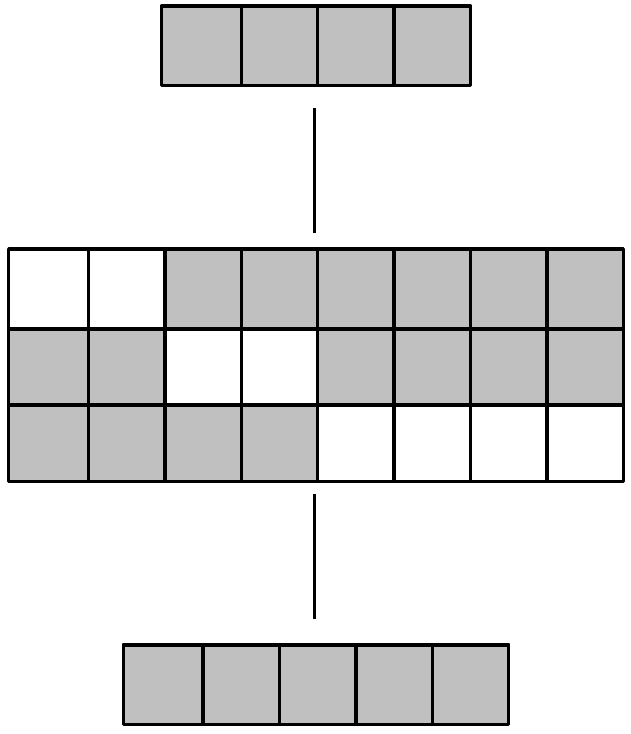} 
\qquad\qquad\qquad
\includegraphics[height=4cm]{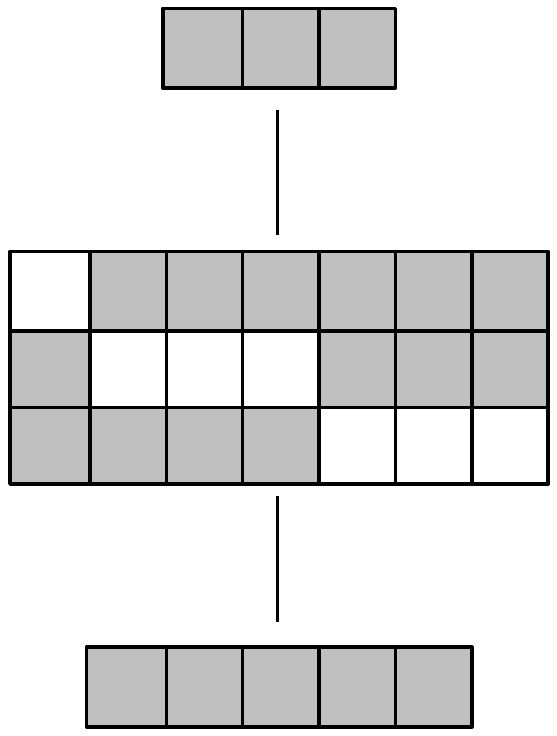} 
\qquad\qquad\qquad
\includegraphics[height=6cm]{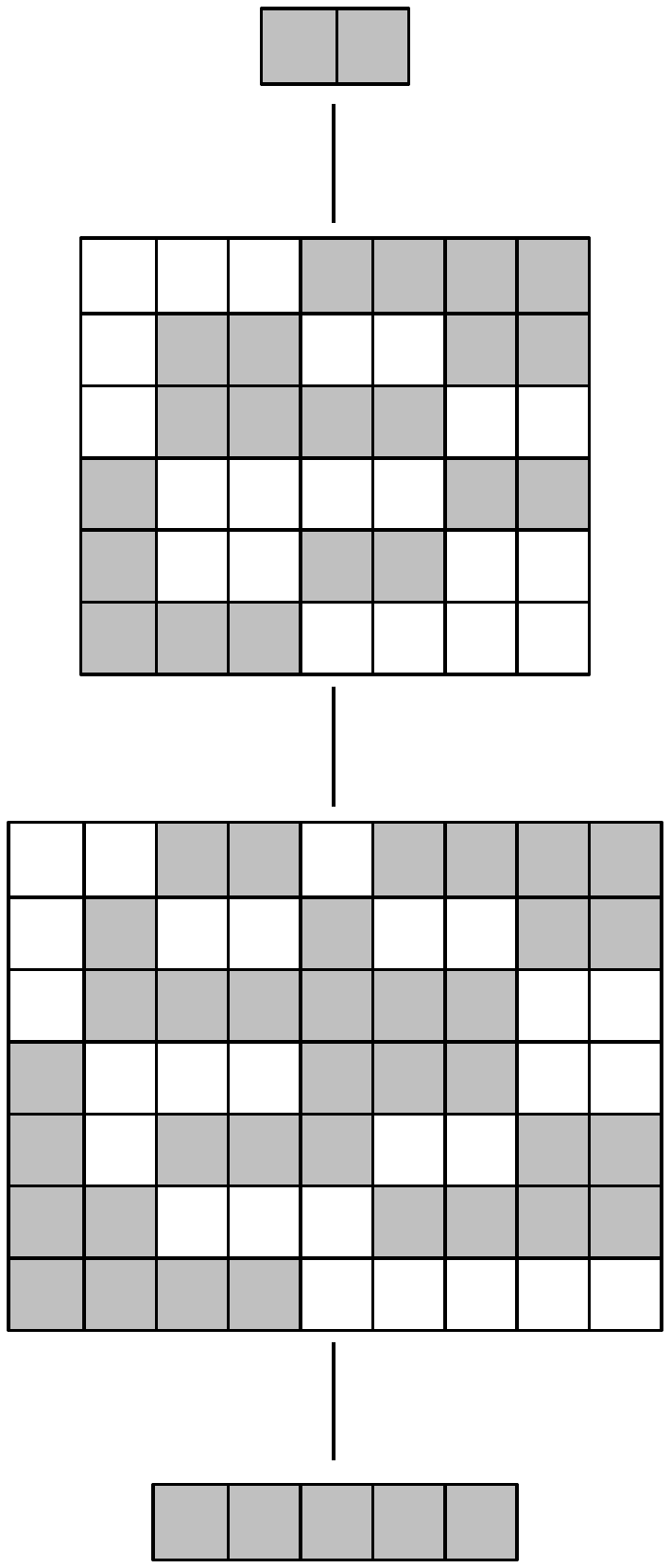} 
\caption[blah]{Egg-box diagrams of $\Reg(\TXal)$, where $|X|=\{1,2,3,4,5\}$ $X/\al=\{\{1\},\{2,3\},\{4,5\}\}$ (left), $X/\al=\{\{1\},\{2\},\{3,4,5\}\}$ (middle) or $X/\al=\{\{1\},\{2\},\{3\},\{4,5\}\}$ (right).}
\label{fig:RXal}
\end{center}
\end{figure}

Figures \ref{fig:TXA}--\ref{fig:RXal} give egg-box diagrams of special cases of the semigroups $\TXA$,~$\TXal$ and their regular subsemigroups; for comparison, Figure \ref{fig:TX} gives egg-box diagrams of $\T_X$ itself for small $|X|$.  These were produced with the aid of the Semigroups package for GAP \cite{GAP}, and may be used to visualise some of the results proved about these semigroups.    

For example, one may compare Figure \ref{fig:TXA} with Corollary \ref{cor:GreenSa}, which describes Green's relations in a principal left ideal (generated by a regular element).  One may also see the ``inflation'' discussed in Remark~\ref{rem:inflation_Sa} by comparing Figures \ref{fig:TX} and Figure \ref{fig:RXA}; each semigroup in Figure \ref{fig:RXA} is an ``inflation'' of a semigroup in Figure \ref{fig:TX}.  Figures \ref{fig:TXal} and \ref{fig:RXal} may be used to visualise the situation for principal right ideals.  The pdf may be zoomed significantly to see more detail in any figure, if required.

\section{Symmetric inverse monoids}\label{sect:I}

We conclude with a short section on symmetric inverse monoids.
Fix a non-empty set $X$, and denote by $\I_X$ the symmetric inverse monoid over $X$, as defined in Section \ref{subsect:trans}.  We also fix an element $a\in\I_X$ with the intention of studying the principal one-sided ideals $\I_Xa$ and $a\I_X$ of~$\I_X$.  Again, since $\I_X$ is regular (indeed, inverse), we may assume that $a$ is an idempotent: i.e., $a=\id_A$ for some $A\sub X$.  It is then easy to see that
\[
\I_Xa=\set{f\in\T_X}{\im(f)\sub A} \AND a\I_X=\set{f\in\I_X}{\dom(f)\sub A}.
\]
Clearly $f\mt f^{-1}$ determines an anti-isomorphism between $\I_Xa$ and $a\I_X$, so it suffices to consider just $\I_Xa$, as results for $a\I_X$ are dual.  In the literature, the semigroup $\I_Xa$ is generally denoted by $\IXA$, and we will continue to use this notation here.

Again, Green's relations and regular elements of $\IXA=\I_Xa$ are determined by the sets
\[
P=\set{f\in\IXA}{f\L af} \AND P'=\set{f\in\IXA}{f\J af}.
\]
Since $\I_X$ is inverse, every element of $\I_X$ (including $a$) is uniquely sandwich-regular, and so Theorem \ref{thm:inverse_P} gives 
\[
P=\Reg(\IXA)=a\I_Xa = \set{f\in\I_X}{\dom(f),\im(f)\sub A},
\]
and it is easy to see that this (local) monoid is isomorphic to $\I_A$; cf.~\cite[Theorem 3.1]{FS2014}.  Thus, any result concerning $\Reg(\IXA)$ reduces to a corresponding result concerning the well-studied inverse monoid $\I_A$.  As for the set $P'$, it is easy to see that for $f\in\IXA$, we have
\[
f\J af \iff \rank(f)=\rank(af) \iff \rank(f)=|A\cap\dom(f)|.
\]

\begin{thm}[cf.~Theorem \ref{thm:Green_TXA}]\label{thm:Green_IXA}
Let $X$ be a non-empty set, let $A$ a subset of $X$, and let $f,g\in\IXA$.  Then in the semigroup~$\IXA$,
\bit
\itemit{i} $f\L g \iff f=g$ or $[\im(f)=\im(g)$ and $\dom(f),\dom(g)\sub A]$,
\itemit{ii} $f\R g \iff \dom(f)=\dom(g)$,
\itemit{iii} $f\H g \iff f=g$ or $[\im(f)=\im(g)$ and $\dom(f)=\dom(g)\sub A]$,
\itemit{iv} $f\D g \iff \dom(f)=\dom(g)$ or $[\rank(f)=\rank(g)$ and $\dom(f),\dom(g)\sub A]$,
\itemit{v} $f\J g \iff \dom(f)=\dom(g)$ or $|A\cap\dom(f)|=\rank(f)=\rank(g)=|A\cap\dom(g)|$.
\eit
Further, ${\D}={\J}$ in $\IXA$ if and only if $A$ is finite or $A=X$.  \epfres
\end{thm}

\begin{rem}
Parts (i)--(v) were proved in \cite[Theorems 3.3, 3.4, 3.6 and 3.7]{FS2014}, but the final statement did not appear in \cite{FS2014}.
\end{rem}

Since $P=\Reg(\IXA)$ is inverse, the idempotent-generated subsemigroup $\bbE(\IXA)=\bbE(P)$ is simply the semilattice of idempotents $E(P)$, which is isomorphic to $E(\I_A)=\set{\id_B}{B\sub A}$; this, in turn, is isomorphic to the power set $\set{B}{B\sub A}$ under intersection.  Its rank and idempotent rank are equal to $1+|A|$ if $A$ is finite, or to $2^{|A|}$ if $A$ is infinite.

\footnotesize
\def\bibspacing{-1.1pt}
\bibliography{biblio}
\bibliographystyle{abbrv}
\end{document}